\documentclass[11pt]{pjm}  
\mathchardef\varSigma="0106  
\usepackage{amssymb,latexsym}

\font\smallbf=ptmbo at 10pt
\font\tensans=cmss10
\font\eightrm=cmr8

\def\ax{a}
\def\cj{{\hskip.4ptr}}
\def\xa{A}
\def\xb{B}
\def\ek{k}
\def\es{s}
\def\mk{m}
\def\ix{x}
\def\ez{z}
\def\bbP{{\mathbf{P}}}
\def\bbR{{\mathbf{R}}}
\def\rto{\text{\bf R}\hskip-.5pt^2}
\def\rtr{\text{\bf R}\hskip-.7pt^3}

\def\rmk{{\mathbf{R}}^{\hskip-.6ptm}}
\def\rn{{\mathbf{R}}^{\hskip-.6ptn}}
\def\bbC{{\mathbf{C}}}

\def\rp{{\mathbf{R}}\text{\rm P}}
\def\cp{{\mathbf{C}}\text{\rm P}}
\def\bbZ{{\mathbf{Z}}}
\def\dzerp{\mathfrak{D}^0_{\hskip-1pt+}}
\def\dzerm{\mathfrak{D}^0_{\hskip-1pt-}}
\def\dzerpm{\mathfrak{D}^0_{\hskip-1pt\pm}}
\def\donep{\mathfrak{D}^1_{\hskip-1pt+}}
\def\donem{\mathfrak{D}^1_{\hskip-1pt-}}
\def\donepm{\mathfrak{D}^1_{\hskip-1pt\pm}}

\def\depspm{\mathfrak{D}^\ve_{\hskip-1pt\pm}}
\def\CW{\text{\rm CW}}

\def\db{\hskip1.3pt\overline{\hskip-1.3pt\partial}}
\def\iff{if and only if}
\def\tr{totally real}
\def\tri{totally real immersion}
\def\tre{totally real embedding}
\def\ac{almost complex}
\def\nb{neighborhood}

\def\rmf{real manifold}
\def\inv{-in\-var\-i\-ant}
\def\feo{dif\-feo\-mor\-phism}
\def\feic{dif\-feo\-mor\-phic}
\def\feicly{dif\-feo\-mor\-phi\-cal\-ly}
\def\rsu{real surface}
\def\csu{complex surface}
\def\acm{almost complex manifold}
\def\acsu{almost complex surface}
\def\acst{almost complex structure}
\def\ie{im\-mer\-sion/em\-bed\-ding}
\def\ies{im\-mer\-sions/em\-bed\-dings}
\def\kb{Klein bottle}
\def\cml{total}
\def\vs{vector space}
\def\rvs{real vector space}
\def\cvs{complex vector space}
\def\vsu{vector subspace}
\def\sc{simply connected}

\def\nmf{$\,n$-man\-i\-fold}
\def\mfd{-man\-i\-fold}
\def\diml{-di\-men\-sion\-al}
\def\psh{pseu\-do\-hol\-o\-mor\-phic}
\def\y{{M}}
\def\tm{{T\hskip-.3pt\y}}
\def\x{\varSigma}

\def\fe{f}
\def\ha{I}
\def\s{S}
\def\py{\varPi}
\def\ly{\varLambda}
\def\gm{\gamma}
\def\vg{\varGamma}
\def\vh{{h}}
\def\vd{\Delta}
\def\sgm{\mu}
\def\mx{\mu}
\def\my{\sigma}
\def\fri{\mathfrak{I}}
\def\ind{\text{\smallbf i}\hskip1.2pt}
\def\jnd{\text{\smallbf j}\hskip1.2pt}
\def\dg{\text{\smallbf d}\hskip.6pt}
\def\hs{\hskip.7pt}
\def\hh{\hskip.4pt}
\def\nh{\hskip-.7pt}
\def\nnh{\hskip-1.5pt}
\def\h{\text{\rm Hom}\hs}

\def\p{P}
\def\q{Q}
\def\er{r}
\def\rx{r}
\def\n{N}
\def\w{W}
\def\k{K}
\def\op{o}
\def\rc{c\hh}
\def\ve{\varepsilon}
\def\spanc{\text{\rm Span}_{\hskip.5pt\bbC\hskip.3pt}}
\def\spanr{\text{\rm Span}_{\hskip.5pt\bbR\hskip.3pt}}
\def\dimr{\dim_{\hskip.4pt\bbR\hskip-1.2pt}}
\def\dimc{\dim_{\hskip.4pt\bbC\hskip-1.2pt}}
\def\detr{\text{\rm det}_\bbR}
\def\detc{\text{\rm det}_\bbC}
\def\eu{\text{\rm e}\hskip.6pt}

\def\df{d\hskip-.8ptf}
\def\pro{\text{\rm pr}\hs}
\def\pr{\hs'\hskip-1pt}
\def\simtri{\sim_{\hskip.4pt\text{\rm tri}}}
\def\emp{\text{\tensans\O}}
\def\compact{closed}
\def\ortwo{_{\text{\rm ord}\hs2}}

\begin{document} 

\renewcommand{\theequation}{\thesection.\arabic{equation}}

\title{Totally real immersions of surfaces}


\author[A. Derdzinski]{Andrzej Derdzinski} 
\address{Department of Mathematics\\
Ohio State University\\
Columbus, OH 43210, USA}
\email{andrzej@math.ohio-state.edu} 
\urladdr{http:/\hskip-1.5pt/www.math.ohio-state.edu/\~{}andrzej}
\author[T. Januszkiewicz]{Tadeusz Januszkiewicz}
\address{Department of Mathematics\\
Ohio State University\\
Columbus, OH 43210, USA\\
and\\
Mathematical Institute, Wroc\l aw University\\
and\\
Mathematical Institute, Polish Academy of Sciences} 

\thanks{The second author was partially supported by the NSF grant no.\ 
DMS-0405825}

\newtheorem{thm}{Theorem}[section] 
\newtheorem{prop}[thm]{Proposition} 
\newtheorem{lem}[thm]{Lemma} 
\newtheorem{cor}[thm]{Corollary} 

\theoremstyle{definition} 
  
\newtheorem{defn}[thm]{Definition} 
\newtheorem{notation}[thm]{Notation} 
\newtheorem{example}[thm]{Example} 
\newtheorem{conj}[thm]{Conjecture} 
\newtheorem{prob}[thm]{Problem} 
  
\theoremstyle{remark} 

\newtheorem{rem}[thm]{Remark} 


\begin{abstract}
Totally real immersions $\,f\,$ of a \compact\ \rsu\ $\,\x\,$ 
in an \acsu\ $\,\y\,$ are completely classified, up to homotopy through 
\tri s, by suitably defined homotopy classes $\,\mathfrak{M}(f)\,$ of mappings 
from $\,\x\,$ into a specific real $\,5$-man\-i\-fold $\,E(\y)$, while 
$\,\mathfrak{M}(f)\,$ themselves are subject to a single cohomology 
constraint. This follows from Gromov's observation that \tri s satisfy the 
$\,h$-prin\-ci\-ple. For the receiving \csu s $\,\bbC^2$, 
$\,\cp^1\nnh\times\cp^1$, $\,\cp^2$ and $\,\cp^2\#\,\mk\hs\overline{\cp^2}$, 
$\,\mk=1,2,\dots,7$, and all $\,\x\,$ (or, $\,\cp^2\#\,8\,\overline{\cp^2}$ 
and all orientable $\,\x$), we illustrate the above nonconstructive result 
with explicit examples of immersions realizing all possible equivalence 
classes. We also determine which equivalence  classes contain \tr\ embeddings, 
and provide examples of such embeddings for all classes that contain them.
\end{abstract} 
  
\subjclass{53C15, 53C42 (primary), 32Q60 (secondary).}
\maketitle 

\setcounter{section}{-1}
\section{Introduction}\label{in}
\setcounter{equation}{0}
Given an {\it \acsu}, that is, an \acm\ $\y$ with 
$\dimr\y\hskip-2pt=\hskip-1.2pt4$, we ask which \compact\ \rsu s $\x$ admit 
\tr\ \ies\ $f$ in $\y$, and how such $\nh f\nh$ can be classified up to the 
equivalence relation $\simtri\nh$ of being homotopic through \tri s.

We study these questions using a two\hs-\nh pronged approach. First, our 
Theorems~\ref{thzer} and \ref{thuno} provide an answer for \tr\ {\it 
immersions}. Theorem~\ref{thuno} states that, when $\,\y\,$ is \sc, the 
$\,\simtri$ equivalence class of a \tri\ $\,f:\x\to\y\,$ is completely 
determined by the {\it Maslov index} $\,\hs\ind\hs=\hs\ind(f)\,$ and {\it 
degree} $\,\hs\dg\hs=\hs\dg(f)\,$ of $\,f$, which in turn form an arbitrary 
element $\,(\ind\hs,\dg)\,$ of a specific set depending on $\,\y\,$ and 
$\,\x$. Theorem~\ref{thzer} classifies such $\,\simtri$ equivalence classes 
for arbitrary $\,\y$, using the {\it Maslov invariant\/} 
$\,\hs\mathfrak{M}(f)$, valued in a certain set of homotopy classes of 
mappings. Our definition of $\,\hs\mathfrak{M}(f)\,$ in Section~\ref{sr} is 
based on that in \cite{arnold}.

What makes Theorems~\ref{thzer} and \ref{thuno} less than completely 
satisfactory is the reliance of their proofs on Gromov's observation \cite[p.\ 
192]{gromov} that \tri s satisfy the $\,h$-prin\-ci\-ple. Consequently, those 
proofs offer little information about how the immersions which are shown to 
exist might actually be constructed. In addition, the two theorems deal only 
with the case of \tr\ {\it immersions}, as opposed to embeddings.

To make up for such shortcomings, we devote most of this paper (starting from 
Section~\ref{se}) to our second approach, which deals with \tri s and 
embeddings $\,\x\to\y\,$ of arbitrary \compact\ \rsu s $\,\x\,$ in the 
``model'' \sc\ \csu s
\begin{equation}\label{lis}
\bbC^2\,,\hskip5pt\cp^2\,,\hskip5pt
\cp^1\!\times\cp^1\hskip8pt\text{\rm and}\hskip10pt
\cp^2\#\,\mk\hs\overline{\cp^2}\nh,\hskip10pt\mk\,\ge\,1\hs,
\end{equation}
$\cp^2\#\,\mk\hs\overline{\cp^2}$ being obtained by blowing up any set of $\,\mk\,$ 
points in $\,\cp^2\nh$. For these $\,\y$, we provide explicit answers to the 
questions stated in the first paragraph. We begin by settling the existence 
question for \tr\ \ies\ in $\,\y$. (Table~\ref{summary} in Section~\ref{sr} 
summarizes the results.) Next, 
we construct examples of \tri s realizing all possible index-degree pairs 
$\,(\ind\hs,\dg)$. When $\,\mk\le7\,$ in (\ref{lis}), or $\,\mk=8\,$ and 
$\,\x\,$ is orientable, this also answers the question about embeddings: 
some of our examples are embeddings, while, as an intersection argument shows, 
none of the other pairs $\,(\ind\hs,\dg)\,$ comes from a \tre.

The constructions in question are summarized in Theorems~\ref{maire} 
and~\ref{maiem}.

In our second approach we make a point of not using Theorems~\ref{thzer} or 
\ref{thuno}. Instead, our nonexistence results are directly derived from 
obstructions involving characteristic classes and cohomology operations (for 
immersions), as well as intersection numbers (for embeddings); the existence 
assertions are in turn established by explicit constructions.

To prove Theorem~\ref{maiem} for $\,\cp^2\#\,\mk\hs\overline{\cp^2}$ with 
$\,\mk=6,7,8\,$ we study the case where a \tr\ \ie\ $\,\x\to\y\,$ of a \compact\ 
oriented \rsu\ $\,\x\,$ in an \acsu\ $\,\y\,$ can be obtained by deforming its 
exact opposite, that is, a {\it pseu\-do\-hol\-o\-mor\-phic\/} \ie\ 
$\,f:\x\to\y$. This turns out to be possible if and only if the normal bundle 
$\,\nu\,$ of $\,f\,$ and the tangent bundle $\,\tau\,$ of $\,\x\,$ are 
anti-isomorphic as complex line bundles (Theorem~\ref{thqnd}), while the 
required (small) deformation is achieved by moving in the direction of a 
section $\,\psi\,$ of $\,\nu\,$ such that $\,\db\hs\psi\ne0\,$ 
everywhere, $\,\db\,$ being a Cau\-chy-Rie\-mann operator. 
Although less constructive than the other existence proofs in the second part 
of the paper, such a technique is still more explicit than a general 
$\,h$-prin\-ci\-ple argument: it requires only finding a section $\,\psi\,$ of 
$\,\nu\,$ such that $\,\db\hs\psi\,$ trivializes the line 
bundle $\,\text{\rm Hom}_\bbC(\overline\tau,\nu)\,$ (assumed trivial to begin 
with). 

Our interest in the subject was sparked by the paper \cite{forstneric}.

\section{Preliminaries}\label{pr}
\setcounter{equation}{0}
By `planes' and `lines' we always mean {\em vector\/} spaces. All manifolds 
are of class $\,C^\infty$ and connected, except when explicitly stated 
otherwise.

An {\it \acm\/} is a \rmf\ $\,\y\,$ with an {\it \acst}, that is, a 
$\,C^\infty$ bundle morphism $\,J:\tm\to\tm\,$ such that 
$\,J^2=\,-\,\text{\rm Id}$. The tangent bundle $\,\tm\,$ then becomes a 
complex vector bundle, in which $\,J\,$ is the multiplication by $\,i$. We 
usually write $\,iv\,$ rather than $\,Jv$ for $\,v\in T_x\y\,$ and 
$\,x\in\y$.

A real vector subspace $\,W$ of a complex \vs\ $\,V\hs$ is said to be {\it 
totally real\/} if $\,\,W\cap\,iW=\{0\}$. A {\it totally real 
im\-mer\-sion/em\-bed\-ding\/} of a \rmf\ $\,\x\,$ (with or without 
boundary) in an \acm\ $\,\y\,$ is an \ie\ $\,f:\x\to\y\,$ such that the image 
of the differential $\,\df_x\,$ at any $\,x\in\x\,$ is a \tr\ subspace of 
$\,T_{f(x)}\y$. If $\,f\,$ is a \tre, the image $\,f(\x)\,$ is called a {\it 
totally real submanifold\/} of $\,\y$. (See \cite{chen-ogiue}, \cite{gromov} 
and Remark~\ref{deftr}.)

Given an \acm\ $\,\y\,$ with $\,\dimc\y=n$, we define $\,E^+(\y)\,$ and 
$\,E(\y)\,$ to be the unit circle bundles of the determinant bundle 
$\,\hs\detc\hs\tm\,=\,[\tm]^{\wedge n}$ and, respectively, of its square 
$\,[\detc\hs\tm]^{\otimes2}$. Thus, $\,E(\y)\,$ is the $\,\rp^1$ bundle over 
$\,\y\,$ associated with $\,\hs\detc\hs\tm$. Both $\,E=E(\y)\,$ and 
$\,E=E^+(\y)\,$ 
are the total spaces of principal $\,\hs\text{\rm U}\hs(1)$-bun\-dles over 
$\,\y$, leading to the homotopy exact sequences
\begin{equation}\label{seq}
\pi_2E\mathop{\,-\!\!\!-\!\!\!-\!\!\!-\!\!\!-\!\!\!\longrightarrow\,}
\limits^{\text{\rm injective}\,}
\pi_2\y\mathop{\,-\!\!\!-\!\!\!-\!\!\!-\!\!\!-\!\!\!-\!\!\!-\!\!\!\!
\longrightarrow\,}\limits^{\text{\rm connecting}\,}
\pi_1[\hs\text{\rm U}\hs(1)]=\bbZ
\,\longrightarrow\,\pi_1E
\mathop{\,-\!\!\!-\!\!\!\!\!\longrightarrow\,}\limits^{\text{\rm onto}\,}
\pi_1\y\hs.
\end{equation}
One also has an obvious twofold covering projection 
\begin{equation}\label{tcp}
E^+(\y)\,\to\,E(\y)\,=\,E^+(\y)/\bbZ_2\,,
\end{equation}
equivariant relative to the homomorphism 
$\,\text{\rm U}\hskip.6pt(1)\ni z\mapsto z^2\in\text{\rm U}\hskip.6pt(1)$. 
Thus, 
\begin{equation}
\label{pie}
\pi_1[E^+(\y)]\,\subset\,\pi_1[E(\y)]\,,\hskip15pt
\pi_1[E(\y)]\,/\,\pi_1[E^+(\y)]\,=\,\bbZ_2\,,
\end{equation}
If, in addition, $\,\y\,$ is \sc, (\ref{pie}) and exactness of (\ref{seq}) 
give
\begin{equation}\label{pzq}
\pi_1[E(\y)]\nh=\bbZ_q
\hskip5pt\text{\rm for\ some}\hskip5ptq\in\{2,4,\dots,\infty\}\hs,\hskip5pt
\text{\rm where\ we\ set}\hskip5pt\bbZ_\infty\hskip-1pt=\hs\bbZ\hs.
\end{equation}
(More on $\,q\,$ can be found in Section~\ref{tq}.) For any manifold $\,\x\,$ 
and any Abelian group $\,G$, we have natural isomorphic identifications
\begin{equation}\label{coh}
H^1(\x,G)\,=\,\,\h(\pi_1\x,G)\,
=\,\h(H_1(\x,\bbZ),G)\hs.
\end{equation}
Thus, according to (\ref{coh}) with $\,G=\bbZ_2$,
\begin{equation}\label{ori}
w_1(\x)\in H^1(\x,\bbZ_2)\hskip5pt\text{\rm is\ the\ orientation\ 
homomorphism}\hskip6pt\pi_1\x\to\bbZ_2\hs.
\end{equation}
We will also use Wu's formula \cite{milnor-stasheff}, valid whenever $\,\x\,$ 
is a \compact\ \rsu:
\begin{equation}\label{wus}
w_1(\x)\nh\smallsmile w_1(\x)\,=\,[\chi(\x)\,\text{\rm mod}\hskip4pt2\hs]
\,\in\,H^2(\x,\bbZ_2)\,=\,\bbZ_2\hs.
\end{equation}
\begin{rem}\label{reduc}In terms of (\ref{coh}), the homomorphism 
$\,H^1(\x,G)\to H^1(\x,G')\,$ of coefficient reduction, corresponding to a 
homomorphism $\,h:G\to G'$ of Abel\-ian groups, sends a homomorphism 
$\,\varphi:H_1(\x,\bbZ)\to G\,$ to $\,h\circ\varphi$.
\end{rem}
\begin{rem}\label{wones}The following properties of \compact\ \rsu s $\,\x\,$ 
are both well known and easily derived from Remark~\ref{srfcs} and 
(\ref{hss}.i) in Section~\ref{ch}: the torsion subgroup of $\,H_1(\x,\bbZ)\,$ is 
trivial when $\,\x\,$ is orientable, and isomorphic to $\,\bbZ_2$ otherwise; 
while, if $\,\chi(\x)\,$ is odd, the torsion subgroup is {\it not\/} contained 
in the kernel of $\,w_1(\x)$.
\end{rem}
\begin{rem}\label{deftr}Our definition of a \tr\ subspace $\,W\subset V\hs$ 
differs from that in \cite{gromov}, where $\,\hs\spanc W\,$ is required to 
have the maximum possible dimension $\,\hs\text{\rm min}\hs(\ek,n)$, for 
$\,\ek=\dimr W\,$ and $\,n=\dimc V\nh$. The two definitions agree if 
$\,\ek\le n\hs$, and are both devoid of content if $\,\ek\ge2n-1$ (as one of 
them then makes every subspace \tr, and the other allows no such subspace 
unless $\,n=\ek=1$); however, $\,\ek\le n\,$ or $\,\ek\ge2n-1$ when $\,n=2\,$ 
(and $\,\ek\le2n$), which is the case of our main interest.
\end{rem}

\section{Statement of the results}\label{sr}
\setcounter{equation}{0}
Let $\,E(\y)\,$ and $\,E^+(\y)\,$ be as in Section~\ref{pr} for an \acm\ $\,\y$ 
with $\,\dimc\y=n$. In view of (\ref{pie}), there exists a unique homomorphism 
\begin{equation}\label{wun}
\mathfrak{w}_1:\pi_1[E(\y)]\to\bbZ_2\qquad\text{\rm with}\qquad
\text{\rm Ker}\hskip2.6pt\mathfrak{w}_1\,=\,\pi_1[E^+(\y)]\,.
\end{equation}
Thus (cf.\ (\ref{coh})), 
$\,\hs\mathfrak{w}_1\in H^1(E(\y),\bbZ_2)\smallsetminus\{0\}\,$ is the first 
Stief\-el\hs-Whit\-ney class of the real line bundle over $\,E(\y)\,$ 
associated with the $\,\bbZ_2$ bundle (\ref{tcp}).

If $\,f:\x\to\y\,$ now is a \tri\ of a \rmf\ $\,\x\,$ in an \acm\ $\,\y\,$ 
with $\,\dimr\x=\dimc\y=n$, we define the {\it Maslov invariant\/} 
$\,\hs\mathfrak{M}(f)\,$ of $\,f\,$ to be the homotopy class of the mapping 
$\,\x\to E(\y)\,$ that sends $\,x\in\x\,$ to the real line in 
$\,[T_{f(x)}\y]^{\wedge n}$ spanned by the (complex) inner product 
$\,v_1\wedge\hs\ldots\hs\wedge v_n$, where $\,v_j=\df_xe_j$ for any basis 
$\,e_1,\dots,e_n$ of $\,T_x\x$. (See also Section~\ref{zt} below and 
\cite{arnold}.)

Obviously, $\,\hs\mathfrak{M}(f)\,$ depends only on the $\,\simtri$ 
equivalence class of $\,f$, for $\,\simtri$ defined in the Introduction. In 
Section~\ref{pz} we prove the following result.
\begin{thm}\label{thzer}Given an \acsu\ $\,\y\,$ and a 
\compact\ \rsu\ $\,\x$, the assignment\/ 
$\,f\mapsto\hs\mathfrak{M}(f)\,$ establishes a one-to-one correspondence 
between the set of all\/ $\,\simtri$ equivalence classes of \tri s 
$\,f:\x\to\y\,$ and the set of those homotopy classes 
of mappings $\,\varTheta:\x\to E(\y)\,$ for which
\begin{equation}\label{wow}
\varTheta^*\mathfrak{w}_1\nnh=w_1(\x)\hskip5pt\text{\rm in}\hskip4ptH^1(\x,\bbZ_2)\hs,
\hskip3pt\text{\rm with}\hskip6pt\mathfrak{w}_1\nnh\in\nh H^1(E(\y),\bbZ_2)
\text{\rm\ as\ in\ (\ref{wun}).}
\end{equation}
\end{thm}
If $\,\y\,$ in Theorem~\ref{thzer} happens to be \sc, $\,\hs\mathfrak{M}(f)\,$ 
may be replaced by a pair of more tangible invariants: the {\it Maslov 
index\/} $\,\hs\ind(f)\,$ of a \tri\ $\,f:\x\to\y$, and its {\it degree\/} 
$\,\hs\dg(f)$, described below. 

Specifically, discussing mappings $\,f\,$ from a real \nmf\ $\,\x\,$ into an 
\acm\ $\,\y$, we will assume that an orientation of $\,\x\,$ has been 
selected, as long as one exists; that is, $\,\x\,$ {\it is either oriented, or 
nonorientable.} 
We then define a group $\,\bbZ_{[2]}$ associated with $\,\x\,$ by
\begin{equation}\label{ztw}
\text{\rm$\bbZ_{[2]}=\,\bbZ\,\,$ if $\,\,\x\,\,$ is oriented, 
$\,\,\,\bbZ_{[2]}=\,\bbZ_2\,$ if $\,\,\x\,\,$ is not orientable.}
\end{equation}
Let $\,[\x\hs]\in H_n(\x,\bbZ_{[2]})\,$ be the fundamental homology class of 
$\,\x$. We set
\begin{equation}\label{dgf}
\dg(f)=f_*[\x\hs]\in H_n(\y,\bbZ_{[2]})\hskip5pt\text{\rm with}\hskip6ptn
=\dimr\x\hskip5pt\text{\rm and}\hskip6pt\bbZ_{[2]}\hskip5pt\text{\rm as\ in\ 
(\ref{ztw}).}
\end{equation}
If $\,f\,$ is a \tri\ of a \rmf\ $\,\x\,$ in an \acm\ $\,\y\,$ and 
$\,\dimr\x=\dimc\y=n$, we define $\,\hs\ind(f)\,$ to be the homomorphism 
$\,\pi_1\x\to\pi_1[E(\y)]\,$ induced by $\,\hs\mathfrak{M}(f)$. Thus (cf.\ 
(\ref{pzq}), (\ref{coh})),
\begin{equation}\label{ifh}
\ind(f)\,\in\,H^1(\x,\,\bbZ_q)\,.
\end{equation}

Rather than being arbitrary elements of the 
(co)homology groups in question, 
$\,\hs\ind(f)\hs\,$ and $\,\hs\dg(f)\hs\,$ are both confined to specific 
subsets. Namely (see Lemma~\ref{ifidf}(b)), for a \tri\ $\,f\,$ of a \compact\ 
\rsu\ $\,\x\,$ in a \sc\ \acsu\ $\,\y$, and $\,q\,$ as in (\ref{pzq}),
\begin{equation}\label{ifi}
\ind(f)\in\fri_q(\x)\,\subset\,H^1(\x,\bbZ_q)\hskip8pt\text{\rm and}
\hskip9pt\dg(f)\hs\in\,\depspm(\y)\,\subset\,H_2(\y,\bbZ_{[2]})\hs.
\end{equation}
Here $\,\ve\,$ and $\,\pm\,$ are 
$\,\bbZ_2$-val\-ued parameters, determined by $\,\x\,$ as follows:
\begin{equation}\label{dfd}
\ve=1\,\text{\rm\ if }\,\x\,\text{\rm\ is orientable,}\ \ \ve=0\,
\text{\rm\ if\ it\ is\ not, and}\hskip5pt(-1)^{\chi(\x)}=\pm1\hs,
\end{equation}
while $\,\hs\dzerpm(\y)\subset H_2(\y,\bbZ_2)\,$ and 
$\,\hs\donepm(\y)\subset H_2(\y,\bbZ)\,$ are defined by
\begin{equation}\label{dmk}
\alignedat2
\dzerp(\y)\,&=\,\,\text{\rm Ker}\,[w_2(\y)]\,,
\hskip16pt&&\dzerm(\y)
=H_2(\y,\bbZ_2)\smallsetminus\hs\text{\rm Ker}\,[w_2(\y)]\hs,\\
\donep(\y)&=\,\text{\rm Ker}\,[c_{\hs1}(\y)]\hs,\hskip16pt
&&\donem(\y)=\,\emp\hs,\endalignedat
\end{equation}
with $\,c_{\hs1}(\y):H_2(\y,\bbZ)\to\bbZ\,$ and 
$\,w_2(\y):H_2(\y,\bbZ_2)\to\bbZ_2$. Therefore, $\,\hs\depspm(\y)$, if 
nonempty, is a coset of a subgroup in $\,H_2(\y,\bbZ_2)\,$ or 
$\,H_2(\y,\bbZ)$. Finally,
$\,\hs\fri_q(\x)\subset H^1(\x,\bbZ_q)\,$ is given by
\begin{equation}\label{iqs}
\fri_q(\x)\,=\,\,\{\lambda\in H^1(\x,\bbZ_q):
[\lambda\hskip2.5pt\text{\rm mod}\hskip4pt2\hskip1pt]\,=\,w_1(\x)\}\,,
\end{equation}
where $\,H^1(\x,\bbZ_q)\ni\lambda\mapsto[\lambda\hskip2.5pt\text{\rm mod}
\hskip4pt2\hskip1pt]\in H^1(\x,\bbZ_2)\,$ denotes the 
$\,\hs\text{\rm mod}\hskip4pt2\hs$ 
reduction homomorphism corresponding to the unique nonzero homomorphism 
$\,\bbZ_q\to\bbZ_2$. (Recall that $\,q\in\{2,4,6,\dots,\infty\}\,$ in 
(\ref{pzq}).) In Sections \ref{sj} and \ref{ij} we establish the following 
theorem.
\begin{thm}\label{thuno}Given a \sc\ \acsu\ $\,\y\,$ and a 
\compact\ \rsu\/ $\,\x$, the assignment\/ 
$\,f\mapsto(\ind(f)\hs,\,\dg(f))\,$ defines a bijective correspondence 
between the set of $\hs\simtri\nh$ equivalence classes of totally real 
immersions $\hs f:\x\to\y\,$ and a specific subset\/ $\,\mathcal{Z}(\x,\y)\,$ 
of the Cartesian product\/ $\,\hs\fri_q(\x)\times\depspm(\y)$, in the notation 
of\/ {\rm(\ref{pzq})} and\/ {\rm(\ref{dfd})} -- {\rm(\ref{iqs})}.
\end{thm}
The set $\,\mathcal{Z}(\x,\y)\,$ coincides with 
$\,\hs\fri_q(\x)\times\depspm(\y)\,$ except in 
the case where, simultaneously, $\,\x\,$ is nonorientable and $\,\chi(\x)\,$ 
is even, while $\,q\,$ is finite and divisible by $\,4$. In this latter case, 
$\,\mathcal{Z}(\x,\y)\,$ has half of the (finite) number of elements of 
$\,\hs\fri_q(\x)\times\depspm(\y)$, and consists of those $\,(\ind\hs,\dg)$ in 
$\,\hs\fri_q(\x)\times\depspm(\y)\,$ which satisfy the following condition: if 
the image of the unique torsion element in $\,H_1(\x,\bbZ)\,$ under 
$\,\hs\ind:H_1(\x,\bbZ)\to\bbZ_q$ is zero (or, nonzero), then 
$\,\hs\dg\in H_2(\y,\bbZ_2)\,$ is (or, respectively, is not) the 
$\,\hs\text{\rm mod}\hskip4pt2$ reduction of an element of 
$\,\hs\text{\rm Ker}\,[c_{\hs1}(\y)]$. (Cf.\ (\ref{coh}) and 
Remark~\ref{wones}.)

In contrast with Theorems~\ref{thzer} and \ref{thuno}, which make use of the 
$\,h$-prin\-ci\-ple, all the results listed below are established by explicit 
geometric arguments. First, in Section~\ref{pf} we prove the following six 
statements.
\begin{thm}\label{thdue}Let\/ $\,\y\,$ be any \acsu. The class of 
\compact\ \rsu s $\,\x\,$ admitting a totally real 
embedding {\rm(\hskip-1.1pt}or, immersion{\rm)} in\/ $\,\y$ then includes 
the\/ $\,2$-torus $\,T^2$ and \kb\ $\,K^2$ {\rm(\hskip-.9pt}and, 
for immersions, also the\/ $\,2$-sphere $\,S^2${\rm)}, and is closed under the 
mapping $\,\x\,\mapsto\,\x\,\#\,T^2\#\,K^2$ {\rm(\hskip-.9pt}and, for 
immersions, under the con\-nect\-ed-sum operation\/ 
$\,(\x,\x\hs')\,\mapsto\,\x\,\#\,\x\hs'${\rm)}. 
\end{thm}
\begin{cor}\label{cotre}Any \compact\ \rsu\ with an even Euler characteristic 
admits a \tri\ in every \acsu.
\end{cor}
\begin{cor}\label{coqtr}The $\,2$-torus $\,T^2$ and all nonorientable 
\compact\ \rsu s $\,\x\,$ with 
$\,\chi(\x)\equiv0\hskip6pt\text{\rm mod}\hskip4pt4\,$ admit totally real 
embeddings in every \acsu. 
\end{cor}
\begin{cor}\label{cocin}Let\/ $\,\y\,$ be an \acsu. If 
there exists a totally real immersion\/ $\,\rp^2\to\y$, then 
every \compact\ \rsu\ $\,\x\,$ admits a totally real immersion in\/ 
$\,\y$. 
\end{cor}
\begin{cor}\label{cosei}Every \compact\ \rsu\ admits totally real 
immersions in\/ $\,\cp^2$ and in all \csu s obtained from\/ $\,\cp^2$ by 
blowing up any finite number of points.
\end{cor}
\begin{cor}\label{coset}The torus $\,T^2$, sphere $\,S^2$, and all 
nonorientable \compact\ surfaces admit totally real embeddings in 
$\,\cp^2\#\,\mk\hs\overline{\cp^2}$ for every integer $\,\mk\ge2$. 
\end{cor} 
Here and in the sequel, given a \csu\ $\,\y$, the connected sum 
$\,\y\,\#\,\,\es\,\overline{\cp^2}$ stands for {\it any\/} \csu\ obtained by 
blowing up $\,\es\,$ distinct points in $\,\y$.

Corollaries~\ref{cotre}, \ref{coset} and Theorem~\ref{thdue} 
lead in turn to the following three results, the detailed proofs of which are 
given in Section~\ref{pf} and Section~\ref{os}:
\begin{cor}\label{conov}For any \acsu\ $\,\y\,$ which is a 
spin manifold, the \compact\ \rsu s $\,\x\,$ that admit totally real 
immersions in\/ $\,\y$ are precisely those having even Euler 
characteristics. This is, for instance, the case for\/ $\,\y\nh=\bbC^2$ and\/ 
$\,\y\hs=\,\cp^1\!\times\cp^1\nnh$.
\end{cor}
\begin{cor}\label{codie}The class of \compact\ \rsu s admitting 
totally real embeddings in\/ $\,\cp^1\!\times\cp^1$
consists of the torus $\,T^2\nnh$, the sphere $\,S^2\nnh$, and all 
nonorientable \compact\ surfaces with even Euler characteristics. 
\end{cor}
\begin{cor}\label{cound}For any fixed\/ $\,\mk\in\{2,3,\dots,9\}$, the 
class of \compact\ \rsu s admitting a totally real embedding in\/ 
$\,\cp^2\#\,\mk\hs\overline{\cp^2}$ consists of the torus\/ $\,T^2\nnh$, the 
sphere\/ $\,S^2\nnh$, and all nonorientable \compact\ surfaces.
\end{cor} 
In each of the last three corollaries an existence statement based on an 
explicit elementary construction is coupled with a nonexistence assertion that 
uses elementary topological obstructions: an intersection-number relation 
(\ref{dot}) or (\ref{dtf}), and a condition involving either the first Chern 
class (for orientable surfaces $\,\x$), or Stief\-el\hs-Whit\-ney classes 
(for nonorientable $\,\x$).

The Stief\-el\hs-Whit\-ney class obstruction in (\ref{wco}.a) fails, however, 
to detect that some nonorientable \compact\ surfaces $\,\x\,$ do not admit 
\tre s in $\,\bbC^2$, $\,\cp^2$ or $\,\cp^2\#\,\overline{\cp^2}\nh$. Instead, 
we use Massey's formula (\ref{msf}) involving 
$\,\hs\text{\rm mod}\hskip4pt4\,$ intersection numbers and Pontryagin squares. 
This leads to a proof, in Section~\ref{ns}, of the next three corollaries. 
\begin{cor}\label{codod}A \compact\ \rsu\/ $\,\x\,$ admits a 
totally real embedding in\/ $\,\bbC^2$ \iff\/ $\,\x\,$ is either 
\feic\ to the torus $\,T^2$, or\/ $\,\x\,$ is nonorientable and\/ 
$\,\chi(\x)\equiv0\hskip6pt\text{\rm mod}\hskip4pt4$. 
\end{cor}
\begin{cor}\label{cotrd}The class of \compact\ \rsu s that admit a 
totally real embedding in\/ $\,\cp^2$ consists of the torus $\,T^2$ and all 
nonorientable \compact\ surfaces $\,\x\,$ with $\,\chi(\x)\equiv0\,$ or 
$\,\chi(\x)\equiv1\hskip6pt\text{\rm mod}\hskip4pt4$. 
\end{cor}
\begin{cor}\label{coqtd}The \compact\ \rsu s admitting totally real embeddings 
in the \csu\/ $\,\,\cp^2\#\,\overline{\cp^2}$ are\/{\rm:} the torus 
$\,T^2\nnh$, and all nonorientable \compact\ surfaces whose Euler 
characteristics are odd or divisible by\/ $\,4$.
\end{cor}
Corollaries~\ref{cosei} and~\ref{conov} --~\ref{coqtd} are summarized in 
Table~\ref{summary}. For a conclusion similar to but weaker than 
Corollaries~\ref{cotrd} and~\ref{coqtd}, see Proposition~\ref{trchi}.
\begin{table}[ht]
\vskip-6pt
\caption{Totally-real immersibility/embeddability of \compact\ \rsu s 
$\,\x\,$ in the \csu s (\ref{lis}) with $\,\mk\le9$. Here t.r.\ means 
`totally real' and $\,\hs\chi_4\in\{0,1,2,3\}\hs\,$ stands for 
$\,\hs\chi(\x)\,\,\text{\rm mod}\hskip2pt4$.}
\label{summary}
\vskip-6pt
\hskip6pt
\vbox{\offinterlineskip
\hrule
\halign
{&\vrule#&\strut
\hskip1.4pt\hfil#\hfil\hskip1.4pt\cr
&\hfil\vbox{\vskip-.8pt\hbox{$_{\text{\rm the \csu}}$}\vskip-2pt}\hfil&
&\vbox{\vskip-1pt\hbox{$_{\y=\,\cp^2}$}\vskip-1.1pt}&
&\vbox{\vskip-.8pt\hbox{$_{\y=\,\cp^2\#\,\mk\hs\overline{\cp^2}}$}\vskip-2.3pt}&
&\vbox{\vskip-.8pt\hbox{$_{\y=\,\cp^1\!\times\cp^1}$}\vskip-1.7pt}&
&\vbox{\vskip-1pt\hbox{$_{\y=\,\bbC^2}$}\vskip-1.1pt}&\cr
\noalign{\hrule}
&\hfil$_{\text{\rm which\ orientable\ }\x\text{\rm\ are}}$\hfil&
&\vbox{\vskip7.2pt\hbox{$_{T^2\hskip4pt\text{\rm only}}$}\vskip-7.2pt}&
&\hfil$_{T^2\text{\rm\ \ only\ if\ }\,\mk\hs=\hs1\hs;}$\hfil&&
\vbox{\vskip7.2pt\hbox{$_{S^2\nh,\hs\,T^2\text{\rm\ only}}$}\vskip-7.2pt}&
&\vbox{\vskip7.2pt\hbox{$_{T^2\text{\rm\ only}}$}\vskip-7.2pt}&\cr
&\hfil$^{\text{\rm t.r.\ embeddable\ in\ }\y}$\hfil&&&
&$^{S^2\nh,\hs\,T^2\text{\rm\ if\ }\hs2\,\le\,\mk\,\le\,9}$&&&&&\cr
\noalign{\hrule}
&\hfil$_{\text{\rm t.r.\ embeddabil-}}$\hfil&&&&\hfil$_{\chi_4\,\ne\,\hs2\,
\text{\rm\ if\ }\,\mk\,=\,1\hs;}$\hfil&&&
&$_{\chi(\x)}$&\cr
&\hfil\vbox{\vskip0pt
\hbox{\eightrm ity\ condition\ for}
\vskip-1pt}\hfil&&\hfil$_{\hskip3pt\chi_4\hs\in\hs\{0,1\}}$
\hfil&&\hfil\vbox{\vskip-1pt
\hbox{$_{\text{\rm all\ }\nh\x\hs\text{\rm\ t.r.\ embed-}}$}\vskip1pt}\hfil&
&\hfil$_{\chi(\x)\,\text{\rm\ even}}$\hfil&
&\hfil{\eightrm divisible}\hfil&\cr
&\hfil$^{\text{\rm nonorientable\ }\x}$\hfil&&&&$^{\text{\rm dable\ if\ }
\hs2\,\le\,\mk\,\le\,9}$&&&
&\hfil$^{\text{\rm by\ }4}$\hfil&\cr
\noalign{\hrule}
&\hfil$_{\text{\rm when\ a\hskip4.5ptt.r.\hskip3ptimmer-}}$\hfil&
&\multispan3\hfil$_{\text{\rm exists\ for\ all\ }\,\x}$\hfil&
&\multispan3\hfil$_{\text{\rm\ \iff\ }\chi(\x)\text{\rm\ is\ 
even\ }}$\hfil&\cr
&\hfil$^{\text{\rm sion\ }\,\x\,\to\,\y\text{\rm\ exists}}$\hfil&
&\multispan3\hfil$^{\text{\rm(since\ it\ does\ for\ }\,\x\,=\,\rp^2\text{\rm)}}
$\hfil&
&\multispan3\hfil\ $^{\text{\rm(just\ because\ }\,\y\,\text{\rm\ is\ spin)}}$
\hfil&\cr
\noalign{\hrule}
}}   
\vskip-5pt
\end{table}

Further such results can be derived from the following more general existence 
theorem, proved in Section~\ref{dc}\ via another explicit argument.
\begin{thm}\label{thott}Let\/ $\,\y'$ be the \csu\ obtained by blowing up 
$\,\ek$ distinct points, $\,\ek\ge1$, in a given \csu\ $\,\y$. The class of 
\compact\ real surfaces admitting a totally real embedding in\/ $\,\y'$ then 
includes
\begin{enumerate}
  \def\theenumi{{\rm\alph{enumi}}}
\item the\/ $\,2$-sphere\/ $\,S^2$, if\/ $\,\ek\ge2$,
\item the connected sum\/ $\,\x\,\#\,s\hskip2pt\overline{\rp^2}$, whenever 
$\,s\in\{0,1,\dots,\ek\}\,$ and\/ $\,\x\,$ is any totally real \compact\ 
surface embedded in\/ $\,\y\,$ that contains at least\/ $\,s\,$ of the 
$\,\ek\,$ blown-up points.
\end{enumerate}
\end{thm}
\begin{thm}\label{maire}If\/ $\,\y\,$ is one of the \csu s {\rm(\ref{lis})} 
and\/ $\,\x\,$ is a \compact\ \rsu, then\/ $\,\mathcal{Z}(\x,\y)\,$ defined in 
the lines following Theorem~\ref{thuno} coincides with the set of all pairs\/ 
$\,(\ind,\dg)\,$ such that\/ $\,\hs\ind(f)=\,\ind\hs\,$ and\/ 
$\,\hs\dg(f)=\,\dg$ for some \tri\ $\,f:\x\to\y$.
\end{thm}
The last result, as stated, is a special case of Theorem~\ref{thuno}; however, 
we establish it separately in Section~\ref{su}, using -- in contrast with our 
proof of Theorem~\ref{thuno} in Sections~\ref{sj} and~\ref{ij} -- only 
explicit geometric constructions.

Geometric constructions are also used to prove, in Section~\ref{rc}, a similar 
theorem about \tre s. First, according to (\ref{dtf}) and (\ref{psq}),
\begin{equation}\label{ddh}
\alignedat2
\text{\rm i)}\quad&\dg\,\cdot\,\dg\,\,=\,-\hs\chi(\x)\,\in\,\bbZ\,,\quad&&
\text{\rm if\ }\,\x\,\text{\rm\ is\ orientable,}\\
\text{\rm ii)}\quad&\dg^{\hskip.2pt2}\,=\,\,[\hs\chi(\x)\hskip2pt\text{\rm mod}
\hskip4pt4\hs]\,\in\,\bbZ_4\,,\quad&&\text{\rm if\ }\,\x\,\text{\rm\ is\ not\ 
orientable.}
\endalignedat
\end{equation}
Here $\,\hs\dg\,=\,\dg(f)\,$ is the degree of any \tre\ of a \compact\ \rsu\ 
$\,\x\,$ in an \acsu\ $\,\y\,$ and $\,\cdot\,$ denotes the intersection form 
in $\,H_2(\y,\bbZ)$. In (\ref{ddh}.ii), $\,\y\,$ is assumed, in addition, to 
be either \feic\ to $\,\bbC^2$ (and then $\,\hs\dg^{\hskip.2pt2}$ 
stands for $\,0\in\bbZ_4$), or compact and \sc\ (and then we set 
$\,\hs\dg^{\hskip.2pt2}=
[\hs\xi\smallsmile\xi\hskip6pt\text{\rm mod}\hskip4pt4\hs]
\in H^4(\y,\bbZ_4)=\bbZ_4$ for any $\,\xi\in H^2(\y,\bbZ)\,$ such that the 
$\,\hs\text{\rm mod}\hskip4pt2\,$ reduction of $\,\xi\,$ is the class 
$\,\sgm\in H^2(\y,\bbZ_2)$, Poin\-ca\-r\'e-dual to 
$\,\hs\dg\in H_2(\y,\bbZ_2)$).

Let $\,\y\,$ be one of the \csu s (\ref{lis}). If $\,\x\,$ is any \compact\ 
\rsu, the index-degree pairs $\,(\ind\hs,\dg)\,$ of \tri s $\,\x\to\y\,$ are 
precisely the elements of $\,\mathcal{Z}(\x,\y)$. (See Theorem~\ref{maire} or 
Theorem~\ref{thuno}.) That is no more the case when `immersion' is replaced by 
`embedding' since, as we just saw, (\ref{ddh}) then must hold as well. For the 
surfaces (\ref{lis}) with $\,1\le\mk\le7$, condition (\ref{ddh}) is, in fact, 
the only additional obstruction:
\begin{thm}\label{maiem}Let $\hs\y\nnh$ be one of the 
\csu s \hbox{{\rm(\ref{lis})}\hskip3pt{\it 
for}\hskip3.5pt$1\nnh\le\nh\mk\nh\le\nh8$,} 
and let\/ $\,\x\,$ be a \compact\ \rsu, orientable if\/ $\,\mk=8$. Any 
$\,(\ind,\dg)\,$ in the set\/ $\,\mathcal{Z}(\x,\y)\,$ described immediately 
after Theorem~\ref{thuno}, such that\/ $\,\dg\,$ satisfies\/ {\rm(\ref{ddh})}, 
then equals $\,(\ind(f),\dg(f))\,$ for some \tre\ $\,f:\x\to\y$. 
\end{thm}
Our proof of Theorem~\ref{maiem} for $\,\cp^2\#\,\mk\hs\overline{\cp^2}$ with 
$\,\mk=6,7,8\,$ uses a fact, stated below as Corollary~\ref{codot}, and 
obtained by deforming pseu\-do\-hol\-o\-mor\-phic \ies\ to \tr\ ones (cf.\ the 
end of the Introduction). The following four results are proved in Sections 
\ref{dp} -- \ref{po}.
\begin{thm}\label{thqnd}Given a pseu\-do\-hol\-o\-mor\-phic 
im\-mer\-sion/em\-bed\-ding\/ $\,\,f\,$ of a \compact\ oriented \rsu\/ 
$\,\x\,$ in an \acsu\/ $\,\y$, the following three conditions 
are equivalent\/{\rm:}
\begin{enumerate}
  \def\theenumi{{\rm\alph{enumi}}}
\item the tangent bundle\/ $\,\tau\,$ of\/ $\,\,\x$ and the normal bundle\/ 
$\,\nu\,$ of\/ $\,\,f\,$ are anti-isomorphic as complex line bundles,
\item $\,f^*[\detc\hs\tm]\,$ is trivial, that is, $\,f^*c_{\hs1}(\y)\,=\,0\,$ 
in\/ $\,H^2(\x,\bbZ)$,
\item $\,f\,$ is homotopic through \ies\ $\,\x\to\y\,$ to 
a totally real im\-mer\-sion/em\-bed\-ding $\,f':\x\to\y$. 
\end{enumerate}
Moreover, $\,f'$ in {\rm(c)} then can be chosen arbitrarily $\,C^1$-close to 
$\,f$.
\end{thm}
Here we call an immersion $\,f:\x\to\y\,$ {\it \psh} \cite{gromov} if, for 
each $\,x\in\x$, the image $\,\df_x(T_x\x)\,$ is a complex line in 
$\,T_{f(x)}\y\,$ and the isomorphism $\,\df_x:T_x\x\to\df_x(T_x\x)\,$ is 
o\-ri\-en\-ta\-tion-pre\-serv\-ing.

Let $\,\y'$ be the \csu\ obtained from a given \csu\ $\,\y\,$ by blowing up 
any (ordered) $\,\ek$-tuple of distinct points. We will use the isomorphic 
identification
\begin{equation}\label{htw}
H_2(\y'\nnh,\bbZ)\,=\,H_2(\y,\bbZ)\times\bbZ^\ek\nh,
\end{equation}
with the convention that, for any oriented \compact\ \rsu\ $\,\x\,$ embedded 
in $\,\y'\nh$, the $\,\bbZ^\ek$ component of its homology class $\,[\x\hs]\,$ in 
the decomposition (\ref{htw}) consists of the ordered $\,\ek$-tuple formed by 
the intersection numbers of $\,\x\,$ with the resulting $\,\ek\,$ exceptional 
divisors in $\,\y'\nh$.
\begin{thm}\label{thsed}Let\/ $\,\x\,$ be a one\diml\ compact complex 
submanifold of a \csu\ $\,\y$, and let\/ $\,c\,$ be the integral of 
$\,c_{\hs1}(\y)\,$ over $\,\x$. If\/ $\,c\ge0\,$ and\/ $\,\y'$ is the \csu\ 
obtained from $\,\y\,$ by blowing up any $\,\ek$-tuple of distinct points, where 
$\,\ek\ge c\hh$, then $\,\x\,$ admits a \tre\ $\,f\,$ in $\,\y'$ such that, 
under the identification {\rm(\ref{htw})}, 
$\,f_*[\x\hs]=([\x\hs],1,\dots,1,0,\dots,0)\,$ with $\,c\ge0\,$ occurrences 
of\/ $\,1\hs$ and\/ $\,\ek-c\hs\ge\hs0$ occurrences of\/ $\,0$.
\end{thm}
The identification $\,H_2(\y'\nh,\bbZ)=\bbZ^{\ek+1}$ used in the next two 
corollaries is nothing else than (\ref{htw}) for $\,\y=\hs\cp^2\nnh$, with 
$\,H_2(\cp^2\nnh,\bbZ)=\bbZ$.
\begin{cor}\label{codss}Given integers\/ $\,d,\ek\,$ with $\,\ek\ge3d\ge3$, let\/ 
$\,\y'$ be the \csu\ obtained from $\,\cp^2$ by blowing up any ordered set 
of $\,\ek$ points. Then the \compact\ orientable surface $\,\x$ of genus 
$\,(d-1)(d-2)/2\,$ admits a \tre\ $\,f:\x\to\y'$ with 
$\,f_*[\x\hs]=(d\hh,\hs1,1,\dots,1,0,\dots,0)$, where $\,1\hs$ occurs 
$\,3d\,$ times and\/ $\,0\,$ occurs $\,\ek-3d\ge0\,$ times.
\end{cor}
\begin{cor}\label{codot}Given positive integers\/ $\,d,j\,$ with $\,j\ge2d+2$, 
let\/ $\,\y\,$ be the \csu\ obtained from $\,\cp^2$ by blowing up any 
$\,j$-tuple of distinct points. Then there exists a totally real oriented\/ 
$\,2$-sphere $\,\x\,$ embedded in\/ $\,\y\,$ such that\/ 
$\,[\x\hs]=(d\hh,\hs d\nh-\nnh1,1,\dots,1,0,\dots,0)$, with $\,1\hs$ occurring 
$\,2d+1$ times and\/ $\,0\hs$ occurring $\,j-2d-2\,$ times.
\end{cor}

\section{A $\,\bbZ_2$ cohomology constraint}\label{zt}
\setcounter{equation}{0}
Given a \cvs\ $\,V\hs$ with $\,\dimc V=n\ge1$, let $\,\hs\text{\rm TR}(V)\,$ 
(or, $\,\hs\text{\rm TR}^+(V)$) denote the set of all \tr\ (or, 
respectively, oriented totally real) vector subspaces of real dimension 
$\,n\,$ in $\,V\nh$. (See Section~\ref{pr}.) Also, let $\,\rp(W)$ (or, 
$\,\hs\text{\rm S}\hs(W)$) be the real projective space (or, sphere) of all 
real lines through $\,0\,$ (or, respectively, rays emanating from $\,0$) in 
any given \rvs\ $\,W$. We have natural mappings 
\begin{equation}\label{ltr}
\mathfrak{L}:\text{\rm TR}(V)\,\to\,\rp(V^{\wedge n})\,,\qquad\quad
\mathfrak{L}^+:\text{\rm TR}^+(V)\,\to\,\,\text{\rm S}\hs(V^{\wedge n})
\end{equation}
sending each $\,W$ in $\,\text{\rm TR}(V)\,$ or $\,\text{\rm TR}^+(V)\,$ to 
the real line/ray containing $\,e_1\wedge\ldots\wedge\,e_n$, where 
$\,e_1,\dots,e_n$ is any basis (or, positive-oriented basis) of $\,W\nnh$, and 
$\,V^{\wedge n}$ is the $\,n$th {\it complex} exterior power of $\,V$. Thus, 
$\,V^{\wedge n}$ is a complex line, while $\,\rp(V^{\wedge n})\,$ 
and $\,\hs\text{\rm S}\hs(V^{\wedge n})\,$ are circles. For a \tri\ 
$\,f\,$ of a \rmf\ $\,\x\,$ in an \ac\ manifold $\,\y\,$ with 
$\,\dimr\x=\dimc\y=n$, the Maslov invariant 
$\,\mathfrak{M}(f)\in\,[\x,E(\y)]\,$ was defined, in Section~\ref{sr}, to be the 
homotopy class of the mapping 
\begin{equation}\label{haf}
\varTheta(f):\x\,\to\,E(\y)
\end{equation}
which sends each $\,x\in\x\,$ to the real line 
$\,\mathfrak{L}(\df_x(T_x\x))\,$ in 
$\,[T_{f(x)}\y]^{\wedge n}$, with $\,\mathfrak{L}\,$ as in (\ref{ltr}) for 
$\,V\nnh=T_{f(x)}\y$. 
If, in addition, $\,\x\,$ is orientable, $\,\mathfrak{M}(f)\,$ can be lifted 
to $\,E^+(\y)$, that is, there exists a naturally distinguished homotopy class 
$\,\hs\mathfrak{M}^+(f)\in[\x\,,\,E^+(\y)]$, whose composite with (\ref{tcp}) 
is $\,\mathfrak{M}(f)$. A representative of $\,\hs\mathfrak{M}^+(f)\,$ is 
$\,\varTheta^+(f):\x\to E^+(\y)\,$ obtained by fixing an orientation of 
$\,\x$, which does not affect the resulting homotopy class, and then assigning 
to any $\,x\in\x\,$ the ray $\,\mathfrak{L}^+(\df_x(T_x\x))\,$ in 
$\,[T_{f(x)}\y]^{\wedge n}\,$ (notation of (\ref{ltr})), with 
$\,\df_x(T_x\x)\,$ oriented via $\,\df_x$. 

Being a homotopy class of mappings, $\,\mathfrak{M}(f)\,$ induces a 
homomorphism 
\begin{equation}\label{hfa}
[\varTheta(f)]_*:\pi_1\x\,\to\,\pi_1[E(\y)]
\end{equation}
of the fundamental groups (with fixed base points). Under the identifications 
(\ref{pzq}) and (\ref{coh}), the homomorphism (\ref{hfa}) coincides with 
$\,\ind(f)\,$ in (\ref{ifh}).
\begin{prop}\label{coobs}Let\/ $\,f:\x\to\y\,$ be a \tri\ of a real\/ 
$\,n$\mfd\/ $\,\x\,$ in an \acm\/ $\,\y\,$ with\/ $\,\dimc\y=n$. Condition 
{\rm(\ref{wow})} then holds for the mapping\/ 
$\,\varTheta=\varTheta(f)\,$ appearing in {\rm(\ref{haf})}.
\end{prop}
\begin{proof}Let $\,\Gamma=\,\text{\rm Ker}\,[w_1(\x)]$, so that $\,\Gamma\,$ 
is the subgroup of index $\,1\,$ or $\,2$ in $\,\pi_1\x\,$ formed by all 
homotopy classes of loops $\,\gamma:S^1\to\x\,$ for which $\,\gamma^*(T\x)\,$ 
is orientable. Orientability of $\,\gamma^*(T\x)\,$ means that the composite 
$\,\varTheta(f)\circ\gamma\,$ can be lifted to a loop in $\,E^+(\y)$, that is, 
its homotopy class lies in $\,\pi_1[E^+(\y)]=\,\text{\rm Ker}\,\mathfrak{w}_1$ 
(see (\ref{wun})). Thus, 
$\,\Gamma=\,\text{\rm Ker}\,(\mathfrak{w}_1\nh\circ[\varTheta(f)]_*)$. Hence 
$\,\mathfrak{w}_1\circ\hs[\varTheta(f)]_*\hs=\,w_1(\x)$, and so (\ref{coh}) 
gives (\ref{wow}) with $\,\varTheta^*=[\mathfrak{M}(f)]^*\nnh$.
\end{proof}
\begin{rem}\label{ncsuf}Let $\,\x,\y\,$ and $\,E(\y)\,$ be a \compact\ 
\rmf, an \acm\/ and, respectively, the principal 
$\,\hs\text{\rm U}\hs(1)$-bun\-dle over $\,\y\,$ defined in Section~\ref{pr}. 
Furthermore, let $\,f:\x\to\y\,$ be a continuous mapping.
\begin{enumerate}
  \def\theenumi{{\rm\roman{enumi}}}
\item Continuous lifts $\,\varTheta:\x\to E(\y)\,$ of $\,f\,$ satisfying 
(\ref{wow}) are nothing else than real-line subbundles of 
$\,f^*[\detc\hs\tm]\,$ isomorphic to $\,\detr\hs T\x$.
\item A lift/subbundle in (i) exists if and only if the complex line bundle 
$\,f^*[\detc\hs\tm]\,$ over $\,\x\,$ is isomorphic to the complexification 
$\,[\detr\hs T\x]^\bbC\nnh$.
\end{enumerate}
In fact, a lift $\,\varTheta\,$ of $\,f\,$ to $\,E(\y)\,$ selects a real line 
$\,\varTheta_x$ in $\,[\detc\hs\tm]_{f(x)}\,$ for each $\,x\in\x$, that is, 
forms a real-line subbundle of $\,f^*[\detc\hs\tm]$, while (\ref{wow}) 
states that this subbundle is isomorphic to $\,\detr\hs T\x$.
\end{rem}
 
\section{Proof of Theorem~\ref{thzer}}\label{pz}
\setcounter{equation}{0}
Let $\,\mathcal{H}\,$ be the set of all homotopy classes of mappings 
$\,\varTheta:\x\to E(\y)$ with (\ref{wow}). Given 
$\,(x,\rho)\in\x\times E(\y)$, let $\,D[x,\rho]\,$ be the set of all injective 
real-lin\-e\-ar operators $\,A:T_x\x\to T_y\y\,$ such that 
$\,y\in\y\,$ is the image of $\,\rho\,$ under the bundle projection 
$\,E(\y)\to\y\,$ and the image $\,A(T_x\x)\,$ is a \tr\ subspace of 
$\,T_y\y\,$ satisfying the condition $\,\mathfrak{L}(A(T_x\x))=\rho\,$ for 
$\,\mathfrak{L}\,$ as in (\ref{ltr}) with $\,V\nnh=T_y\y$. The Lie group 
$\,G_{x,\rho}\,$ of all complex automorphisms $\,B\,$ of $\,T_y\y\,$ with 
$\,\detc B\in\bbR\hs$, for $\,y\,$ as above, now acts on $\,D[x,\rho]\,$ 
simply transitively by the left multiplication, giving rise to a homotopy 
equivalence $\,D[x,\rho]\,\approx\,\,\text{\rm SU}(2)\times\bbZ_2$.

In view of Proposition~\ref{coobs}, the assignment 
$\,f\mapsto\hs\mathfrak{M}(f)\,$ descends to a mapping 
$\,\mathcal{T}\to\mathcal{H}$, where $\,\mathcal{T}\,$ is the set of all 
$\,\simtri$ equivalence classes of \tri s $\,\x\to\y$. To show that 
$\,\mathcal{T}\to\mathcal{H}\,$ is surjective, let $\,\varTheta:\x\to E(\y)\,$ 
be continuous and satisfy (\ref{wow}), and let $\,\varLambda\,$ be (the total 
space of) the bundle over $\,\x\,$ with the fibres 
$\,\varLambda_x=D[x,\varTheta(x)]$, for $\,x\in\x$.

Then $\,\varLambda\,$ has two connected components: they are bundles over 
$\,\x$, and their fibres, homo\-topy-e\-quiv\-a\-lent to 
$\,\hs\text{\rm SU}(2)$, are connected components of the fibres of 
$\,\varLambda$. In fact, as stated at the very end of the last section, we may 
choose a vec\-tor-bun\-dle isomorphism $\,F\,$ between $\,\detr\hs T\x\,$ and 
$\,\varTheta\,$ treated as a real-line bundle over $\,\x$. A selection of one 
connected component of $\,\varLambda_x$, varying continuously with $\,x\in\x$, 
consists of those $\,A\in\varLambda_x$ inducing real-line isomorphisms 
$\,[T_x\x\hs]^{\wedge2}\to\varTheta_x$ equal to $\,F_x$ times a {\it 
positive\/} factor.

Using a $\,\hs\CW$-de\-com\-po\-si\-tion of the surface $\,\x\,$ we now see 
that (a connected component of) $\,\varLambda\,$ admits a continuous section, 
and so surjectivity of $\,\mathcal{T}\to\mathcal{H}\,$ follows from the 
$\,h$-prin\-ci\-ple for \tri s \cite[p.\ 192]{gromov}.

Finally, to show that $\,\mathcal{T}\to\mathcal{H}\,$ is injective, consider 
two \tri s $\,f,f':\x\to\y\,$ with 
$\,\mathfrak{M}(f)=\mathfrak{M}(f')$. We may thus choose a homotopy 
$\,[\hs0,1]\ni t\mapsto \varTheta_t$ between the mappings 
$\,\varTheta_0=\varTheta(f)\,$ and $\,\varTheta_1=\varTheta(f')\,$ defined as 
in the line following (\ref{haf}). 
If $\,\varXi\,$ now is the bundle over $\,\x\times[\hs0,1]\,$ with the fibres 
$\,D[x,\varTheta_t\hskip-.4pt(x)]$, for $\,(x,t)\in\x\times[\hs0,1]$, a 
$\,\hs\CW$-de\-com\-po\-si\-tion argument shows, as before, that (a connected 
component of) $\,\varXi\,$ admits a continuous section which coincides with 
$\,\df\,$ on $\,\x\times\{0\}\,$ and with $\,\df'$ on $\,\x\times\{1\}$. More 
precisely, if an extension of $\,\df\,$ to a section of $\,\varXi\,$ thus 
obtained has a restriction to $\,\x\times\{1\}\,$ lying in a different 
connected component than that containing $\,\df'\nh$, a correction can be made 
by extending the original homotopy to one pa\-ram\-e\-trized by 
$\,t\in[\hs0,2]$, with $\,\varTheta_t=e^{2(t-1)\pi i}\varTheta_1$ for 
$\,t\in[1,2]$ (multiplication in the principal 
$\,\hs\text{\rm U}\hs(1)$-bun\-dle $\,E(\y)$).

Hence, according to the $\,h$-prin\-ci\-ple \cite[p.\ 192]{gromov}, 
$\,f\,$ and $\,f'$ are homotopic through \tri s, which proves 
Theorem~\ref{thzer}.

Note that, as $\,\hs\text{\rm SU}(3)\,$ is $\,2$-connected, the surjectivity 
part of the above argument is still valid when the dimension $\,n=2\,$ is 
replaced by $\,n=3$.

\section{Topological obstructions}\label{tp}
\setcounter{equation}{0}
Every \acm\ carries a natural orientation. Specifically, an $\,n$\diml\ \cvs\ 
$\,V\hs$ ($1\le n<\infty$) becomes an oriented \rvs\ if one 
declares the real basis $\,e_1,ie_1,\dots,e_n,ie_n$ to be positive oriented 
for some (or any) complex basis $\,e_1,\dots,e_n$. With this convention, the 
effect on the orientation of the direct sum operation for complex spaces 
agrees with that for oriented real spaces. For the oriented totally real 
subspace $\,W=\,\spanr\{e_1,\dots,e_n\}\,$ of $\,V\nnh$, with the orientation 
of $\,iW$ obtained using the di\-rect-sum requirement (that 
$\,V=W\oplus\,iW$ as oriented spaces), the isomorphism $\,W\to iW$ of 
multiplication by $\,i\,$ ``multiplies'' the orientation by the sign factor 
$\,(-1)^{n(n-1)/2}$.

Let $\,f\,$ now be a \tri\ of a real $\,n$\mfd\ $\,\x\,$ in an \acm\ $\,\y\,$ 
with $\,\dimc\y=n$. The multiplication by $\,i$ provides an isomorphic 
identification between the tangent bundle $\,\tau\,$ of $\,\x$ (treated as a 
subbundle of $\,f^*\tm$) and the normal bundle $\,\nu\,$ of $\,f$. Thus,
\begin{equation}\label{tnu}
\tau\,=\,(-1)^{n(n-1)/2}\hskip.8pt\nu\,,\hskip9pt\text{\rm where}\hskip8pt
\tau\,=\,\df(T\x)\hskip8pt\text{\rm and}\hskip8pt
\nu\,=\,[f^*\tm]/\tau\,.
\end{equation}
The factor $\,(-1)^{n(n-1)/2}$ represents the orientation if $\,\x\,$ is 
oriented, and is to be ignored otherwise. In the former case, the isomorphism 
$\,\tau\approx\nu\,$ in (\ref{tnu}) is o\-ri\-en\-ta\-tion-pre\-serv\-ing if 
and only if $\,n\equiv0\,$ or $\,n\equiv1\hskip6pt\text{\rm mod}\hskip4pt4$. 

Let $\,f:\x\to\y\,$ be a \tri\ of a $\,\ek$\diml\ \rmf\ $\,\x\,$ in an \acm\ 
$\,\y\,$ with $\,\dimc\y=n$. Thus, $\,\ek\le n\,$ and we have an obvious 
isomorphic identification $\,\hs\spanc\tau\,=\,[T\x]^\bbC$ of complex vector 
bundles over $\,\x$, with $\,\tau\subset f^*\tm\,$ as in (\ref{tnu}), and 
$\,[\hskip5pt]^\bbC$ denoting complexification. In fact, since $\,f\,$ is an 
immersion, $\,\tau\,$ is isomorphic to $\,T\x$, while 
$\,\hs\spanc\tau\,=\,\tau\oplus\,i\tau\,\approx\,\tau^\bbC\,$ as $\,f\,$ is 
\tr. 

If, in addition, $\,\x\,$ is \compact\ and $\,\dimr\x=\dimc\y$, we can rewrite 
the relation $\,\hs\spanc\tau\,=\,[T\x]^\bbC$ as 
$\,\,\spanc[\df(T\x)]=f^*\tm$, which, followed by the operation $\,\detc$, 
leads to natural isomorphic identifications
\begin{equation}\label{ftm}
\text{\rm i)}\hskip7ptf^*\tm\,=\,[T\x]^\bbC\nh,\hskip30pt\text{\rm ii)}
\hskip7ptf^*[\detc\hs\tm]\,=\,[\detr\hs T\x]^\bbC\nh.
\end{equation} 
Here and in the sequel, given a \rmf\ $\,\x\,$ with $\,\dimr\x=n\,$ 
(or, an \acm\ $\,\y\,$ with $\,\dimc\y=n$), we will denote 
by $\,\,\detr\hs T\x\,=\,[T\x]^{\wedge n}$ or 
$\,\,\detc\hs\tm\,=\,[\tm]^{\wedge n}$ the {\it determinant bundle\/} of the 
tangent bundle, that is, its highest real/com\-plex exterior power. Taking 
$\,w_2$ (or $\,c_{\hs1}$) of both sides in (\ref{ftm}.ii) and noting that 
$\,\detr\hs T\x\,$ is trivial if $\,\x\,$ is orientable, we obtain
\begin{equation}\label{wco}
\aligned
\text{\rm a)}\hskip6pt&f^*[w_2(\y)]\,=\,w_1(\x)\smallsmile w_1(\x)\quad
\text{\rm in}\quad H^2(\x,\bbZ_2)\,,\\
\text{\rm b)}\hskip6pt&f^*[c_{\hs1}(\y)]\,=\,0\quad\text{\rm in}\quad 
H^2(\x,\bbZ)\quad\text{\rm whenever}
\quad\x\quad\text{\rm is\ orientable.}\endaligned
\end{equation}
If, in addition, $\,f:\x\to\y\,$ is a totally real {\it embedding} and the 
\compact\ manifold $\,\x\,$ with $\,\dimr\x=\dimc\y=n\,$ is orientable, we 
have 
\begin{equation}\label{dot}
f_*[\x\hs]\,\cdot\,f_*[\x\hs]\,=\,(-1)^{n(n-1)/2}\hskip.8pt\chi(\x)\,,
\end{equation}
where $\,f_*[\x\hs]\in H_n(\y,\bbZ)\,$ corresponds to either fixed orientation 
of $\,\x\,$ and the dot $\,\,\cdot\,\,$ denotes the $\,\bbZ$-val\-ued 
intersection form. In fact, the Euler class $\,\eu(\nu)\,$ of the normal 
bundle $\,\nu=\nu_f$ of any embedding $\,f:\x\to\y$, integrated over $\,\x$, 
yields $\,f_*[\x\hs]\,\cdot\,f_*[\x\hs]$, while for {\it totally real\/} 
embeddings $\,f$, (\ref{tnu}) gives 
$\,\int_\x\eu(\nu)=(-1)^{n(n-1)/2}\chi(\x)$. When $\,\x\,$ is not assumed 
orientable, instead of (\ref{dot}) the same argument gives
\begin{equation}\label{dtf}
f_*[\x\hs]\,\cdot\,f_*[\x\hs]\,
=\,[\chi(\x)\,\text{\rm mod}\hskip4pt2\hs]\,,
\end{equation}
where, this time, $\,f_*[\x\hs]\in H_n(\y,\bbZ_2)\,$ and $\,\,\cdot\,\,$ takes 
values in $\,\bbZ_2$. 
\begin{prop}\label{trchi}Let a \compact, orientable manifold\/ $\,\x\,$ of 
even real dimension\/ $\,n\ge2\,$ admit a totally real embedding in\/ 
$\,\y\,=\,\cp^n$ or in the complex manifold\/ 
$\,\y\,=\,\cp^n\,\#\,\,\overline{\cp^n}$ obtained from\/ $\,\cp^n$ by blowing 
up a point.
\begin{enumerate}
  \def\theenumi{{\rm\roman{enumi}}}
\item If\/ $\,\y\,=\,\cp^n$, then\/ $\,(-1)^{n/2}\chi(\x)\ge0\,$ and either\/ 
$\,\chi(\x)\equiv0\hskip6pt\text{\rm mod}\hskip4pt4$ or\/ 
$\,\chi(\x)\equiv1\hskip6pt\text{\rm mod}\hskip4pt4$. 
\item If\/ $\,\y\,=\,\cp^n\,\#\,\,\overline{\cp^n}$, then\/ $\,\chi(\x)\,$ 
is either odd, or divisible by\/ $\,4$. 
\end{enumerate}
\end{prop}
\begin{proof}The (quadratic) intersection form in $\,H^n(\y,\bbZ)\,$ is 
algebraically equivalent to $\,\bbZ\ni p\mapsto p^2$ or 
$\,\bbZ\oplus\bbZ\ni(p,q)\mapsto p^2-q^2$. As $\,p^2=pp\,$ and 
$\,p^2-q^2=(p+q)(p-q)$, with both factors even or both odd, our claim follows 
from (\ref{dot}).
\end{proof}
\begin{cor}\label{notrs}If\/ $\,n\,$ is even, $\,S^n$ admits no totally real 
embedding in\/ $\,\cp^n$ or\/ $\,\cp^n\,\#\,\,\overline{\cp^n}$. 
\end{cor}
As shown in \cite{ahern-rudin}, an analogous claim fails for $\,n=3$.

Finally, given a continuous mapping $\,f:\x\to\y\,$ of a \compact\ \rsu\ 
$\,\x\,$ into an \acsu\ $\,\y$, conditions (\ref{wco}) are not only {\em 
necessary} for $\,f\,$ to be homotopic to a totally real immersion 
$\,\x\to\y$, but also {\em sufficient}. This is clear from Theorem~\ref{thzer} 
and Remark~\ref{ncsuf}, since $\,c_{\hs1}$ and $\,w_2$ classify complex line 
bundles over surfaces \cite[p.~798]{whitney-topo}.

\section{Degrees modulo $\hs\,q\hs\,$ for circle\hs-val\-ued 
mappings}\label{dq}
\setcounter{equation}{0}
Given $\,\,q\in\{1,2,3,\dots,\infty\}$, a manifold $\,\,\x$, and a continuous 
mapping \hbox{$g:\x\to\,\text{\rm U}\hs(1)$,} we define 
$\,[\hs g\,\text{\rm mod}\hskip2.5ptq\hs]\in H^1(\x,\bbZ_q)
=\hs\h(\pi_1\x,\bbZ_q)\,$ to be the composite $\,\pi_1\x\to\bbZ\to\bbZ_q$ of 
the action of $\,g\,$ on the fundamental groups with the projection 
$\,\bbZ\to\bbZ_q$, where $\,\bbZ_1=\{0\}\,$ and $\,\bbZ_\infty=\bbZ\hs$.
\begin{enumerate}
  \def\theenumi{{\rm\roman{enumi}}}
\item $\,[\hs g\,\text{\rm mod}\hskip2.5ptq\hs]$, as a homomorphism 
$\,\pi_1\x\to\bbZ_q$, sends the homotopy class of any loop $\,S^1\to\x\,$ to 
the remainder modulo $\,q\,$ of the degree of the composite 
$\,S^1\to\x\to\,\text{\rm U}\hs(1)$, in which the loop is followed by $\,g$.
\item $\,[\hs g\,\text{\rm mod}\hskip2.5ptq\hs]=0\,$ if and only if either 
$\,q=\infty\,$ and $\,g\,$ has a lift $\,\x\to\bbR$ to the universal 
covering of $\,\hs\text{\rm U}\hs(1)$, or $\,q<\infty\,$ and the $\,q$th 
root of $\,g$ treated as a complex-val\-ued function, with 
$\,\hs\text{\rm U}\hs(1)=S^1\subset\bbC\hs$, has a single-val\-ued continuous 
branch $\,\x\to\,\text{\rm U}\hs(1)\subset\bbC\hs$.
\end{enumerate}
In fact, (i) is obvious, and easily implies (ii).

Let $\,E\,$ now be any principal 
$\,\hs\text{\rm U}\hs(1)$-bun\-dle over a \sc\ manifold $\,\y$. As (\ref{seq}) 
is exact, $\,\pi_1E=\bbZ_q$ for some $\,q\in\{1,2,3,\dots,\infty\}$, with 
$\,\bbZ_1=\{0\}\,$ and $\,\bbZ_\infty=\bbZ\hs$. Given a manifold $\,\x\,$ 
and a continuous mapping $\,\varTheta:\x\to E$, let us define 
$\,\jnd(\varTheta)\in H^1(\x,\bbZ_q)\,$ to be the homomorphism of the 
fundamental groups induced by $\,\varTheta\,$ (cf.\ (\ref{coh})). Then, for 
any continuous mapping $\,g:\x\to\,\text{\rm U}\hs(1)$, with 
$\,[\hs g\,\text{\rm mod}\hskip2.5ptq\hs]\,$ defined above, we have
\begin{equation}\label{eqv}
\jnd(g\varTheta)\,=\,\hs\jnd(\varTheta)\,
+\,[\hs g\,\text{\rm mod}\hskip2.5ptq\hs]\hs,
\end{equation}
$g\varTheta\,$ being the value\-wise product. (This is clear from (i) in 
Section~\ref{dq}, since a principal $\,\hs\text{\rm U}\hs(1)$-bun\-dle over 
$\,S^1$ is trivial.)

\section{Proof of Theorem~\ref{thuno}\hs: surjectivity}\label{sj}
\setcounter{equation}{0}
In view of Theorem~\ref{thzer}, the assertion of Theorem~\ref{thuno} amounts 
to bijectivity of the assignment given, in the notations of Section~\ref{dq} and 
(\ref{dgf}), by
\begin{equation}\label{tid}
[\varTheta]\mapsto(\ind\hs,\dg)\hskip9pt\text{\rm with}\hskip7pt\ind
=\jnd(\varTheta)\hskip6pt\text{\rm and}\hskip6pt\dg=\dg(\pi\circ\varTheta)\hs.
\end{equation}
Explicitly, (\ref{tid}) sends the homotopy class of any mapping 
$\,\varTheta:\x\to E(\y)\,$ with (\ref{wow}) to the pair $\,(\ind,\dg)\,$ 
formed by $\,\ind\in H^1(\x,\bbZ_q)=\hs\h(\pi_1\x,\bbZ_q)\,$ which is the 
action of $\,\varTheta\,$ on the fundamental groups, and $\,\dg\,$ which is 
the image of the fundamental homology class 
$\,[\x\hs]\in H_2(\x,\bbZ_{[2]})\,$ under $\,\pi\circ\varTheta$. (Here 
$\,\bbZ_{[2]}$ is defined as in (\ref{ztw}), and $\,\pi:E(\y)\to\y\,$ stands 
for the bundle projection.) Recall that $\,\ind\,$ and $\,\dg\,$ in 
(\ref{tid}) are the Maslov index and degree of a \tri\ $\,f\,$ of $\,\x\,$ in 
$\,\y\,$ whose $\,\simtri$ equivalence class corresponds to $\,[\varTheta]\,$ 
as in Theorem~\ref{thzer}, since $\,\pi\circ\varTheta\,$ then is homotopic to 
$\,f$.

In this section we prove the surjectivity part of Theorem~\ref{thuno}. 
Injectivity will be established in Section~\ref{ij}. 
We begin with a lemma:
\begin{lem}\label{ifidf}Given a \sc\ \acsu\ $\,\y\,$ and a \compact\ real 
surface $\,\x$, let\/ $\,q,\hs\fri_q(\x),\hs\depspm(\y),\ve\,$ and\/ $\,\pm\,$ 
be defined by {\rm(\ref{pzq})} and\/ {\rm(\ref{dfd})} -- {\rm(\ref{iqs})}.
\begin{enumerate}
  \def\theenumi{{\rm\alph{enumi}}}
\item $(\ind\hs,\dg)\in\,\fri_q(\x)\times\depspm(\y)\,$ if\/ 
$\,(\ind\hs,\dg)\,$ corresponds as in {\rm(\ref{tid})} to a continuous mapping 
$\,\varTheta:\x\to E(\y)\,$ satisfying {\rm(\ref{wow})}.
\item $\ind(f)\in\fri_q(\x)\subset H^1(\x,\bbZ_q)\,$ and\/ 
$\,\dg(f)\in\depspm(\y)\subset H_2(\y,\bbZ_{[2]})\,$ whenever $\,f:\x\to\y\,$ 
is a \tri\ and\/ $\,\hs\ind(f),\hs\dg(f)\,$ are its Maslov index and degree.
\end{enumerate}
\end{lem}
\begin{proof}As $\,\pi_1[E^+(\y)]=2\bbZ_q\subset\bbZ_q\,$ by (\ref{pie}) -- 
(\ref{pzq}), $\,\hs\mathfrak{w}_1$ in (\ref{wun}) must be the unique nonzero 
homomorphism $\,\bbZ_q\to\bbZ_2$. For $\,\varTheta\,$ as in (a), relation 
(\ref{wow}) states that $\,\mathfrak{w}_1\circ\varTheta_*\hs=\,w_1(\x)$, and 
so $\,w_1(\x)\,$ is the $\,\hs\text{\rm mod}\hskip4pt2\,$ reduction of 
$\,\hs\ind\hs\,$ (see Remark~\ref{reduc}). Hence $\,\hs\ind\in\fri_q(\x)$. 
Next, $\,\hs\dg(f)\hs\in\,\depspm(\y)\,$ in view of Remark~\ref{ncsuf}(ii) 
for $\,f=\pi\circ\varTheta\,$ and (\ref{wus}), which proves (a). Finally, 
(b) follows from (a) and Proposition~\ref{coobs}, completing the proof.
\end{proof}
Let $\,\mathcal{Z}\,$ be the image of the mapping (\ref{tid}). Thus, 
$\,\mathcal{Z}\subset\fri_q(\x)\times\depspm(\y)$ by Lemma~\ref{ifidf}(a). To 
prove surjectivity in Theorem~\ref{thuno}, we show that $\,\mathcal{Z}\,$ is 
the set $\,\mathcal{Z}(\x,\y)\,$ defined immediately after 
Theorem~\ref{thuno}. First, if $\,\depspm(\y)=\hh\emp\hh$, our claim follows 
as $\,\mathcal{Z}=\fri_q(\x)\times\depspm(\y)=\hh\emp\hh$.

From now on, we may thus assume that $\,\depspm(\y)\ne\hh\emp\hh$. Let us now 
fix any $\,\dg\in\,\depspm(\y)$, and describe the set of all 
$\,\ind\in\,\fri_q(\x)\,$ with $\,(\ind\hs,\dg)\in\mathcal{Z}$.

Since $\,\y\,$ is \sc, $\,\hs\dg\hs\,$ (or, in fact, any class in 
$\,H_2(\y,\bbZ_{[2]})$) is realized by a mapping $\,S^2\to\y\,$ and, with the 
aid of a degree $\,1\,$ map $\,\x\to S^2\nnh$, also by a mapping 
$\,f:\x\to\y$. However, $\,\hs\dg\hs=\hs\dg(f)\in\depspm(\y)$, and so 
relations (\ref{dmk}), (\ref{dfd}) and (\ref{wus}) show that either 
$\,\x\,$ is orientable and $\,f^*[\detc\hs\tm]\,$ is trivial, or $\,\x\,$ is 
nonorientable and the line bundles $\,f^*[\detc\hs\tm]\,$ and 
$\,[\detr\hs T\x]^\bbC$ have the same $\,w_2$. Thus, in either case, 
$\,f^*[\detc\hs\tm]\,$ and $\,[\detr\hs T\x]^\bbC$ are isomorphic 
\cite[p.~798]{whitney-topo}; hence, by Remark~\ref{ncsuf}(ii), $\,f\,$ admits 
a continuous lift $\,\varTheta:\x\to E(\y)\,$ with (\ref{wow}).

Once such $\,\varTheta\,$ is fixed, those $\,\ind\in\,\fri_q(\x)\,$ for which 
$\,(\ind\hs,\dg)\in\mathcal{Z}\,$ are precisely the values 
$\,\jnd(g\varTheta)\,$ for all continuous mappings 
$\,g:\x\to\,\text{\rm U}\hs(1)\,$ (notation of Section~\ref{dq}) with the 
property that the homomorphism 
$\,[\hs g\,\text{\rm mod}\hskip2.5ptq\hs]:\pi_1\x\to\bbZ_q$ defined in 
Section~\ref{dq} is valued in the {\it even\/} subgroup $\,2\bbZ_q$ of 
$\,\bbZ_q$. In fact, $\,\hs\dg\hs\,$ determines the homotopy class of 
$\,f:\x\to\y\,$ uniquely \cite{whitney-class}, which, combined with an obvious 
homotopy-lifting argument, shows that the elements of $\,\mathcal{Z}\,$ 
having the form $\,(\ind\hs,\dg)\,$ are images under (\ref{tid}) of products 
$\,g\varTheta$. However, as $\,\varTheta$ satisfies (\ref{wow}), condition 
(\ref{wow}) for $\,g\varTheta\,$ (rather than $\,\varTheta$) is equivalent to 
even-val\-ued\-ness of $\,[\hs g\,\text{\rm mod}\hskip2.5ptq\hs]$, as 
$\,\varTheta^*$ in (\ref{wow}) is dual to $\,\jnd(\varTheta)$.

The set $\,\{\ind\in\,\fri_q(\x):(\ind\hs,\dg)\in\mathcal{Z}\}\,$ is thus 
nonempty (as it contains $\,\jnd(\varTheta)$), and hence, by (\ref{eqv}), it 
is a coset, in $\,\hs\h(\pi_1\x,2\bbZ_q)$, of the subgroup 
$\,\mathcal{G}\,$ consisting of those elements of 
$\,\hs\h(\pi_1\x,2\bbZ_q)\,$ which have the form 
$\,[\hs g\,\text{\rm mod}\hskip2.5ptq\hs]$, where 
$\,g:\x\to\,\text{\rm U}\hs(1)\,$ is continuous. Therefore $\,\hs\fri_q(\x)\,$ 
itself is nonempty and, by (\ref{iqs}) and Remark~\ref{reduc}, it is a coset 
of $\,\hs\h(\pi_1\x,2\bbZ_q)\,$ in $\,\hs\h(\pi_1\x,\bbZ_q)=H^1(\x,\bbZ_q)$.

First, suppose that $\,\x\,$ is orientable, or $\,q=\infty$, or $\,\chi(\x)\,$ 
is odd, or $\,q$ is finite but not divisible by $\,4$. In these four cases, 
$\,\mathcal{G}=\,\h(\pi_1\x,2\bbZ_q)$ (and hence 
$\,\mathcal{Z}=\,\fri_q(\x)\times\depspm(\y)\,$ according to the last 
paragraph). In fact, for every \compact\ surface $\,\x$, continuous mappings 
$\,g:\x\to\,\text{\rm U}\hs(1)\,$ realize all homomorphisms 
$\,\pi_1\x\to\bbZ=\pi_1[\hs\text{\rm U}\hs(1)]$, while, in each of the four 
cases, a homomorphism $\,H_1(\x,\bbZ)\to\bbZ_q$ valued in $\,2\bbZ_q$ is 
necessarily the composite $\,H_1(\x,\bbZ)\to\bbZ\to\bbZ_q$ of some 
homomorphism $\,H_1(\x,\bbZ)\to\bbZ\,$ and the projection $\,\bbZ\to\bbZ_q$. 
In the first two cases, this is obvious since $\,H_1(\x,\bbZ)$ is free, or, 
respectively, $\,\bbZ_\infty\nh=\bbZ$. In the last two cases, $\,\x\,$ may 
thus be assumed nonorientable, and $\,q\,$ finite; then $\,H_1(\x,\bbZ)\,$ has 
just one nontrivial element $\,\xi\,$ of finite order, namely, of order $\,2$. 
(See Remark~\ref{wones}.) Case three now implies case four: if $\,\chi(\x)\,$ 
is odd, $\,\xi\notin\,\text{\rm Ker}\,[w_1(\x)]$, and so the image of 
$\,\xi\,$ under any homomorphism that lies in $\,\hs\fri_q(\x)\,$ necessarily 
equals $\,q/2$, due to its being odd in $\,\bbZ_q$ and of order two; hence 
$\,q/2\,$ is odd as $\,\fri_q(\x)\ne\hh\emp\hh$. Next, in the fourth case, any 
homomorphism $\,H_1(\x,\bbZ)\to2\bbZ_q$ sends $\,\xi\,$ to $\,0$, since 
$\,q/2$, the only nontrivial element of order $\,2\,$ in $\,\bbZ_q$, is odd.

Finally, let us assume that $\,\x\,$ is nonorientable, $\,\chi(\x)\,$ is 
even, and $\,q\,$ is a finite multiple of $\,4$. A homomorphism 
$\,\varphi:H_1(\x,\bbZ)\to2\bbZ_q$ sends the torsion element 
$\,\xi\in H_1(\x,\bbZ)\,$ either to $\,0\,$ or to $\,q/2$, and only those 
homomorphisms $\,\varphi\,$ with $\,\varphi(\xi)=0\,$ have factorizations 
$\,H_1(\x,\bbZ)\to\bbZ\to\bbZ_q$ as above, that is, lie in $\,\mathcal{G}$. 
Hence $\,\mathcal{G}\,$ is an index $\,2\,$ subgroup of 
$\,\hs\h(\pi_1\x,2\bbZ_q)$, which shows that, in this case, $\,\mathcal{Z}\,$ 
has half the number of elements of $\,\hs\fri_q(\x)\times\depspm(\y)$. 
The relation $\,\mathcal{Z}=\mathcal{Z}(\x,\y)\,$ will now follow once we show 
that the image of $\,\xi$ under some (or any) $\,\hs\ind\hs\,$ with 
$\,(\ind\hs,\dg)\in\mathcal{Z}\,$ equals $\,0\,$ if and only if our 
$\,\hs\dg\hs\,$ is the $\,\hs\text{\rm mod}\hskip4pt2\,$ reduction of an 
element of $\,\hs\text{\rm Ker}\,[c_{\hs1}(\y)]\subset H_2(\y,\bbZ)$.

To establish the `only if' part of this last statement, let some such 
$\,\hs\ind\hs\,$ send $\,\xi\,$ to $\,0$. For a suitably chosen embedded 
circle $\,\varGamma\subset\x\,$ representing $\,\xi\,$ in homology, 
$\,\x\smallsetminus\varGamma\,$ is the interior of a compact orientable 
surface with a boundary formed by two circles. (Cf.\ case (b) of 
Remark~\ref{srfcs} in Section~\ref{ch}). Capping the two circles with two copies of 
a $\,2$-disk 
$\,D\,$ we obtain a \compact\ orientable surface $\,\x\hs'$ such that $\,\x\,$ 
is homeomorphic to $\,\x\hs'\#\,K^2\nnh$, where $\,K^2$ is the \kb. A mapping 
$\,\varTheta:\x\to E(\y)\,$ with 
(\ref{wow}) that realizes $\,(\ind\hs,\dg)\,$ as in (\ref{tid}) now gives rise 
to a mapping $\,\varTheta\hs':\x'\nh\to E(\y)\,$ equal to $\,\varTheta$ on 
$\,\x\smallsetminus\varGamma\,$ and obtained on both copies of $\,D\,$ by 
extending $\,\varTheta\,$ from $\,\varGamma\,$ to $\,D$, which is possible as 
$\,\hs\ind\hs\,$ sends $\,\xi\,$ to $\,0$. In addition, 
$\,\varTheta\hs'$ still satisfies (\ref{wow}): every element of 
$\,H_1(\x\hs'\nnh,\bbZ)\,$ is represented by a loop 
$\,\gamma:S^1\to\x\smallsetminus\varGamma$, for which $\,\gamma^*(T\x)\,$ 
is orientable (since so is $\,\gamma^*(T\x\hs')$), and, therefore, 
$\,\varTheta\hs'\nnh\circ\gamma=\varTheta\circ\gamma$ can be lifted to a 
loop in $\,E^+(\y)$. Thus, by Lemma~\ref{ifidf}(a), 
$\,\hs\dg\hs'=\hs\dg(\pi\circ\varTheta)$ 
lies in $\,\donep(\y)=\,\text{\rm Ker}\,[c_{\hs1}(\y)]$, and 
the $\,\hs\text{\rm mod}\hskip4pt2\,$ reduction of $\,\hs\dg\hs'$ is 
$\,\hs\dg\hs$, due to mutual cancellation of the contibutions from the two 
copies of $\,D$.

For the `if' part, let $\,\hs\dg\hs\,$ be the 
$\,\hs\text{\rm mod}\hskip4pt2\,$ reduction of 
$\,\hs\dg\hs'\in\,\hs\text{\rm Ker}\,[c_{\hs1}(\y)]$, and let $\,\x\hs'$ be a 
\compact\ orientable \rsu\ with $\,\x\hs'\#\,K^2\nh=\x$. For reasons given in 
the third paragraph after the proof of Lemma~\ref{ifidf} (but now applied to 
$\,\hs\dg\hs'\nnh$, $\,\x\hs'$ and 
$\,\donep(\y)=\,\text{\rm Ker}\,[c_{\hs1}(\y)]\,$ rather than $\,\hs\dg\hs$, 
$\,\x\,$ and $\,\depspm(\y)$), one can realize $\,\hs\dg\hs'$ by a continuous 
mapping $\,f':\x\hs'\to\y$, and any such $\,f'$ admits a continuous lift 
$\,\varTheta\hs':\x\hs'\to E(\y)\,$ satisfying (\ref{wow}). Next, let 
a mapping $\,\varTheta\hs''\nnh:K^2\to E(\y)\,$ from the \kb\ into $\,E(\y)\,$ 
be obtained as the composite of the bundle projection $\,K^2\to S^1$ 
followed by a homeomorphism of $\,S^1$ onto a fibre, intersecting 
$\,\varTheta\hs'(\x\hs')$, of the principal $\,\hs\text{\rm U}\hs(1)$-bun\-dle 
$\,E(\y)$. Without changing the homotopy class of $\,\varTheta\hs'$ or 
$\,\varTheta\hs''\nnh$, we may further assume that they map some 
small nonempty open sets $\,\,U'\nh\subset\x\hs'$ and 
$\,\,U''\nnh\subset K^2$ onto a single point of $\,E(\y)$. Modifying both 
mappings in closed $\,2$-disks $\,D\hs'\subset\hs U'$ and 
$\,D\hs''\nh\subset\hs U''\nnh$, we obtain a mapping 
$\,\varTheta:\x=\x\hs'\#\,K^2\nh\to E(\y)$ that is constant on the tube
connecting $\,\x\hs'\nnh\smallsetminus D\hs'$ to 
$\,K^2\nh\smallsetminus D\hs''$ in $\,\x\hs'\#\,K^2\nnh$, while 
$\,\varTheta=\varTheta\hs'$ on $\,\x\hs'\nnh\smallsetminus D\hs'$ and 
$\,\varTheta=\varTheta\hs''$ on $\,K^2\nh\smallsetminus D\hs''\nnh$. It is now 
clear that $\,\varTheta\,$ sends $\,\xi$ to $\,0\,$ and satisfies 
(\ref{wow}): $\,\xi\,$ is represented by a fibre of the bundle projection 
$\,K^2\to S^1\nnh$, on which $\,\varTheta\,$ is constant, while $\,\xi\,$ and 
the homology classes of loops lying in $\,\x\hs'\nnh\smallsetminus D\hs'$ 
generate the kernel of 
$\,w_1(\x):H_1(\x,\bbZ)\to\bbZ_2$, so that the 
$\,\varTheta$-im\-ages of these loops can be lifted to loops in $\,E^+(\y)\,$ 
due to orientability of $\,\x\hs'\nh$. However, no such lifts exist for a 
circle in $\,K^2$ forming a section of the bundle $\,K^2\to S^1\nnh$, as its 
$\,\varTheta$-im\-age is a fibre of $\,E(\y)\,$ (and so any lift is a {\it 
semicircle\/} in a fibre of $\,E^+(\y)$).

\section{Proof of Theorem~\ref{thuno}\hs: injectivity}\label{ij}
\setcounter{equation}{0}
To prove injectivity of (\ref{tid}), let two mappings 
$\,\varTheta,\varTheta\hs':\x\to E(\y)$, both satisfying (\ref{wow}), have 
$\,\jnd(\varTheta)=\,\jnd(\varTheta\hs')=\,\ind\hs\,$ and 
$\,\dg(\pi\circ\varTheta)=\,\dg(\pi\circ\varTheta\hs')=\,\dg\hs$. As 
$\,\hs\dg(f)=f_*[\x\hs]\in H_2(\y,\bbZ_{[2]})$ uniquely determines the 
homotopy class of $\,f:\x\to\y\,$ (see \cite{whitney-class}), we can lift a 
fixed homotopy between $\,f=\pi\circ\varTheta\,$ and 
$\,f'=\pi\circ\varTheta\hs'$ to the principal 
$\,\hs\text{\rm U}\hs(1)$-bun\-dle $\,E(\y)\,$ over $\,\y$, obtaining a 
homotopy between $\,\varTheta\hs'$ and some lift $\,g\varTheta\,$ of 
$\,f=\pi\circ\varTheta$ to $\,E(\y)$, where 
$\,g:\x\to\,\text{\rm U}\hs(1)\,$ is a suitable continuous mapping. By 
(\ref{eqv}), $\,[\hs g\,\text{\rm mod}\hskip2.5ptq\hs]=0$. It now suffices to 
show that $\,\varTheta\,$ and $\,g\varTheta\,$ are homotopic. We may assume 
that $\,q<\infty$, for otherwise $\,g\,$ is homotopic to a constant mapping 
((ii) in Section~\ref{dq}).

The required homotopy $\,\varXi:\x\times[\hs0,1]\to E(\y)\,$ with 
$\,\varXi(\,\cdot\,,0)=\varTheta\,$ and $\,\varXi(\,\cdot\,,1)=g\varTheta\,$ 
will be built on successive skeleta of a specific 
$\,\hs\CW$-de\-com\-po\-si\-tion of $\,\x$, namely, one resulting from a 
surjective continuous mapping $\,\phi:D\to\x$, where $\,D\,$ is an oriented 
$\,2$-disk with the accordingly oriented boundary circle $\,\partial D$. We 
choose $\,\phi\,$ so that $\,\partial D\,$ is partitioned by some 
$\,2\ek$-el\-e\-ment subset, $\,\ek\ge1$, into $\,2\ek\,$ compact {\it boundary 
segments\/} and the only identifications of points of $\,D\,$ under $\,\phi\,$ 
are those provided by some $\,\ek\,$ homeomorphisms between pairs of boundary 
segments. Replacing $\,\varXi\,$ by 
$\,\hat\varXi:D\times[\hs0,1]\to E(\y)\,$ with 
$\,\hat\varXi(z,t)=\varXi(\phi(z),t)\,$ for all 
$\,(z,t)\in D\times[\hs0,1]$, we see that, at any stage, instead of 
$\,\varXi\,$ we may just construct $\,\hat\varXi\,$ which depends on 
$\,(z,t)\,$ only through $\,(\phi(z),t)$.

First, we choose a single-val\-ued continuous function 
$\,h:\x\to\,\text{\rm U}\hs(1)\subset\bbC$ with $\,h^q=g$. (See (ii) in 
Section~\ref{dq}.) Now, for any point $\,x\,$ in the $\,0$-skel\-e\-ton of 
$\,\x\,$ (the $\,\phi$-im\-age of the $\,2\ek$-el\-e\-ment partitioning set), we 
let $\,t_j=j/q\,$ and $\,\varXi(x,t)=[h(x)]^j\gamma_x(t-t_j)\varTheta(x)\,$ 
for $\,t\in[t_j,t_{j+1}]\,$ and $\,j\in\{0,\dots,q-1\}$, where 
$\,\gamma_x:[\hs0,1/q\hs]\to\,\text{\rm U}\hs(1)\,$ is any fixed curve joining 
$\,1\,$ to $\,h(x)$.

An extension of $\,\hat\varXi\,$ from the $\,0$-skel\-e\-ton to the 
$\,1$-skel\-e\-ton is in turn obtained separately on each rectangle 
$\,R=S\times[\hs0,1]$, where $\,S\subset\partial D\,$ is one of the $\,2\ek\,$ 
boundary segments. Our $\,\hat\varXi\,$  is already defined on the (oriented) 
boundary $\,\partial R\,$ of $\,R$, so that 
$\,\hat\varXi:\partial R\to E(\y)\,$ represents a free homotopy class 
$\,\alpha\,$ of loops in $\,E(\y)$. As $\,\pi_1[E(\y)]\nh=\bbZ_q$ is Abelian, 
such free homotopy classes can be meaningfully multiplied, and form a group 
isomorphic to $\,\pi_1[E(\y)]$. We now show that $\,\alpha=\beta\hh^q$ for the 
free homotopy class $\,\beta\,$ of some loop $\,\varPi$, and so $\,\alpha\,$ 
is trivial, which provides an extension of $\,\hat\varXi\,$ from 
$\,\partial R$ to $\,R$, thus concluding the $\,1$-skel\-e\-ton step. 
Specifically, $\,\varPi:\partial R\hh'\to E(\y)$ is defined on the boundary of 
$\,R\hh'=S\times[\hs0,1/q\hs]\,$ by 
$\,\varPhi(\,\cdot\,,1/q)=(h\varTheta)\circ\phi$ on $\,S\times\{1/q\}\,$ and 
$\,\varPi=\hat\varXi\,$ everywhere else. Also, 
$\,\alpha=\beta_0\ldots\beta_{q-1}$, with $\,\beta_j$ denoting the free 
homotopy class of the loop $\,\varPi_j:\partial R\hh'\to E(\y)\,$ given by 
$\,\varPi_j(\,\cdot\,,t)=(h^j\circ\phi)\hs\varPi(\,\cdot\,,t)$, as one sees 
noting that for each $\,j\in\{0,\dots,q-2\}$ the contribution to $\,\beta_j$ 
from $\,S\times\{1/q\}\,$ cancels the contribution to $\,\beta_{j+1}$ from 
$\,S\times\{0\}$. On the other hand, $\,\beta_j=\beta\,$ for all $\,j$, since 
the segment $\,S\,$ is contractible and so we may choose a homotopy between 
the constant mapping $\,1:S\to\,\text{\rm U}\hs(1)\,$ and 
$\,h:S\to\,\text{\rm U}\hs(1)$, which leads to a free homotopy between 
$\,\varPi_j$ and $\,\varPi_{j+1}$, $\,j=0,\dots,q-2$.

Now that $\,\hat\varXi\,$ is already defined on the boundary $\,\partial C\,$ 
of the solid cylinder $\,C=D\times[\hs0,1]$, our  $\,2$-skel\-e\-ton step amounts to extending it from $\,\partial C\,$ to $\,C$. 
Let $\,\sigma\in\pi_2[E(\y)]\,$ and $\,\xi\in H_2(E(\y),\bbZ)\,$ be the 
homotopy class with any fixed base point, and the homology class, of 
$\,\hat\varXi:\partial C\to E(\y)\,$ (for the standard orientation of the 
$\,2$-sphere $\,\partial C$). There are two cases.

If $\,\x\,$ is orientable, $\,\sigma=0$. In fact, $\,\xi\,$ then equals the 
difference of the homology classes of $\,\varTheta\,$ and $\,g\varTheta$, for 
a suitable orientation of $\,\x$. (Contributions from $\,\hat\varXi\,$ 
restricted to $\,\partial D\times[\hs0,1]\,$ undergo pairwise 
cancellations, as all identifications, under $\,\phi$, of pairs of boundary 
segments in $\,\partial D\,$ are o\-ri\-en\-ta\-tion-re\-vers\-ing.) Since 
$\,\pi\circ\varTheta=\pi\circ(g\varTheta)$, we thus get 
$\,\pi_*\hh\xi=0\,$ in $\,H_2(\y,\bbZ)$, so that, in view of the Hurewicz 
isomorphism, $\,\pi_*\hh\sigma=0\,$ in $\,\pi_2\y$. Injectivity of the first 
homomorphism in (\ref{seq}) now gives $\,\sigma=0$, which provides the 
required extension of $\,\hat\varXi\,$ from $\,\partial C\,$ to the $\,3$-disk 
$\,C$.

If $\,\x\,$ is not orientable, $\,\sigma\,$ may depend on how one chose the 
extension on $\,\varXi\,$ from the $\,0$-skel\-e\-ton to the 
$\,1$-skel\-e\-ton of $\,\x$. Namely, if we choose that extension differently, 
we can modify the resulting $\,\sigma\,$ so as to add to it any prescribed 
even element of $\,\pi_2[E(\y)]$, `even' meaning divisible by $\,2$ in 
$\,\pi_2[E(\y)]$. In fact, at least one identification under $\,\phi\,$ of a 
pair  $\,S,S\hs'$ of boundary segments in $\,\partial D\,$ is now 
o\-ri\-en\-ta\-tion-pre\-serv\-ing. Let us fix such $\,S,S\hs'\subset\partial D$, 
with $\,\phi(S)=\phi(S\hs')$. Any 
given element $\,\rho\,$ of $\,\pi_2[E(\y)]\,$ can be represented by a mapping 
$\,F\,$ from the $\,2$-sphere obtained when two separate copies of the 
rectangle $\,R=S\times[\hs0,1]\,$ are glued together by identifying their 
boundaries $\,\partial R$, and $\,F\,$ may be chosen so that the restriction 
of $\,F\,$ to one copy of $\,R\,$ is the extension of $\,\hat\varXi\,$ from 
$\,\partial R\,$ to $\,R\,$ which we used in the $\,1$-skel\-e\-ton step. If 
we now replace that extension by a new one, namely, by the restriction of 
$\,F\,$ to {\it the other copy\/} of $\,R$, the corresponding element 
$\hs\sigma$ of $\,\pi_2[E(\y)]\,$ will be replaced by $\,\sigma+2\rho$. 
(As the identification of $\,S\,$ with $\,S\hs'$ under $\,\phi\,$ is 
o\-ri\-en\-ta\-tion-pre\-serv\-ing, the restriction of either version of 
$\,\hat\varXi$ to $\,S\times[\hs0,1]\,$ contributes twice to the 
homotopy or homology class.)

Therefore, to show that $\,\sigma=0\,$ for a suitably chosen extension to the 
$\,1$-skel\-e\-ton, we just need to verify that the original $\,\sigma\,$ 
is an even element of $\,\pi_2[E(\y)]$. To this end, first note that 
 $\,\pi_*\hs\xi\,$ is even in $\,H_2(\y,\bbZ)$, since for some pairs $\,S,S\hs'$ 
of boundary segments identified under $\,\phi\,$ the contributions of 
$\,\hat\varXi:S\times[\hs0,1]\to E(\y)\,$ to homology {\it count twice\/} 
(namely, the pairs whose identification is 
o\-ri\-en\-ta\-tion-pre\-serv\-ing), while for the remaining pairs the 
contributions {\it cancel each other}, which is also the case for the 
contributions of $\,\pi\circ\hat\varXi:D\times\{t\}\to\y\,$ for $\,t=0\,$ and 
$\,t=1\,$ (due to their having the same image, with opposite orientations). 
The Hurewicz isomorphism theorem now implies that $\,\pi_*\hh\sigma\,$ is even 
in $\,\pi_2\y$. If $\,\sigma\,$ itself were {\it not\/} even in 
$\,\pi_2[E(\y)]$, an element of $\,\pi_2\y\,$ whose double is 
$\,\pi_*\hh\sigma\,$ would project onto a nontrivial element of order $\,2\,$ 
in the quotient of $\,\pi_2\y\,$ over the isomorphic image of 
$\,\pi_2[E(\y)]\,$ under $\,\pi_*$ (cf.\ (\ref{seq})). This would contradict 
the fact that, in view of exactness of (\ref{seq}), the quotient group in 
question is isomorphic to a subgroup of $\,\pi_1[\hs\text{\rm U}\hs(1)]=\bbZ$, 
and hence torsion-free. Thus, (\ref{tid}) is injective, which proves 
Theorem~\ref{thuno}.

\section{The simplest examples}\label{se}
\setcounter{equation}{0}
Here begins the second part of the paper, devoted to constructing explicit 
examples of \tri s and embeddings. The simplest such constructions are 
described in this section.
\begin{example}\label{trsub}A real subspace $\,W$ of a \cvs\ $\,V\hs$ 
with $\,\dim V<\infty\,$ is \tr\ \iff\ $\,\hs\spanc W\,$ in $\,V\hs$ 
has the complex dimension $\,\dimr W$. This means that some (or, every) 
$\,\bbR$-basis of $\,W$ is also linearly independent over $\,\bbC\,$ in $\,V$. 
Thus, given complex-val\-ued $\,C^1$ functions $\,f_1,\dots,f_n$ on a nonempty 
open set $\,\,U\,$ in $\,\rn$, the mapping $\,f=(f_1,\dots,f_n)\,$ is a \tri\ 
$\,\,U\to\bbC^n$ \iff\ $\,\mathcal{J}(f_1,\dots,f_n)\ne0\,$ at every 
$\,x\in U$, where $\,\mathcal{J}(f_1,\dots,f_n)=\,\det\,\mathfrak{F}\,$ for 
the complex $\,n\times n\,$ Jacobian matrix 
$\,\hs\mathfrak{F}\,=\,\mathfrak{}F(x)\,$ with the entries 
$\,\hs\partial f_j/\partial x_\ek$. 
\end{example}
\begin{example}\label{trimg}Obviously, an embedding $\,f\,$ of a \rmf\ 
$\,\x\,$ in an \acm\ $\,\y\,$ is \tr\ \iff\ so is the image $\,f(\x)\,$ as a 
submanifold of $\,\y$.
\end{example}
\begin{example}\label{trgph}Let 
$\,\x=\{(z,v)\in U\nnh\times V:v=\varphi(z)\}\,$ be the graph of a 
$\,C^\infty$ mapping $\,\varphi:U\to V\,$ from a nonempty connected open set 
$\,\,U\subset\bbC$ into a complex vector space $\,V\hs$ with 
$\,\dim V<\infty$. Then $\,\x\,$ is a \tr\ submanifold of 
$\,\,U\nnh\times V\,$ \iff\ $\,\varphi\,$ satisfies, at each point of $\,\,U$, 
the {\it Cau\-chy-Rie\-mann inequality\/} $\,\varphi_{\bar z}\ne0$, where 
$\,\varphi_{\bar z}=(\varphi_x\nh+i\varphi_y)/2\,$ with 
$\,x=\,\text{\rm Re}\,z$, $\,y=\,\text{\rm Im}\,z$, and the subscripts stand 
for the partial derivatives.

In fact, $\,\x\,$ is \tr\ in $\,\,U\nnh\times V\,$ \iff\ the graph embedding 
$\,\,U\to\,U\nnh\times V\,$ given by $\,z\mapsto(z,\varphi(z))\,$ is \tr\ 
(Example~\ref{trimg}), and our claim follows in view of Example~\ref{trsub}.
\end{example}
\begin{example}\label{totre}Given a \cvs\ $\,V\hs$ with 
$\,\dimc V=2$, any real subspace $\,W\subset V\,$ with 
$\,\dimr W=2\,$ which is not \tr\ must, obviously, be a complex 
$\,1$\diml\ subspace of\/ $\,V$. Suppose now that $\,\x\,$ is 
a submanifold of an \acm\ $\,\y\,$ and 
$\,\dimr\x=2$. Removing from $\,\x\,$ all {\it complex points}, that 
is, those $\,x\in\x\,$ for which $\,T_x\x\,$ is a complex line in 
$\,T_x\y$, we obtain an open subset $\,\,U\,$ of $\,\x\,$ and, if $\,\,U\,$ is 
nonempty, its connected components are \tr\ submanifolds of $\,\y$. 
\end{example}
Totally real submanifolds naturally arise in many other common situations. For 
instance, a real submanifold $\,\x\subset\y\,$ of an \acm\ $\,\y\,$ is \tr\ in 
each of the following obvious cases:
\begin{enumerate}
  \def\theenumi{{\rm\roman{enumi}}}
\item $\dim\x=1$.  
\item $\x\,$ is a connected component of the fix\-ed-point set 
$\,\{x\in\y:\overline x=x\}$ of any $\,C^\infty$ involution 
$\,\y\ni x\mapsto\overline{x}\in\y\,$ reversing the \acst. One then has 
$\,\dimr\x=\dimc\y$. 
\item $N\,$ is an \acm\ admitting an involution 
$\,x\mapsto\overline{x}$ that reverses the almost complex structure and 
$\,\y\,$ is the product \acm\ $\,N\times\hskip.4ptN$, while 
$\,\x\subset\y\,$ is the {\it an\-ti-di\-ag\-o\-nal submanifold\/} 
$\,\{(\overline x,x):x\in N\}$, \feic\ to $\,N\nh$. (This is clear from 
(ii) applied to the involution $\,(x,y)\mapsto(\overline{y},\overline{x})\,$ 
of $\,N\times\hskip.4ptN$.)
\item $\,\x=\x\hs'\hskip-.3pt\times \x\hs''\nnh$, where $\,\x\hs'\subset\y'$ 
and $\,\x\hs''\subset\y''$ are \tr\ submanifolds, and 
$\,\y\nh=\y'\hskip-.3pt\times\y''\nh\,$ is the product \acm. (For a more 
general construction, see Lemma~\ref{trbdl}.) 
\item By (i) and (iv), embedded closed curves $\,K,K'\subset\bbC\,$ give 
rise to a Clifford-like \tr\ embedded $\,2$-torus 
$\,\x=K\times K'\subset\y=\bbC^2\nnh$. Iterating this produces \tr\ 
embedded $\,n$-tori in $\,\bbC^n$. 
\item The standard {\it real form\/} $\,\x=\rp^n\subset\y\,=\,\cp^n$ is a 
\tr\ embedded submanifold, which follows from of (ii) for the 
involution $\,\cp^n\ni[x_0,\dots,x_n]
\mapsto[\hs\overline x_0,\dots,\overline x_n\hs]\in\cp^n$ (in projective 
coordinates).
\item Denoting by $\,x\mapsto\overline{x}\,$ any an\-ti\-hol\-o\-mor\-phic 
involution of $\,S^2=\cp^1$ (e.g., the complex conjugation in $\,\bbC\,$ 
extended to the Riemann sphere), we see that, by (iii), the an\-ti-di\-ag\-o\-nal 
$\,2$-sphere $\,\x=\{(\overline{x},x):x\in S^2\}\,$ is \tr\ in 
$\,\y=\,S^2\times\,S^2=\cp^1\!\times\cp^1\nnh$.
\end{enumerate} 

\section{Zooming}\label{zg}
\setcounter{equation}{0}
We now establish what might be called a {\it zooming principle}.
\begin{prop}\label{zoopr}If a compact manifold\/ $\,\x\,$ admits a \tr\ 
im\-mer\-sion/em\-bed\-ding\/ $\,f\,$ in\/ $\,\bbC^n$, then it admits a \tr\ 
im\-mer\-sion/em\-bed\-ding\/ $\,h\,$ in every almost complex manifold\/ 
$\,\y\,$ of complex dimension\/ $\,n$. We may choose\/ $\,h\,$ to be the 
composite of\/ $\,f\,$ with a suitable $\,C^\infty$-\feic\ embedding 
in\/ $\,\y\,$ of an open ball in\/ $\,\bbC^n$ containing\/ $\,f(\x)$.
\end{prop}
\begin{cor}\label{tretn}For every integer\/ $\,n\ge1$, the\/ $\,n$-torus\/ 
$\,T^n$ admits a totally real embedding in every almost complex manifold of 
complex dimension $\,n$. Such an embedding may be chosen so that its image 
lies in any prescribed open subset \feic\ to a ball.
\end{cor}
Proposition~\ref{zoopr} trivially follows from Lemma~\ref{cptgr} below (as 
explained in Remark~\ref{rezoo}), and Corollary~\ref{tretn} then is immediate 
from (v) in Section~\ref{se}. First, we need a definition and another lemma.

Given an immersion $\,\varPhi:\x\to V\,$ of a real $\,n$\diml\ manifold 
$\,\x\,$ in a real or \cvs\ $\,V\nnh$, the {\it Gauss mapping}
\begin{equation}\label{gau}
\text{\rm G}_\varPhi:\x\to\,\text{\rm Gr}_n(V)
\end{equation}
of $\,\varPhi\,$ assigns to each $\,x\in\x\,$ the image 
$\,d\hs\varPhi_x(T_x\x)$. Here $\,\hs\text{\rm Gr}_n(V)\,$ is the Grassmann 
manifold of all $\,n$\diml\ real \vsu s of $\,V\nh$. 
\begin{lem}\label{cpttr}Let\/ $\,J\,$ be a complex structure in a \rvs\ 
$\,V\nnh$, that is, a linear operator\/ $\,V\to V\,$ with 
$\,J^2=\,-\,\text{\rm Id}\hs$, and let\/ $\,\,Y\hs$ be some given set of\/ 
$\,J$-\tr\ subspaces of a fixed dimension\/ $\,n\ge0$. If\/ 
$\,\hs\dim\hs V<\infty$ and\/ $\,Y\hs$ is compact as a subset of\/ 
$\,\hs\text{\rm Gr}_n(V)$, then all\/ $\,W\in Y\,$ are \tr\ relative to every 
complex structure that lies in a suitable \nb\/ $\,\varOmega\,$ of\/ 
$\hs J\,$ in\/ $\,\hs\text{\rm Hom}_\bbR(V,V)$.
\end{lem}
In fact, otherwise there would exist a sequence 
$\,J_\ek\in\,\text{\rm Hom}_\bbR(V,V)\,$ of complex structures and sequences 
$\,W_\ek\in Y\,$ and $\,u_\ek\in V\,$ such that $\,J_\ek\to J$ as 
$\,\ek\to\infty$, while $\,u_\ek\in W_\ek\nnh\cap J_\ek W_\ek$ and 
$\,|u_\ek|=1\,$ for some fixed Euclidean norm $\,|\,\,|\,$ in $\,V$. Using 
compactness of $\,Y$ and the unit sphere, we could pass to subsequences for 
which $\,W_\ek\to W\,$ and $\,u_\ek\to u\,$ with some $\,W\in Y\,$ and 
$\,u\in V\nnh$, so that $\,u\in W\nnh\cap JW$ and $\,|u|=1$, contradicting the 
assumption that $\,\,W\nnh\cap JW=\{0\}\,$ for all $\,W\in Y$. 
\begin{lem}\label{cptgr}Let there be given an \acm\/ 
$\,\y$, a point\/ $\,y\,$ in\/ $\,\y$, a \rmf\/ $\,\x$, a \nb\/ 
$\,\,U'$ of\/ $\,0\,$ in $\,T_y\y$, as well as $\,C^\infty$ mappings 
$\,\varPhi:\x\to T_y\y\,$ and\/ $\,F:U'\to\y\,$ such that\/ $\,\varPhi\,$ is a 
totally real im\-mer\-sion/em\-bed\-ding of\/ $\,\x\,$ in the \cvs\ 
$\,T_y\y$, while $\,F(0)=y\,$ and\/ $\,dF_0$ is the identity mapping of\/ 
$\,T_y\y$. If, in addition, the image\/ $\,Y\,$ of the Gauss mapping 
{\rm(\ref{gau})}, with $\,V=T_y\y$, is compact, and 
$\,\ve:T_y\y\to T_y\y\,$ denotes the multiplication by 
$\,\ve\in\bbR$, then, for some \nb\/ $\,\,U\,$ of\/ $\,0\,$ 
in\/ $\,T_y\y\,$ contained in\/ $\,\,U'\nnh$, and all sufficiently small\/ 
$\,\ve>0$, the composite\/ 
$\,F\nh\circ\ve\circ\varPhi:\varPhi^{-1}(\ve^{-1}U)\to\y\,$ is a totally real 
im\-mer\-sion/em\-bed\-ding in\/ $\,\y\,$ of the open subset\/ 
$\,\varPhi^{-1}(\ve^{-1}U)$ of\/ $\,\x$.
\end{lem} 
\begin{proof}Let $\,\,U'\ni x\mapsto J(x)\in\,\text{\rm Hom}_\bbR(V,V)$, with 
$\,V=T_y\y$, be the $\,F$-pullback to $\,\,U'$ of the original almost 
complex structure in $\,\y$, and let $\,\,U\subset U'$ be a \nb\ of 
$\,0\,$ in $\,V\hs$ such that and $\,J(x)\in\varOmega\,$ for all $\,x\in U$, 
with $\,\varOmega\,$ obtained by applying Lemma~\ref{cpttr} to our $\,Y\hs$ 
and $\,J=J(0)$. Obviously, $\,d\hskip.5pt\ve_v(W)=W\,$ whenever $\,\ve\ne0$, 
$\,v\in V\,$ and $\,W$ is a real vector subspace of $\,V=T_vV=T_{\ve v}V$. 
Hence $\,\ve\circ\hskip.4pt\varPhi:\varPhi^{-1}(\ve^{-1}U)\to U\,$ is a 
totally real \ie\ relative to the almost complex structure in the receiving 
manifold $\,\,U\nh$, pulled back from $\,\y\,$ via $\,F$.
\end{proof}
\begin{rem}\label{rezoo}Lemma~\ref{cptgr} becomes 
particularly simple when $\,\x\,$ is compact (and hence so is $\,Y$), as one 
then has $\,\ve(\varPhi(\x))\subset U$ for sufficiently small $\,\ve>0$, and 
so $\,F\nh\circ\ve\circ\varPhi\,$ is a \tr\ \ie\ of $\,\x\,$ in $\,\y$.
\end{rem}

\section{Immersions of spheres}\label{is}
\setcounter{equation}{0}
The immersed spheres described in Proposition~\ref{sfinw} go back to Whitney 
\cite{whitney-self}. See also \cite{weinstein}.
\begin{lem}\label{multi}Let\/ $\,|\,\,|\,$ be a fixed Euclidean norm in a 
totally real subspace\/ $\,W$ of a \cvs\/ $\,V\hs$ with 
$\,\dimr W=\dimc V\nnh$, and let\/ $\,\,S\,$ denote the sphere of some radius 
$\,\ax>0\,$ in $\,W\nnh$, centered at $\,0$. Furthermore, let\/ 
$\,K\subset\bbC\smallsetminus\{0\}\,$ be an embedded $\,C^\infty$ curve with 
$\,0\notin K+K$, where\/ $\,K+K=\{b+c:b,c\in K\}$. Then the mapping\/ 
$\,f:K\nh\times S\to V\nh$, given by\/ $\,f(c\hh,v)=c\hs v$, is a 
totally real embedding.
\end{lem}
\begin{proof}The differential of $\,f\,$ at $\,(c\hh,v)\,$ sends 
$\,(\dot c\hh,\dot v)\in\,T_c K\times T_vS\,\subset\,\bbC\times W$ to $\,\dot c\hs v+c\hs\dot v$, and so it transforms an $\,\bbR$-basis 
of $\,T_c K\times T_vS$, each of whose vectors 
$\,(\dot c\hh,\dot v)\,$ has either $\,\dot c=0\,$ or $\,\dot v=0$, 
onto a $\,\bbC$-basis of $\,V\nh$. Thus, $\,f\,$ is a \tri\ (see 
Example~\ref{trsub}). Injectivity of $\,f\,$ follows since the relation 
$\,c\hs v=b\hs w\,$ with $\,b,c\in K\,$ and $\,v,w\in S\,$ implies 
$\,v=c^{-1}b\hs w\,$ and so, as $\,W$ is \tr\ and $\,|v|=|w|=\ax$, we have 
$\,b=\pm\hs c\,$ (cf.\ Example~\ref{trsub} again), while $\,b\ne-\hs c\,$ as 
$\,0\notin K+K$.
\end{proof} 
Any fixed Euclidean norm in a \tr\ subspace $\,W$ of a \cvs\ $\,V\hs$ with 
$\,\dimr W=\dimc V=n\,$ gives rise to the group 
$\,\hs\text{\rm SO}(W)\approx\hs\text{\rm SO}(n)\,$ of all 
o\-ri\-en\-ta\-tion-pre\-serv\-ing linear isometries of $\,W\nnh$, with the 
inclusion $\,\hs\text{\rm SO}(W)\subset\text{\rm GL}\hs(V)\,$ obtained by 
extending operators $\,W\to W\,$ com\-plex-lin\-e\-arly to $\,V\nh$.
\begin{prop}\label{sfinw}Let\/ $W\nh$ and\/ $L\hs$ be a \tr\ subspace and a 
complex subspace of a \cvs\/ $V$ with\/ $W\nh\cap L\ne\{0\}$, 
$\hs\dimr W\nnh=\dimc V\nh=n\hs$ and\/ $\dimc L\nh=1$. Next, let\/ 
$\,\varGamma\subset L\,$ be a compact set such that\/ $\,0\in\varGamma\,$ 
and\/ $\,\varGamma\smallsetminus\{0\}\,$ is a\/ $\,1$\diml\/ $\,C^\infty$ 
submanifold of\/ $\,L\,$ with\/ $\,v+v\hh'\nh\ne0\,$ whenever 
$\,v,v\hh'\in\varGamma\smallsetminus\{0\}$, while the intersection of\/ 
$\,\varGamma\,$ with some \nb\ of\/ $\,0\,$ in\/ $\,L\,$ consists of 
two non-par\-al\-lel real-line segments emanating from $\,0$. Finally, let\/ 
$\,\q\,=\,\hs\text{\rm SO}(W)\hskip1pt\varGamma$, so that\/ $\,\q\,$ is the 
set of all\/ $\,Ax\,$ with\/ $\,A\in\,\text{\rm SO}(W)\,$ and\/ 
$\,x\in\varGamma$, where\/ 
$\,\hs\text{\rm SO}(W)\subset\hs\text{\rm GL}\hs(V)\,$ is defined as above for 
any fixed Euclidean norm\/ $\,|\,\,|\,$ in\/ $\,W$.

Then\/ $\,\q\,$ is a totally real\/ $\,n$-sphere immersed in\/ $\,V$. It has 
just one self-in\-ter\-sec\-tion in the form of a double point at\/ $\,0$, and 
its two tangent spaces $\,T,T\pr$ at\/ $\,0\,$ are related by\/ 
$\,T\pr\,=\,z\hskip.8ptT\,$ for some\/ $\,z\in\bbC\smallsetminus\bbR\hs$.
\end{prop}
\begin{proof}Fix $\,u\in W\nh\cap L\smallsetminus\{0\}$. Now 
$\,\varGamma\smallsetminus\{0\}=Ku\,$ for some embedded $\,C^\infty$ curve 
$\,K\subset\bbC\smallsetminus\{0\}\,$ with $\,0\notin K\nh+K$, which 
approaches $\,0\,$ along two real-line segments $\,(0,1)\hs b\,$ and 
$\,(0,1)\hs c$, where $\,b,c\in\bbC\smallsetminus\{0\}\,$ and 
$\,b/c\notin\bbR\hs$. For the sphere $\,S\,$ of radius $\,|u|\,$ in 
$\,W\nnh$, centered at $\,0$, we have $\,\q\smallsetminus\{0\}=KS$, as 
$\,\q\smallsetminus\{0\}=\,\text{\rm SO}(W)Ku
=K\,\text{\rm SO}(W)\hs u$. Thus, in view of Lemma~\ref{multi}, 
$\,\q\smallsetminus\{0\}\,$ is a totally real submanifold of $\,V\nnh$, 
\feic\ to $\,K\times S$, that is, to $\,S^n$ minus two points. However, 
because of how $\,K\,$ approaches $\,0\,$ in $\,\bbC\hs$, a \nb\ of 
$\,0\,$ in $\,\q\,$ is the union of two open $\,n$-balls centered at $\,0$ 
in the \tr\ subspaces $\,T=bW$ and $\,T\pr=cW$ spanned by 
$\,\hs\text{\rm SO}(W)bu$ and $\,\hs\text{\rm SO}(W)c\hh u$. 
\end{proof}
As an obvious consequence of Propositions ~\ref{zoopr} and ~\ref{sfinw}, we 
obtain
\begin{cor}\label{sfinm}Given an \acm\/ $\,\y\,$ with\/ 
$\,\hs\dimc\y=n\ge2$, a point\/ $\,y\in\y$, and a \nb\/ $\,\,U\,$ 
of\/ $\,y\,$ in\/ $\,\y$, there exists a \tri\ of the\/ $\,n$-sphere in\/ 
$\,\,U\,$ which has a double point with a transverse self-in\-ter\-sec\-tion 
at\/ $\,y$, and no other multiple points. 
\end{cor}

\section{Totally real blow-ups}\label{bu}
\setcounter{equation}{0}
Let $\,\y'$ be the complex manifold, \feic\ to 
$\,\y\hs\#\,\overline{\cp^n}\nh$, obtained by blowing up a point in a given 
complex manifold $\,\y\,$ with $\,n=\dimc\y\nh$.

Totally real \ies\ in $\,\y\,$ lead to \tr\ \ies\ in $\,\y'\nnh$. We discuss 
three such  constructions, two of which are summarized in the following lemma; 
the third one is presented, for $\,n=2$ only, in Section~\ref{rb} (see 
Remark~\ref{third}).
\begin{lem}\label{blwup}Given a complex manifold\/ $\,\y$, a point\/ 
$\,y\in\y$, and a $\,\ek$\diml\ totally real submanifold\/ $\,\x\subset\y\,$ 
which is closed as a subset of\/ $\,\y$ and carries the subset topology, 
let the complex manifold\/ $\,\y'=\,\y\,\#\,\,\overline{\cp^n}\nnh$, where\/ 
$\,n=\dimc\y$, be the result of blowing up\/ $\,y\,$ in $\,\y$.
\begin{enumerate}
  \def\theenumi{{\rm\alph{enumi}}}
\item If\/ $\,y\notin\x$, then\/ $\,\x\,$ is also totally 
real as a submanifold of\/ $\,\y'\nnh$.
\item If\/ $\,y\in\x$, then the closure $\,\hat\x$, in\/ $\,\y'\nnh$, of the 
pre\-im\-age of\/ $\,\x\smallsetminus\{y\}$ under the blow-down projection 
$\,\y'\to\y$, is a totally real submanifold of\/ $\,\y'\nnh$, \feic\ to the 
manifold\/ $\,\x\hs'\nh\approx\hs\x\,\,\#\,\,\rp^\ek$ obtained by the real 
blow-up of\/ $\,y\,$ in\/ $\,\x$.
\item If $\,n=\ek=2\,$ and\/ $\,y\in\x$, while $\,\bbP\subset\y'$ is the 
divisor created by the blow-up, then the surface $\,\hat\x\,$ mentioned in\/ 
{\rm(b)} can be deformed, by arbitrarily small isotopies supported in any 
given open subset of\/ $\,\y'$ containing $\,\bbP$, to a surface embedded in 
$\,\y'$ and having a single, transverse intersection with $\,\bbP$.
\end{enumerate} 
\end{lem}
\begin{proof}(a) is obvious. For (b), let us fix a point 
$\,y\hh'\in\pi^{-1}(y)\subset\x\hs'\nh$, where $\,\pi:\x\hs'\to\x\,$ is the 
real blow-down projection. We may choose hol\-o\-mor\-phic coordinates $\,z_a$ 
in $\,\y$, $\,a=1,\dots,n\,$ and $\,C^\infty$ coordinates $\,x_j$ in $\,\x$, 
$\,j=1,\dots,\ek$, both defined near $\,y$, such that $\,z_a=x_j=0\,$ and 
$\,\,\partial z_a/\partial x_j=\delta_{aj}$ at $\,y\,$ for all $\,a,j$, while 
$\,y\hh'\nh$, as a real line through $\,0\,$ 
in $\,T_y\x$, is tangent to the $\,x_1$ coordinate axis. (This is easily 
achieved by an affine coordinate change, cf.\ Example~\ref{trsub}.) In 
suitable local coordinates $\,\ix,\xi_2,\dots,\xi_\ek$ and 
$\,\ez,\zeta_2,\dots,\zeta_n$ for $\,\x\hs'$ and $\,\y'\nh$, the blow-down 
projections $\,\x\hs'\to\x\,$ and $\,\y'\to\y\,$ are given by 
$\,(x_1,\dots,x_\ek)=(\ix,\ix\xi_2,\dots,\ix\xi_\ek)\,$ and 
$\,(z_1,\dots,z_n)=(\ez,\ez\zeta_2,\dots,\ez\zeta_n)$, with $\,\ix\in\bbR\,$ 
and $\,\ez\in\bbC\,$ both varying in \nb s of $\,0$. The integral form of the 
first-order Taylor formula now gives 
$\,z_a=\sum_{j=1}^\ek x_jh_{aj}(x_1,\dots,x_\ek)\,$ for some $\,C^\infty$ 
functions $\,h_{aj}$ with $\,h_{aj}=\,\partial z_a/\partial x_j=\delta_{aj}$ 
at $\,x_1=\ldots=x_\ek=0$. In terms of $\,\ix,\xi_2,\dots,\xi_\ek$ this 
becomes $\,\ez=z_1=\ix\mu(\ix,\xi_2,\dots,\xi_\ek)$ for a $\,C^\infty$ 
function $\,\mu\,$ with $\,\mu(0,\xi_2,\dots,\xi_\ek)=1$. Similarly, 
$\,z_a/\ix\,$ is, for each $\,a=2,\dots,n$, a $\,C^\infty$ function of 
$\,\ix,\xi_2,\dots,\xi_\ek$ (where $\,\ix\,$ varies around $\,0\,$ in 
$\,\bbR$), and hence so is $\,\zeta_a=z_a/\ez=z_a/(\ix\mu)$, while
$\,(\ez,\zeta_2,\dots,\zeta_n)\to(0,\xi_2,\dots,\xi_\ek,0,\dots,0)\,$ as 
$\,\ix\to0$, so that $\,\,\partial\zeta_a/\partial\xi_j=\delta_{aj}$ wherever 
$\,\ix=0$. The inclusion mapping 
$\,\x\smallsetminus\{y\}\to\y\smallsetminus\{y\}\,$ thus has a $\,C^\infty$ 
extension $\,f:\x\hs'\to\y'$ represented by our assignment 
$\,(\ix,\xi_2,\dots,\xi_\ek)\mapsto(\ez,\zeta_2,\dots,\zeta_n)$, 
which is again injective, as it acts by 
$\,(0,\xi_2,\dots,\xi_\ek)\mapsto(0,\xi_2,\dots,\xi_\ek,0,\dots,0)$ on the 
added $\,\rp^{\ek-1}$ represented by $\,\ix=0$. 

At $\,\ix=0\,$ we have $\,\mu=1\,$ and hence 
$\,\,\partial\mu/\partial\xi_j=0$, while, as we just saw, 
$\,\,\partial\zeta_a/\partial\xi_j=\delta_{aj}$ if 
$\,\ix=0$. On the other hand, the relation $\,\ez=\ix\mu\,$ gives 
$\,\,\partial\ez/\partial \ix=\mu=1\,$ and $\,\,\partial\ez/\partial\xi_j=0\,$ 
when $\,\ix=0$. Consequently, at points with $\,\ix=0$, the matrix 
$\,\,\mathfrak{A}\,=[\partial\zeta_a/\partial\xi_j]\,$ with $\,1\le a\le n\,$ 
and $\,1\le j\le\ek$ (where $\,\xi_1=\ix$, $\,\zeta_1=\ez$) has a nonzero 
$\,\ek\times\ek\,$ subdeterminant obtained by restricting 
$\,a\,$ to $\,\{1,\dots,\ek\}$. Therefore, 
$\,\,\text{\rm rank}\,\mathfrak{A}\,=\ek\,$ and so $\,f\,$ is a totally real 
immersion (see Example~\ref{trsub}), which proves (b).

Now let $\,n=\ek=2$, and let $\,f:\x\hs'\to\y'$ be the \tre\ with 
$\,f(\x\hs')=\hat\x$, described above. Thus, $\,\vg=\hat\x\cap\bbP\,$ is the 
$\,f\nh$-im\-age of the circle $\,\rp^1$ in $\,\x\hs'$ created by the real 
blow-up of $\,y\in\x$. A tubular \nb\ $\,\,U\,$ of $\,\vg\,$ in $\,\hat\x\,$ 
is \feic\ to the M\"obius strip, so that we may fix $\,\ve>0\,$ and a two-fold 
covering map $\,F:S^1\times(-\hs\ve,\ve)\to U$ invariant under the 
involution of $\,S^1\times(-\hs\ve,\ve)\,$ sending 
$\,(z,t)\,$ to $\,(-\hs z,-\hs t)$, with $\,S^1\nh=\{z\in\bbC:|z|=1\}$. The 
push-for\-ward under $\,F\,$ of the standard unit vector field on 
$\,S^1\times(-\hs\ve,\ve)\,$ tangent to the $\,(-\hs\ve,\ve)\,$ factor is a 
double-val\-ued vector field $\,\pm\hs v\,$ tangent to $\,\,U\nh$, defined 
only up to a sign. If $\hs\ve\hs$ is made sufficiently small and a Riemannian 
metric is chosen on $\,\y$, the union of suitable short geodesic segments 
emanating from points of $\,\,U\,$ in the direction of $\,\pm\hs iv\,$ is a 
$\,3$\diml\ open submanifold $\,\n\hs$ of $\,\y$, and $\,F$ has an extension 
to a two-fold covering map $\,H:S^1\times D_\ve\to\n\nh$, where 
$\,D_\ve=\{w\in\bbC:|w|<\ve\}$, such that $\,H\,$ is invariant under the 
involution $\,(z,w)\mapsto(-\hs z,-\hs w)$. Making $\,\ve\,$ even smaller, if 
necessary, we may also assume that $\,H(S^1\times D_\ve)\,$ intersects 
$\,\bbP\,$ only along $\,\vg$. (In fact, since $\,\hat\x$ is \tr, at each 
point of $\,\vg\,$ the complex line spanned by $\,\pm\hs v\,$ has a trivial 
intersection with the tangent complex line of $\,\bbP$.) Finally, let us fix 
a $\,C^\infty$ function $\,\varphi:[\hs0,\infty)\to\bbR\,$ equal to $\,1\,$ 
near $\,0\,$ and vanishing outside the interval $\,(-\hs\ve^2/2,\ve^2/2)$. We 
now obtain the required small deformations of $\,\hat\x\,$ by replacing 
$\,\,U\subset\hat\x\,$ with the $\,H$-im\-age of the involution\inv\ surface 
$\,\{(z,t+r\hh i\hs\varphi(t^2)\,\text{\rm Re}\,z):z\in S^1\nh,\hskip4pt
t\in(-\hs\ve,\ve)\}$, depending on a real parameter $\,r\hs$ close to $\,0$, 
and transverse to the curve $\,S^1\times\{0\}$.
\end{proof}
In the following corollary, by the {\it nonorientable \compact\ surface of 
genus\/} $\,s\ge1\,$ we mean, as usual, the connected sum $\,s\hs\rp^2$ of 
$\,\,s\,$ copies of $\,\rp^2$.
\begin{cor}\label{trerp}The surface\/ $\,3\hskip.8pt\rp^2\,
=\,\,\rp^2\#\,\rp^2\#\,\rp^2$ admits totally real embeddings in\/ 
$\,\cp^2\,\#\,\,\mk\hs\overline{\cp^2}$ for all\/ $\,\mk\ge1$.

More generally, for any integers $\,\mk,s\,$ with $\,\mk\ge s-2\ge0$, the 
nonorientable \compact\ surface\/ $\,s\hs\rp^2$ of genus\/ $\,s\hs$ admits a 
totally real embedding in the \csu\ $\,\cp^2\,\#\,\,\mk\hs\overline{\cp^2}$ 
obtained by blowing up the points of a suitable $\,\mk$-el\-e\-ment set in\/ 
$\,\cp^2\nh$. 
\end{cor}
This is clear if one blows up $\,\hs\mk\,$ distinct points of $\,\cp^2\nh$, 
of which $\,s-2\,$ lie in a given \tr\ torus embedded in 
$\,\bbC^2\subset\cp^2$ (cf.\ (v) in Section~\ref{se}), and then uses 
Lemma~\ref{blwup}. 

\section{Removability of complex points by blow-up}\label{rb}
\setcounter{equation}{0}
Given a point $\,y\,$ of a \rsu\ $\,\x\,$ embedded in a \csu\ $\,\y$, we will 
say that $\,\x\,$ contains $\,y\,$ as a {\em complex point removable by 
blow-up\/} if, for some \nb\ $\,\,U\,$ of $\,y\,$ in $\,\x\,$ and some 
\tr\ $\,C^\infty$ submanifold $\,\,U'$ of the \csu\ $\,\y'$ obtained from 
$\,\y\,$ by blowing up the point $\,y$, the blow-down projection 
$\,\pi:\y'\to\y\,$ maps $\,\,U'$ \feicly\ onto $\,\,U$. This amounts 
to requiring that $\,\pi^{-1}(x)\,$ have a limit $\,y'\in\y'$ as 
$\,x\in U\smallsetminus\{y\}\,$ approaches $\,y$, and that 
$\,\,U'=\{y'\}\cup\pi^{-1}(U\smallsetminus\{y\})\,$ be a \tr\ $\,C^\infty$ 
submanifold of $\,\y'$ transverse, at $\,y'\nh$, to the divisor 
$\,\pi^{-1}(y)$. Note that $\,y\,$ then must actually be an (isolated) complex 
point of $\,\x$, since $\,T_y\x\,$ coincides with the image of the 
differential of $\,\pi\,$ at $\,y'$ (which is the complex line in 
$\,T_y\y$ corresponding to $\,y'$), while 
$\,\pi:\y'\smallsetminus\pi^{-1}(y)\to\y\smallsetminus\{y\}\,$ is a 
biholomorphism (and so 
$\,\,U\smallsetminus\{y\}=\pi(U'\smallsetminus\{y'\})\,$ is \tr\ in $\,\y$).
\begin{rem}\label{third}Removability by blow-up leads to our third blow-up 
procedure. Namely, let a \rsu\ $\,\x\,$ embedded in a \csu\ $\,\y\,$ be \tr\ 
except for a finite number of complex points removable by blow-up, and let 
$\,\y'$ be the \csu\ obtained from $\,\y\,$ by blowing up those points. Then 
$\,\y'$ contains a \tr\ embedded surface $\,\x\hs'\nnh$, which the blow-down 
projection $\,\y'\to\y\,$ maps \feicly\ onto $\,\x$.
\end{rem}
\begin{example}\label{rmcgr}Let $\,\x=\{(z,w)\in U\times\bbC:w=\vh(z)\}\,$ 
be the graph of a $\,C^\infty$ function $\,\vh:U\to\bbC\,$ defined on a 
\nb\ $\,\,U\,$ of $\,0\,$ in $\,\bbC\hs$, and let $\,y=(0,\vh(0))$. 
Then the following two conditions are equivalent:
\begin{enumerate} 
  \def\theenumi{{\rm\alph{enumi}}}
\item $\x\,$ treated as a real surface embedded in $\,\y=\bbC^2$ contains 
$\,y\,$ as a complex point removable by blow-up, while 
$\,\x\smallsetminus\{y\}\,$ is \tr;
\item $\vh(z)=\vh(0)+z\varphi(z)\,$ for all $\,z\in U$, where 
$\,\varphi:U\to\bbC\,$ is a $\,C^\infty$ function satisfying at each point the 
inequality $\,\varphi_{\bar z}\ne0\,$ of Example~\ref{trgph}.
\end{enumerate} 
In fact, let $\,\y'$ be the \csu\ obtained from $\,\bbC^2$ by blowing up 
$\,y$. The blow-down projection $\,\pi:\y'\to\bbC^2$ is given by 
$\,(\zeta,\eta)\mapsto(z,w)=(\zeta,\zeta\eta+\vh(0))\,$ in suitable 
hol\-o\-mor\-phic local coordinates $\,(\zeta,\eta)\,$ for $\,\y'\nnh$. The relation 
$\,w=\vh(z)\,$ for $\,z\ne0\,$ now reads $\,\eta=[\vh(\zeta)-\vh(0)]/\zeta\,$ 
whenever $\,\zeta\ne0$. Thus, our claim follows from Example~\ref{trgph}; note 
that a limit $\,L\,$ of $\,[\vh(\zeta)-\vh(0)]/\zeta\,$ as $\,\zeta\to0$, if 
it exists, must be finite: using real values of $\,\zeta$ we get 
$\,L=\vh_x(0)$, with $\,\vh_x=\hs\partial\vh/\partial x\,$ and 
$\,x=\,\text{\rm Re}\,z$.
\end{example}
\begin{example}\label{parab}The graph surface 
$\,\x=\{(z,w)\in U\times\bbC:w=z\overline z\}\subset\bbC^2$ is a paraboloid 
of revolution in the real subspace $\,\bbC\times\bbR\hs$, and $\,y=(0,\nh0)\,$ is 
its unique complex point, as well as a complex point removable by blow-up. 
This is clear from Example~\ref{rmcgr} with $\,\,U=\bbC\,$ and 
$\,\vh(z)=z\overline z$.
\end{example}
\begin{example}\label{prbcp}The embedded $\,2$-sphere $\,\s\subset\cp^2$ given 
by $\,a\overline a=b\overline c\,$ in the homogeneous coordinates 
$\,[a,b,c\hs]\,$ has just two complex points, $\,x=[\hs0,0,1]\,$ and 
$\,y=[\hs0,1,0\hs]$, both removable by blow-up. This is clear from 
Example~\ref{parab}: if one removes $\,x\,$ (or $\,y$) from $\,\s$, by setting 
$\,b=1\,$ (or $\,c=1$), one gets the paraboloid $\,c=|a|^2$ (or $\,b=|a|^2$) 
in the $\,ac\hh$-plane (or, $\,ab\hh$-plane).
\end{example}
\begin{example}\label{rdsph}In $\,\bbC^2$ with the coordinates 
$\,(a,b)$, the $\,2$-sphere $\,\s\,$ of radius $\,R>0$, given by 
$\,|a|^2+|b|^2=R^2$ and $\,\,\text{\rm Im}\hskip1.5ptb=0$, has just 
two complex points $\,x^\pm=(0,\pm R)$, both removable by blow-up. In fact, 
$\,\s\smallsetminus\{x^+\nnh,x^-\}\,$ is \tr\ since $\,\bbC\times\bbR\,$ 
contains only one complex-line direction: that of the $\,z\,$ axis 
$\,\bbC\times\{0\}$. Removability of $\,x^\pm$ is clear from 
Example~\ref{rmcgr}: in the new coordinates $\,(z,w)=(a,R\mp b)$, the 
components of $\,x^\pm$ are $\,(z,w)=(0,\nh0)$, while $\,\s\,$ is, near 
$\,x^\pm\nnh$, the graph of $\,w=\vh(z)\,$ for $\,\vh(z)=R-\sqrt{R^2-|z|^2}$. 
Next, $\,\varphi(z)=\vh(z)/z\,$ equals $\,\overline z\hs F(z\overline z)$, 
where $\,F(s)=[R-(R^2-s)^{1/2}]/s\,$ is obviously real-analytic at $\,s=0$. 
Finally, $\,\varphi_{\bar z}(0)=F(0)=R/2>0$.
\end{example}

\section{Deformations involving complex points}\label{dc}
\setcounter{equation}{0}
The definition of removability by blow-up given in Section~\ref{rb} has an 
immediate extension to the case of immersions. Namely, if $\,f\,$ is an 
immersion of a \rsu\ $\,\x\,$ in a \csu\ $\,\y$, by a {\em complex point of\/} 
$\,f\,$ {\it removable by blow-up\/} we mean any $\,x\in\x\,$ such that 
$\,f(x)\,$ is a complex point, removable by blow-up, for the surface 
$\,f(\x\hs')\subset\y$, where $\,\x\hs'$ is some connected \nb\ of 
$\,x\,$ in $\,\x\,$ with the property that $\,f\,$ restricted to $\,\x\hs'$ is 
injective. As in Section~\ref{rb}, $\,x\,$ then is a {\it complex point\/} of 
$\,f$, that is, $\,\df_x(T_x\x)\,$ forms a complex line in $\,T_{f(x)}\y$. 
Cf.\ Example~\ref{totre}.

The following lemma allows us to modify an immersion $\,f\,$ by moving the 
images of any number $\,\es\,$ of its complex points 
$\,x_j$, removable by blow-up, to arbitrary prescribed locations $\,y_j$. This 
is achieved by replacing the $\,f$-im\-age of a small \nb\ of each 
$\,x_j$ with a thin protrusion or ``tentacle'' reaching all the way to $\,y_j$.
\begin{lem}\label{tntcl}Given $\,\es\,$ distinct points $\,x_1,\dots,x_\es$ in 
a \rsu\ $\,\x$, a \csu\ $\,\y$, points $\,y_1,\dots,y_\es\in\y$, and an 
immersion $\,f:\x\to\y\,$ such that\/ $\,x_1,\dots,x_\es$ are complex points 
of\/ $\,f\,$ removable by blow-up, there exists an immersion $\,f':\x\to\y\,$ 
homotopic to $\,f$, with the same complex points as $\,f$, for which\/ 
$\,x_1,\dots,x_\es$ are complex points removable by blow-up and\/ 
$\,f'(x_j)=y_j$ for $\,j=1,\dots,\es$.

If, in addition, $\,y_1,\dots,y_\es$ are all distinct and lie in 
$\,\y\smallsetminus f(\x)$, the above assertion remains valid also when 
instead of an immersion one speaks, in both instances, of an embedding whose 
image is closed as a subset of $\,\y$.
\end{lem}
\begin{proof}We may assume that $\,\es=1$, since the case $\,\es\ge2\,$ will 
then follow via induction on $\,\es$. Namely, $\,f\,$ can be deformed to 
$\,f'$ in two stages: first, using the inductive assumption that our claim 
holds for $\,x_j,y_j$, $\,j=2,\dots,\es$, with $\,\y\,$ replaced by 
$\,\y\smallsetminus\{y_1\}\,$ if $\,f\,$ is an embedding; then, applying our 
claim for $\,\es=1\,$ to the resulting new immersion and the points 
$\,x_1,y_1$.

With $\,\es=1$, writing $\,x,y\,$ for $\,x_1,y_1$ and $\,p=f(x)$, let us also 
assume that
\begin{enumerate}
\item[($*$)] some biholomorphism $\,(z,w)\,$ maps an open set in $\,\y\,$ 
containing $\,p$ and $\,y\,$ onto a product $\,D\times D\hs'$ of disks 
around $\,0\,$ in $\,\bbC\,$ with $\,1\in D\hs'\nnh$, sending $\,p\,$ to 
$\,(0,\nh0)\,$ and $\,y\,$ to $\,(0,1)$, while, if $\,f\,$ is an embedding and 
$\,f(\x)\,$ is a closed set in $\,\y$, the $\,(z,w)$-pre\-im\-age of 
$\,\{0\}\times(0,1]\,$ does not intersect $\,f(\x)$.
\end{enumerate}
To prove our assertion for $\,\es=1$, under the hypothesis ($*$), we may 
further assume $\,dz\,$ at $\,p$, restricted to $\,\df_x(T_x\x)$, to be 
nonzero. (This is achieved by replacing $\,(z,w)\,$ with 
$\,(z+\mu w(w-1),w)\,$ for $\,\mu\in\bbR\,$ close to $\,0$, and making $\,D\,$ 
smaller if necessary.) Let $\,\,U\,$ now be a \nb\ of $\,x\,$ in $\,\x\,$ such 
that $\,f:U\to\y\,$ is a homeomorphic embedding. Thus, $\,(z,w)$ maps some 
\nb\ of $\,p\,$ in $\,f(U)\,$ onto a graph surface $\,w=\vh(z)\,$ as in 
Example~\ref{rmcgr}, with $\,\vh(0)=0$, and, if $\,D\,$ is made even smaller, 
we may in addition require that $\,\vh\,$ be defined on the whole disk $\,D$, 
and (after the coordinate $\,z\,$ has been replaced by $\,e^{i\theta}z\,$ for 
a suitable $\,\theta\in\bbR$), that also $\,\varphi_{\bar z}(0)\notin\bbR\hs$, 
for $\,\varphi\,$ appearing in Example~\ref{rmcgr}. (Notation of 
Example~\ref{trgph}.) Hence $\,\varphi_{\bar z}(z)\notin\bbR\,$ whenever 
$\,z\in\bbC\,$ and $\,|z|<\ve\,$ for some fixed $\,\ve>0$ that is less than 
the radius of $\,D$. Let us now fix $\,\delta\in(0,\ve)\,$ and a $\,C^\infty$ 
function $\,\alpha:\bbR\to[\hs0,1]\,$ with $\,\alpha=1\,$ on 
$\,(-\infty,\delta\hs]\,$ and $\,\alpha=0\,$ on $\,[\hs\ve,\infty)$. Our 
deformation of $\,f\,$ consists in replacing the graph of $\,\vh:D\to D\hs'$ 
by that of the function $\,\tilde\vh:D\to D\hs'$ with 
$\,\tilde\vh(z)=\vh(z)+\alpha(r)\,$ for $\,z\in D\,$ and $\,r=|z|$. This is 
clear from Example~\ref{rmcgr}, as 
$\,\tilde\vh(z)=\tilde\vh(0)+z\tilde\varphi(z)\,$ with 
$\,\tilde\varphi(z)=\varphi(z)+[\alpha(r)-1]/z$, so that the relations 
$\,2\tilde\varphi_{\bar z}=2\varphi_{\bar z}+r^{-1}d\alpha/dr$ (immediate 
since $\,z_{\bar z}=0\,$ and $\,(r^2)_{\bar z}=z$, that is, 
$\,2r_{\bar z}=z/r$), and $\,\varphi_{\bar z}(z)\notin\bbR$ whenever 
$\,|z|<\ve\,$ (see above) give $\,\tilde\varphi_{\bar z}\ne0\,$ everywhere in 
$\,D$, as required in Example~\ref{rmcgr}(b).

Finally, to obtain our assertion for $\,\es=1$, it suffices to prove it, as we 
just did, under the additional assumption ($*$). In fact, as $\,\y\,$ is 
connected, there exist $\,p_0,p_1,\dots,p_m\in\y\,$ with $\,p_0=f(x)$, 
$\,p_m=y\,$ and such that ($*$) holds if $\,p,y\,$ are replaced by 
$\,p_{l-1},p_l$, for any $\,l=1,\dots,m$. Specifically, $\,f\,$ can be 
modified $\,m\,$ times in a row, leading to successive \tri s 
$\,f_0,f_1,\dots,f_m$ of $\,\x\,$ in $\,\y$, with $\,f_l(x)=p_l$ for 
$\,l=0,\dots,m$, all homotopic to $\,f_0=f$, and we may set $\,f'=f_m$. The 
final clause of the lemma now follows as well, provided that one chooses 
$\,p_1,\dots,p_{m-1}$ and the corresponding $\,m\,$ biholomorphisms in ($*$) 
more carefully, requiring the $\,m\,$ pre\-im\-ages of 
$\,\{0\}\times(0,1]\,$ to be disjoint except for the endpoints shared by 
the $\,l$th and $\,(l+1)$st pre\-im\-age, $\,l=1,\dots,m$, and not to 
intersect $\,f(\x)\,$ except at $\,p=f(x)$. (To guarantee such disjointness 
properties in the case of an embedding, the disk $\,D\,$ mentioned in the last 
paragraph needs to be chosen sufficiently small at each stage.) This completes 
the proof.
\end{proof}
\begin{proof}[Proof of Theorem~\ref{thott}]Assertion (b) is obvious from 
Lemma~\ref{blwup}. To establish (a), let $\,\es\ge2\,$ and let 
$\,y_1,y_2\in\y\,$ be any two of the blown-up points ($y_1\ne y_2$). A 
bi\-hol\-o\-mor\-phic identification of an open set in $\,\y\,$ with a \nb\ of 
$\,(0,\nh0)\,$ in $\,\bbC^2$ allows us to treat a small round sphere $\,\s$ of 
Example~\ref{rdsph} as a \rsu\ in $\,\y\smallsetminus\{y_1,y_2\}\,$ having 
only two complex points $\,x^\pm$, both removable by blow-up. Now (a) follows 
from the final clause of Lemma~\ref{tntcl} applied to $\,\es=1$, the inclusion 
mapping $\,f:\x\to\y\,$ of the sphere $\,\x=\s$, our $\,y_1,y_2$, and 
$\,x_1=x^+\nnh$, $\,x_2=x^-\nnh$.
\end{proof}

\section{Other immersed two-spheres}\label{oi}
\setcounter{equation}{0}
The main result of this section, Corollary~\ref{untre}, establishes the 
existence of \tri s $\,f:S^2\to\,\cp^2\,\#\,\,\mk\hs\overline{\cp^2}\nnh$, for 
$\,\mk\in\{1,3\}$, with some special properties. We use such immersions in 
Sections~\ref{su} and~\ref{rc}.
\begin{lem}\label{holse}Let there be given the total space $\,\y\,$ of a 
hol\-o\-mor\-phic line bundle $\,\mathcal{L}\,$ over a complex curve $\,\x$, a 
hol\-o\-mor\-phic section $\,\phi\,$ of $\,\mathcal{L}$, and a $\,C^\infty$ function 
$\,F:\x\to\bbR\hs$. Then, for the product section $\,F\phi\,$ treated as an 
embedding $\,f\,$ of\/ $\,\x\,$ in the \csu\ $\,\y$, 
\begin{enumerate}
  \def\theenumi{{\rm\roman{enumi}}}
\item the complex points of\/ $\,f\,$ are those $\,x\in\x\,$ at which 
$\,\phi=0\,$ or $\,dF=0$.
\end{enumerate}
A complex point\/ $\,x\,$ of\/ $\,f\,$ is removable by blow-up, cf.\ 
Section~{\rm\ref{rb}}, if
\begin{enumerate}
  \def\theenumi{{\rm\alph{enumi}}}
\item $x\,$ is a simple zero of\/ $\,\phi\,$ and\/ $\,dF\ne0\,$ at\/ 
$\,x$, or
\item $\phi(x)\ne0\,$ and there exists a hol\-o\-mor\-phic coordinate $\,z\,$ 
on a \nb\ of\/ $\,x\,$ in $\,\x\,$ with $\,z=0\,$ at\/ $\,x\,$ and\/ 
$\,F=F(x)+z\varphi\,$ for some $\,C^\infty$ function $\,\varphi\,$ such that\/ 
$\,\varphi(x)=0$.
\end{enumerate}
\end{lem}
\begin{proof}Given $\,x\in\x$, let us choose a hol\-o\-mor\-phic coordinate $\,z\,$ 
on a \nb\ $\,\,U\,$ of $\,x\,$ in $\,\x$, with $\,z=0\,$ at $\,x$, and 
a local hol\-o\-mor\-phic section $\,\psi$ trivializing $\,\mathcal{L}\,$ on 
$\,\,U$. In the resulting hol\-o\-mor\-phic local coordinates $\,(z,w)$ for $\,\y$, 
the image $\,f(\x)\,$ becomes a graph surface $\,w=\vh(z)$, with 
$\,\vh=\sigma F$ for the hol\-o\-mor\-phic function $\,\sigma\,$ such that 
$\,\phi=\sigma\psi\,$ on $\,\,U$. Example~\ref{trgph} now yields (i), since 
$\,\sigma_{\bar z}=0\,$ and so $\,\vh_{\bar z}=\sigma F_{\bar z}$ vanishes 
only where $\,\sigma=0\,$ (that is, $\,\phi=0$) or $\,F_{\bar z}=0\,$ (which, 
as $\,F\,$ is real-val\-ued, amounts to $\,dF=0$). Similarly, 
Example~\ref{rmcgr} gives removability of $\,x\,$ by blow-up, both in case 
(ii) (as $\,\sigma/z\,$ then is hol\-o\-mor\-phic on $\,\,U\,$ and 
$\,F_{\bar z}(x)\ne0$, so that $\,(\sigma F/z)_{\bar z}\ne0$ at $\,x$), and 
(iii) (since we may choose $\,\psi=\phi$, that is, $\,\sigma=1$).
\end{proof}
\begin{prop}\label{trplc}An oriented\/ $\,2$-sphere $\,\x\,$ admits an 
immersion $\,f\hs$ in $\,\cp^2$ such that, for some three-el\-e\-ment set\/ 
$\,Y\subset\x$,
\begin{enumerate}
  \def\theenumi{{\rm\alph{enumi}}}
\item $f_*[\x\hs]\,$ is a generator of $\,H_2(\cp^2\nnh,\bbZ)$,
\item $Y\,$ consists of three complex points of $\,f\hs$ removable by blow-up,
\item the immersion $\,f:\x\smallsetminus Y\to\cp^2$ is \tr.
\end{enumerate}
Furthermore, such $\,f\hs$ and\/ $\,Y\,$ can be chosen so that either
\begin{enumerate}
  \def\theenumi{{\rm\roman{enumi}}}
\item $f\,$ is an embedding, or
\item $Y$ is the $\,f$-pre\-im\-age $\,f^{-1}(y)\,$ of a point\/ 
$\,y\in\cp^2$.
\end{enumerate}
\end{prop}
\begin{proof}To realize (a) -- (c) and (i), we define $\,\y\subset\cp^2$ to be 
the complement of a point in $\,\cp^2\nh$, so that $\,\y\,$ is 
bi\-hol\-o\-mor\-phic to the total space of the dual tautological line bundle 
$\,\mathcal{L}\,$ of a projective line $\,\x\subset\y$. Our claim is now 
obvious from Lemma~\ref{holse} if we choose a hol\-o\-mor\-phic section 
$\,\phi\,$ of $\,\mathcal{L}$ with a unique, simple zero and a $\,C^\infty$ 
function $\,F:\x\to\bbR\,$ having just two critical points, such that either 
critical point $\,x\,$ satisfies condition (b) in Lemma~\ref{holse}. 
Specifically, treating $\,\x\,$ as the Riemann sphere $\,\bbC\cup\{\infty\}$, 
we may set $\,F=z\overline z/(z\overline z+1)\,$ for $\,z\in\bbC\subset\x\nh$. 
The role of the coordinate $\,z$ in Lemma~\ref{holse}(b) is now played by 
$\,z\,$ at $\,x=0\in\bbC\hs$, and by $\,1/z\,$ at $\,x=\infty$, while $\,F\,$ 
becomes a standard height function if we identify $\,\bbC\cup\{\infty\}$, via 
the stereographic projection, with a sphere in $\,\rtr=\bbC\times\bbR\hs$. 

To obtain (a) -- (c) and (ii), we fix $\,f\,$ and $\,Y\,$ with (a) -- (c) and 
(i), and then replace $\,f\,$ by $\,f'$ chosen as in Lemma~\ref{tntcl} 
for $\,\y=\cp^2\nnh$, our $\,f,\x$, and $\,\es=3$, with $\,x_j\in Y\,$ and 
$\,y_j=y$, $\,j=1,2,3$.
\end{proof}
Blowing up all points of $\,Z=f(Y)\,$ in Proposition~\ref{trplc}, we get
\begin{cor}\label{untre}Given $\,\mk\nh\in\nh\{1,3\}$ and any 
$\hs\mk$-el\-e\-ment set\/ $\hs Z\nh\subset\nh\cp^2\nnh$, let\/ $\hs\y$ be the 
\csu\ obtained from $\,\cp^2$ by blowing up all points of\/ $\,Z$. The 
oriented\/ $\,2$-sphere $\,S^2$ then admits a totally real immersion 
$\,f:S^2\nh\to\y$ with\/ $\,f_*[S^2]\,\cdot\,[L]=1\,$ for any projective 
line $\,L\subset\cp^2\smallsetminus Z\,$ treated as a submanifold of\/ $\,\y$, 
where $\,\cdot\hs$ is the intersection form in $\,H_2(\y,\bbZ)$. If, in 
addition, $\,\mk=3$, we may choose\/ $\,f\,$ to be an embedding.
\end{cor}

\section{Removal of transverse intersections}\label{rt}
\setcounter{equation}{0}
In this section we describe how, using infinitely many homotopically 
differrent surgeries, one can remove a transverse self-in\-ter\-sec\-tion of a 
\tr\ surface immersed in an \acsu.

Let $\,\p,\q\,$ be oriented real planes with $\,\p\cap \q=\{0\}\,$ in a 
complex plane $\,V\nh$. We define the {\it sign\/} of the (transverse) 
intersection of $\,\p\,$ and $\,\q\,$ to be {\it positive}, or {\it negative}, 
if the di\-rect-sum orientation of $\,V\nnh=\p\oplus\q\,$ does or, 
respectively, does not, agree with the standard orientation of $\,V\nh$ 
described at the beginning of Section~\ref{tp}. The sign remains the same if 
$\,\p\,$ and $\,\q\,$ are switched, or their orientations are both reversed. 
This allows us to extend the concept of the sign to the case of a double point 
$\,y=f(x)=f(x')$ representing a transverse self-in\-ter\-sec\-tion of an 
immersion $\,f:\x\to\y$, in an \acsu\ $\,\y$, of a \rsu\ $\,\x\,$ which is 
orientable (though not necessarily oriented).
\begin{lem}\label{nozdv}Let\/ $\,\p,\q\,$ and\/ $\,\w$ be real vector 
spaces of dimension\/ $\,n\ge2$, and let\/ $\,\gm:\p\times\q\to\w$ be a 
bi\-lin\-ear mapping such that\/ $\,\gm(u,v)\ne0$ \hbox{whenever 
$\,u\in\p\smallsetminus\{0\}\,$ and\/ $\,v\in\q\smallsetminus\{0\}$. 
If\/ $\,u_j$} and\/ $\,v_\ek$ are fixed bases of\/ $\,\p\,$ and\/ $\,\q$, and\/ 
$\,\mathfrak{T}_\xi$ is the $\,n\times n\,$ matrix with the entries 
$\,\xi(\gm(u_j,v_\ek))$, where $\,\xi\in\w^*$ and\/ $\,j,\ek\in\{1,\dots,n\}$, 
then either\/ $\,\det\,\mathfrak{T}_\xi>0\,$ for all\/ 
$\,\xi\in\w^*\nh\smallsetminus\{0\}$, or\/ $\,\det\,\mathfrak{T}_\xi<0\,$ 
for all\/ $\,\xi\in\w^*\nh\smallsetminus\{0\}$.
\end{lem}
In fact, $\,\det\,\mathfrak{T}_\xi\ne0\,$ whenever 
$\,\xi\in\w^*\nh\smallsetminus\{0\}$, for otherwise there would exist 
$\,v\in\q\smallsetminus\{0\}\,$ such that the injective operator 
$\,\gm(\,\cdot\,,v)\,$ sends $\,\p\,$ into the subspace 
$\,\hs\text{\rm Ker}\,\hs\xi\subset\w$ of dimension $\,n-1$. Our claim now 
follows since $\,\w^*\nh\smallsetminus\{0\}\,$ is connected.
\begin{lem}\label{ctout}Let\/ \tr\ oriented planes\/ $\,\p\hs$ and\/ $\,\q\,$ 
in a complex plane\/ $\,V$ have a negative transverse intersection at $\,0$, 
and let\/ $\,\vg$ be the ca\-non\-i\-cal\-ly-o\-ri\-ent\-ed unit circle\/ 
$\,\{x\in\p:\langle x,x\rangle=1\}\,$ for a fixed Euclidean inner product\/ 
$\,\langle\,,\rangle\,$ in $\,\p$. The oriented cylinder $\,\bbR\times \vg$ 
then admits a \tre\/ $\,\fe\hs$ in\/ $\,V$ such that\/ $\,\fe(t,x)=|t|x\,$ 
if\/ $\,(t,x)\in(-\infty,-\cj\hh]\times\vg$ and\/ $\,\fe(t,x)=t\hh \xa x$ if\/ 
$\,(t,x)\in[\hs\cj,\infty)\times\vg\nnh$, for some $\,\cj\in(0,\infty)\,$ and 
some o\-ri\-en\-ta\-tion-pre\-serv\-ing real-lin\-ear isomorphism 
$\,\xa:\p\to\q$.

Any\/ $\,\fe\hs$ as above and\/ $\,x\in \vg$ give rise to the curve 
$\,[-\cj,\cj\hh]\ni t\mapsto\mathfrak{L}^+(\w_{t,x})$ in the circle 
$\,\hs\text{\rm S}\hs(V^{\wedge2})$, where 
$\,\mathfrak{L}^+$ and \/ $\,\hs\text{\rm S}\hs(V^{\wedge2})\,$ are as in\/ 
{\rm(\ref{ltr})}, while\/ $\,\w_{t,x}$ denotes the oriented plane tangent to 
$\,\fe(\bbR\times \vg)\,$ at $\,\fe(t,x)$. The curves arising in this way from 
all such\/ $\,\fe\,$ represent all fix\-ed-end\-point homotopy classes of 
curves joining\/ $\,\mathfrak{L}^+(\p)\,$ to\/ $\,-\mathfrak{L}^+(\q)\,$ in 
$\,\hs\text{\rm S}\hs(V^{\wedge2})$.
\end{lem}
We now proceed to discuss some consequences of Lemma~\ref{ctout} (a proof of 
which is given in Section~\ref{pl}).

Suppose that a \tr\ immersed surface $\,\x\,$ in an \acsu\ $\,\y\hs$ has only 
one self-in\-ter\-sec\-tion, in the form a transverse double point at some 
$\,y\in\y$, and either $\,\x\,$ is nonorientable, or it is orientable and the 
self-in\-ter\-sec\-tion is negative (as defined above). Then a \tr\ surface 
$\,\x\hs'$ smoothly embedded in $\,\y\,$ and \feic\ to $\,\x\,\#\,K^2\nnh$, 
where $\,K^2$ is the \kb, can be obtained from $\,\x\,$ by replacing the 
disjoint union of disk-like \nb s of $\,y\,$ in the two branches of $\,\x\,$ 
through $\,y$ with a cylinder contained in an open subset of $\,\y\,$ \feic\ 
to $\,\bbR^4\nnh$.

To see this, let us identify a \nb\ of $\,y\,$ in $\,\y\,$ with a \nb\ 
$\,\,U\,$ of $\,0\hs$ in a complex plane $\,V$ in such a way that $\,y=0$, the 
almost complex structure $\,J\,$ of $\,\y\,$ coincides at $\,y=0\,$ with the 
standard complex structure of $\,V\nh$, and $\,\x\cap U=(\p\cup\q)\cap U\,$ 
for some oriented totally real planes $\,\p,\q$ through $\,0\,$ in $\,V$ 
having a negative intersection at $\,0$. If $\,\x\,$ is orientable, we fix an 
orientation of $\,\x\,$ and require, in addition, that the orientations of 
$\,\p$ and $\,\q\,$ agree with those of the two branches of $\,\x\,$ through 
$\,y$. For $\hs\fe$ and $\,\cj\,$ chosen as in Lemma~\ref{ctout}, the image 
$\,Y\hs$ of the Gauss mapping $\,\hs\text{\rm G}_{\nh\fe}$ (see (\ref{gau})) 
is a compact subset of $\,\hs\text{\rm Gr}_2(V)$, as 
$\,\text{\rm G}_{\nh\fe}((-\infty,-\cj\hh]\times\vg)=\{\p\}$ and 
$\,\hs\text{\rm G}_{\nh\fe}([\cj,\infty)\times\vg)=\{\q\}$. Next, we choose an 
open set $\,\,U'\nh\subset U\,$ \feic\ to $\,\bbR^4$ such that 
$\,J(x)\in\,\text{\rm Hom}_\bbR(V,V)\,$ lies in $\,\varOmega\,$ whenever 
$\,x\in U'\nnh$, for $\,\varOmega\,$ obtained by applying Lemma~\ref{cpttr} to 
our $\,Y\hs$ and $\,J=J(0)$. We also select $\,\ve\in(0,\infty)\,$ with 
$\,\ve\fe([-\cj,\cj\hh]\times\vg)\subset U'\nnh$. Our claim is now obvious, 
the embedded cylinder being $\,\ve\fe((-\cj,\cj)\times\vg)$.

In Section~\ref{tn} we give another application of Lemma~\ref{ctout}, which 
also uses the assertion about homotopy classes.

\section{Proof of Lemma~\ref{ctout}}\label{pl}
\setcounter{equation}{0}
Totally real planes through $\,0\,$ in $\,V$ form an open set 
$\,\hs\text{\rm TR}\hh(V)\,$ in the $\,4$\mfd\ $\,\hs\text{\rm Gr}_2(V)$. The 
circle $\,\rp^1=\,\text{\rm U}\hs(1)/\bbZ_2$ acts on 
$\,\hs\text{\rm TR}\hh(V)\,$ freely by $\,(\pm\hs z,\q)\mapsto z\q\,$ for 
$\,z\in\,\text{\rm U}\hs(1)\subset\bbC\,$ and $\,\q\in\,\text{\rm TR}\hh(V)$. 
The quotient $\,\hs\text{\rm TR}\hh(V)/\rp^1$ can be \feicly\ identified, as 
described next, with the $\,3$\hh\mfd\ $\,N\hs$ of all plane circles in the 
unit sphere $\,S^2\subset\rtr\nnh$, that is, circles obtained by intersecting 
$\,S^2$ with planes which need not pass through $\,0$. The corresponding 
principal $\,\rp^1$-bun\-dle projection $\,\pi:\text{\rm TR}\hh(V)\to N$ sends 
$\,\q\in\,\text{\rm TR}\hh(V)\,$ to the image 
$\,\pi(\q)\,$ of $\,\,\q\smallsetminus\{0\}\,$ under the standard 
projection $\,V\smallsetminus\{0\}\to\,\text{\rm P}\hs(V)\approx\hs\cp^1\nh
\approx S^2\nnh$.

In fact, each such image $\,\pi(\q)\,$ is a circle in the 
Riemann sphere $\,S^2\approx\bbC\cup\{\infty\}$, as one sees choosing an 
isomorphic identification $\,V\nnh=\bbC^2$ under which $\,\q=\bbR^2$ (and 
so $\,\pi(\q)=\bbR\cup\{\infty\}\subset\bbC\cup\{\infty\}$). Furthermore, if 
$\,\q,\q\hh'\nh\in\,\text{\rm TR}\hh(V)\,$ and $\,\pi(\q)=\pi(\q\hh'\hh)$, 
then $\,zu,av,c\hs(\nh u+v)\in\q\hh'$ for any fixed basis $\,u,v\,$ of 
$\,\q\,$ and some $\,z,a,c\in\bbC\smallsetminus\{0\}$, so that, writing 
$\,c\hs(\nh u+v)\,$ as a real combination of $\,zu\,$ and $\,av\,$ (which span 
$\,\q\hh'$ over $\,\bbR$), we get $\,a/z\in\bbR\hs$, and hence 
$\,\q\hh'\nh=z\q$. Finally, $\,\pi:\text{\rm TR}\hh(V)\to N\hs$ has constant 
rank due to its e\-qui\-var\-i\-ance relative to the transitive actions of the 
linear group $\,\hs\text{\rm GL}\hs(V)\approx\hs\text{\rm GL}\hs(2,\bbC)\,$ on 
$\,\hs\text{\rm TR}\hh(V)\,$ and of the complex automorphism group 
$\,\hs\text{\rm Aut}\hh(S^2)\,$ on $\,N\nh$, for the standard homomorphism 
$\,\hs\text{\rm GL}\hs(V)\to\hs\text{\rm Aut}\hh(S^2)
=\hs\text{\rm PGL}\hs(V)$.

Note that $\,N\nh=\widetilde N\nnh/\bbZ_2$, where $\,\widetilde N\hs$ is the 
manifold of all oriented plane circles in $\,S^2$ (\feic\ to 
$\,S^2\nh\times(0,\pi)\,$ via the cen\-ter-ra\-di\-us 
pa\-ram\-e\-tri\-za\-tion of oriented circles in the oriented $\,2$-sphere), 
and the fix\-ed-point free involution generating the $\,\bbZ_2$ action 
reverses the circle orientations. 

We begin by proving our assertion in the special case where the circles 
$\,\pi(\p),\pi(\q)\subset S^2$ are disjoint. We then choose a smoothly 
embedded curve segment joining $\,\pi(\p)\,$ to $\,\pi(\q)\,$ in $\,N\nh$, 
consisting of pairwise disjoint circles, and lift it: first to the bundle 
space $\,\hs\text{\rm TR}\hh(V)\,$ over $\,N\nh$, and from there to the bundle 
over $\,\hs\text{\rm TR}\hh(V)\,$ formed by the manifold of all $\,\bbC$-bases 
of $\,V\nh$, which we identify with $\,\hs\text{\rm GL}\hs(V)\,$ using a fixed 
basis of $\,\p$. This yields a $\,C^\infty$ embedding 
$\,[-\cj,\cj\hh]\ni t\mapsto\xb_t\in\text{\rm GL}\hs(V)$, with any fixed 
$\,\cj\in(0,\infty)$, such that $\,\xb_{-\cj}=\,\text{\rm Id}\hs$, 
$\,\xb_\cj\p=\q$, and the circles $\,\pi(\xb_t\p)\,$ are pairwise disjoint for 
all $\,t\in[-\cj,\cj\hh]$.

Next, we replace the curve $\,t\mapsto\xb_t$ by a composite in which the 
original curve is preceded by a $\,C^\infty$ function 
$\,\bbR\to[-\cj,\cj\hh]\,$ equal to $\,-\cj\,$ on $\,(-\infty,-\cj\hh]$, to 
$\,\cj\,$ on $\,[\hs\cj,\infty)$, and having a positive derivative everywhere 
in $\,(-\cj,\cj)$. The resulting new curve, still written as $\,t\mapsto\xb_t$ 
(but now defined on $\,\bbR$), leads to the required \tre\ 
$\,\fe:\bbR\times\vg\to V\nh$, given by 
$\,\fe(t,x)=e^{\hh\varphi(t)}\nh\xb_tx\,$ for a $\,C^\infty$ immersion 
$\,\varphi:\bbR\to\bbC\,$ with $\,\varphi(t)=\log|t|$ whenever $\,|t|\ge\cj$. 
The choice of $\,\varphi\,$ on $\,(-\cj,\cj)\,$ is described below.

We denote by $\,\ha\,$ the complex-lin\-e\-ar operator $\,V\to V$ whose 
restriction to the oriented Euclidean plane $\,\p\hs$ is the positive rotation 
by the angle $\,\pi/2$. If $\,\fe\,$ corresponds as above to any given 
$\,C^\infty$ immersion $\,\varphi:\bbR\to\bbC\hs$, then
\begin{equation}\label{dft}
\df_{(t,x)}(1,0)\,=\,e^\varphi(\dot\varphi\xb\hs+\dot\xb)x\hskip8pt\text{\rm 
and}\hskip9pt\df_{(t,x)}(0,\ha x)\,=\,e^\varphi\nh\xb\ha x
\end{equation}
for the basis $\,(1,0),(0,\ha x)\,$ of $\,T_{(t,x)}(\bbR\times \vg)$. (We omit 
$\,t\,$ in our notation, and set $\,(\,\,)\dot{\,}=\,d/dt$.) Thus, $\,\fe\,$ 
is a \tri\ if and only if, for every $\,(t,x)\in\bbR\times \vg\nnh$, the 
($e^\varphi\nh\xb$-im\-ages of the) vectors 
$\,(\dot\varphi+\xb^{-1}\nh\dot\xb)x$ and $\,\ha x\,$ are linearly independent 
over $\,\bbC\hs$. Since $\,\p$ is \tr, we may identify $\,V^{\wedge2}$ with 
$\,\bbC\,$ in such a way that $\,x\wedge \ha x=1\,$ for some $\,x\in \vg\nnh$, 
and so the linear independence condition will follow if 
$\,\dot\varphi\ne-\hs s\xb^{-1}\dot\xb x\wedge \ha x\,$ for all 
$\,(s,t,x)\in[\hs0,1]\times[-\cj,\cj\hh]\times\vg\nnh$. (Note that 
$\,x\wedge \ha x=1\,$ for {\it every\/} $\,x\in \vg\nnh$, as each pair 
$\,x,\ha x\,$ is a positive orthonormal basis of $\,\p$.) The condition holds 
if $\,|t|\ge\cj$, since we then have $\,\dot\varphi\ne0=\dot\xb$.

Such an immersion $\,\varphi$, with $\,\varphi(t)=\log|t|\,$ whenever 
$\,|t|\ge\cj$, exists and may be chosen so that 
$\,\dot\varphi:[-\cj,\cj\hh]\to\bbC\smallsetminus\{0\}$ belongs to any 
prescribed fixed-end\-point homotopy class of curves joining $\,-1/\cj\,$ to 
$\,1/\cj\,$ in $\,\bbC\smallsetminus\{0\}$. In fact, as $\,\dot\xb=0\,$ 
whenever $\,|t|\ge\cj$, there exist constants $\,\kappa>0\,$ and 
$\,\delta\in(0,\cj)$ with $\,|\xb^{-1}\dot\xb x\wedge \ha x\hs|\le\kappa\,$ 
for all $\,(t,x)\in[-\cj,\cj\hh]\times \vg\nnh$, and 
$\hs|\xb^{-1}\dot\xb x\wedge\ha x\hs|<1/\cj$ if, in addition, 
$\,|t|\ge\delta$. 
We now select $\,\varphi\,$ with $\,\varphi(-t)=\varphi(t)\,$ and 
$\,|\dot\varphi(t)|>1/\cj$ whenever $\,|t|\in[\hs\delta,\cj)$, such that one 
also has $\,|\dot\varphi|=\mx\,$ on $\,[-\hs\delta,\delta\hh]\,$ for a 
constant $\,\mx>\kappa$. As a result, $\,\fe\,$ is a \tri, and it is also 
injective due to pairwise disjointness of the circles $\,\pi(\xb_t\p)\,$ for 
$\,t\in[-\cj,\cj\hh]$. A curve $\,\varphi:[-\hs\delta,\delta\hh]\to\bbC\,$ 
with the constant speed $\,\mx\,$ and a fixed length can clearly be chosen so 
that, before returning at time $\,t=\delta\,$ to the initial point 
$\,\varphi(-\hs\delta)\,$ with the velocity 
$\,-\hs\dot\varphi\hh(\nh-\hs\delta)$, it traverses any prescribed number 
$\,q\,$ of times, clockwise or counterclockwise, some circle of radius 
depending on $\,q$. This produces the required homotopy effect: using the 
definition of $\,\mathfrak{L}^+$ (see (\ref{ltr})), our identification 
$\,V^{\wedge2}\nh=\bbC\,$ with $\,x\wedge \ha x=1\,$ for all $\,x\in \vg\nnh$, 
and (\ref{dft}), we see that the homotopy classes in question correspond to 
those of the curves 
$\,t\mapsto 
e^\varphi(\dot\varphi\xb\hs+s\dot\xb)x\,\wedge\,e^\varphi\nh\xb\ha x\,$ with 
$\,s=1$, joining $\,-\cj\,$ to $\,\cj\,$ in $\,\bbC\smallsetminus\{0\}$. However, 
we may replace $\,s=1\,$ by $\,s=0$, since our choice of $\,\varphi\,$ turns 
$\,s\in[\hs0,1]\,$ into a homotopy parameter. The resulting homotopy class of 
$\,t\mapsto e^{2\varphi}\dot\varphi\hs\det\nnh\xb\,$ equals that of 
$\,t\mapsto\dot\varphi\hs\det\nnh\xb\,$ (we deform the factor $\,e^{2\varphi}$ 
to a constant in $\,\bbC\smallsetminus\{0\}\,$ by first deforming 
$\,\varphi\,$ in $\,\bbC$). As the curve $\,t\mapsto\xb_t$ does not depend on 
$\,\varphi$, our claim about homotopy classes follows.

Also, $\,\xb_\cj:\p\to\q\,$ is an {\it o\-ri\-en\-ta\-tion-pre\-serv\-ing\/} 
isomorphism. Namely. using the original orientation of $\,\p$, we now orient 
the planes $\,\q_t=\xb_t\p$, $\,t\in(-\cj,\cj\hh]$, by requiring that 
$\,\xb_t:\p\to\q_t$ be o\-ri\-en\-ta\-tion-pre\-serv\-ing. This might, 
possibly, result in re-orient\-ing $\,\q=\q_\cj$, and we will prove that is 
does not, by showing that at least one new\-ly-orient\-ed $\,\q_t$, 
$\,t\in(-\cj,\cj\hh]$, has a negative intersection with $\,\p\,$ at $\,0$. 
(The words `at least one' then may be replaced with `every' by reasons of 
continuity.) Let us fix $\,t\in(-\cj,\cj\hh]$. The circles $\,\pi(\p)\,$ and 
$\,\pi(\q_t)\,$ are disjoint, so that $\,\gm:\p\times\q_t\to V^{\wedge2}$ 
given by $\,\gm(u,v)=u\wedge v\,$ satisfies the hypotheses of 
Lemma~\ref{nozdv} for $\,n=2$. Identifying $\,V^{\wedge2}$ with $\,\bbC\,$ as 
before, we may apply the conclusion of Lemma~\ref{nozdv} to the bases 
$\,u_j,v_\ek$ with $\,u_1=x,\hs u_2=\ha x\,$ (where $\,x\in\vg$), and 
$\,v_j=\xb_tu_j$ for $\,j=1,2$, along with $\,\xi=\,\text{\rm Re}\hs\,$ and 
$\,\xi=\,\text{\rm Im}\hs$. Taking the limit as $\,t\to\cj$, we have 
$\,\xb_t\to\,\text{\rm Id}\hs$, and so $\,\det\,\mathfrak{T}_\xi>0\,$ for 
$\,\xi=\,\text{\rm Im}\hs\,$ and $\,t>\cj\,$ close to $\,\cj$, since 
$\,\det\,\mathfrak{T}_\xi\to1\,$ as $\,t\to\cj\,$ when 
$\,\xi=\,\text{\rm Re}\hs$. (In fact, the limit, as $\,t\to\cj$, of the 
complex $\,2\times2\,$ matrix with the entries $\,\gm(u_j,v_\ek)\,$ has the rows 
$\,[\hs 0\hskip8pt1\hs]$ and $\,[\hs-1\hskip8pt0\hs]$.) Writing 
$\,[v_1\hskip6ptv_2]=[u_1\hskip6ptu_2]\hskip1pt\mathfrak{B}$, where 
$\,\mathfrak{B}\,$ is a complex $\,2\times2$ transition matrix depending on 
$\,t$, we get the transition formula
\[
[u_1\hskip6ptu_2\hskip8ptv_1\hskip6ptv_2]\,
=\,[u_1\hskip6ptu_2\hskip8ptiu_1\hskip6ptiu_2]\,
\left[\begin{matrix}\text{\rm Id}&\text{\rm Re}\,\mathfrak{B}\,\cr
0&\text{\rm Im}\,\mathfrak{B}\end{matrix}\right]
\]
with a real $\,4\times4\,$ matrix composed of four $\,2\times2\,$ blocks. As 
one easily verifies, 
$\,\det\,\mathfrak{T}_{\text{\rm Im}}=\det\hs\text{\rm Im}\hs\mathfrak{B}$, so 
that $\,\det\hs\text{\rm Im}\hs\mathfrak{B}>0$. Negativity of the intersection 
follows: the basis $\,u_1,u_2,iu_1,iu_2$ represents the opposite of the 
standard orientation of $\,V\nh$, cf.\ Section~\ref{tp}.

Using $\,\varphi\,$ as above, the corresponding $\,\fe$, and $\,\xa=\xb_\cj$, 
we thus get our claim under the additional assumption that 
$\,\pi(\p)\cap\pi(\q)=\hs\emp\hh$. The general assertion, for an arbitrary 
pair $\,\p,\q$, can now be derived from a {\it reduction argument}. Namely, to 
prove it for the given $\,\p,\q$, we just need to establish it for the pair 
$\,R,\q$ with some oriented Euclidean plane $\,R\hs\,$ satisfying the same 
hypotheses as $\,\p\,$ and such that there exist a \tre\ $\,h:\p\to V\nnh$, 
an o\-ri\-en\-ta\-tion-pre\-serv\-ing linear isometry $\,\varPsi:\p\to R$, 
and constants $\,c,\cj\hh'\nnh$, for which $\,\cj\hh'\nh>c>\cj\,$ and 
$\,h(\p)\cap\q=\{0\}$, as well as
\begin{enumerate}
  \def\theenumi{{\rm\alph{enumi}}}
\item $h(y)=\varPsi(y)\,$ whenever $\,y\in\p\,$ and $\,|y|\le c$,
\item $h(y)=y\,$ whenever $\,y\in\p\,$ and $\,|y|\ge\cj\hh'\nnh$,
\end{enumerate}
where $\,\cj$, along with some $\,\xa\,$ and $\,\fe:\bbR\times\vg\to V\nh$, 
satisfies the assertion of the lemma for $\,R,\q\,$ instead of $\,\p,\q$. In 
fact, given such $\,R,h,\varPsi,c,\cj\hh'$, we define the analogues of 
$\,\cj,\xa,\fe\,$ for $\,\p\,$ and $\,\q\,$ to be $\,\cj\hh'\nnh,\xa\varPsi$, 
and the mapping $\,\fe'$ given by $\,\fe'(t,x)=h(|t|x)\,$ if $\,t\le-\cj\,$ 
and $\,\fe'(t,x)=\ve\fe(t/\ve,\varPsi x)\,$ if $\,t\ge-\hs c$, for suitable 
$\,\ve>0$. Both formulae for $\,\fe'$ clearly define $\,V\nh$-val\-ued \tre s 
of appropriate domains. They also agree, yielding $\,|t|\hs\varPsi x$, when 
$\,-\hs c\hs\le t\le\nh-\cj\,$ (and $\,x\,$ lies in the unit circle 
$\,\vg\hh'\nh=\varPsi^{-1}(\vg)\subset\p$); thus, the only property of 
$\,\fe'$ that we need to verify is its injectivity, which amounts to 
$\,\fe'((-\infty,-\cj\hh]\times\vg\hh')\cap\fe'((-\cj,\infty)\times\vg\hh')$. 
However, if $\,\ve\in(0,1)$, the definition of $\,\fe'$ implies that 
$\,\fe'((-\cj,\infty)\times\vg\hh')\,$ is contained in the union of three 
sets: the radius $\,\cj\,$ open disk $\,D\,$ about $\,0\,$ in $\,R$, the 
compact set $\,\ve S$, where $\,S=\fe([-\cj,\cj\hh]\times\vg)$, and $\,\q$. 
None of the three sets intersects 
$\,Y\nh=\fe'((-\infty,-\cj\hh]\times\vg\hh')\,$ if $\,\ve\,$ is chosen so that 
$\,\ve S\subset U\,$ for a \nb\ $\,\,U\,$ of $\,0\,$ in $\,V$ with 
$\,\,U\cap Y\nh=\hs\emp\hh$. In fact, $\,D\cap Y\nh=\hs\emp\,$ since $\,D\,$ 
and $\,Y\hs$ are the $\,h$-im\-ages of $\,\varPsi^{-1}(D)\,$ and 
$\,\p\smallsetminus\varPsi^{-1}(D)$, while $\,\q\cap Y\nh=\hs\emp\,$ as 
$\,Y\nh\subset h(\p\smallsetminus\{0\})$ and $\,h(\p)\cap\q=\{0\}$. Finally, 
such $\,\fe'$ realize all the homotopy classes required in the lemma, since we 
assume here that the same can be achieved by varying the choices of $\,\fe\,$ 
for the pair $\,R,\q$.

We apply the above reduction argument twice. In both cases, using arbitrary 
$\,c,\cj\hh'$ with $\,\cj\hh'\nh>c>\cj$, we obtain $\,R,h,\varPsi\,$ by 
exhibiting a $\,3$\diml\ real vector subspace $\,\w\subset V$ containing both 
$\,\p\,$ and a plane $\,R\,$ with $\,\p\cap R\cap L=R\cap\w\nh\cap\q=\{0\}$, 
where $\,L\,$ is (always) the unique complex line through $\,0\,$ in $\,V$ 
contained in $\,\w\nh$. Note that $\,\dim\hs(\w\nh\cap\q)=1$, as $\,\q\,$ 
intersects $\,\p\,$ trivially and so cannot be contained in $\,\w\nh$. Also, 
if $\,\hs\pro:\w\to\py$ denotes the quotient projection onto the plane 
$\,\py=\w/(\p\cap R)$, it is clear that $\,\pro(L)=\py$, while 
$\,\ly_1=\,\pro(\p)$, $\,\ly_2=\,\pro(R)\,$ and 
$\,\ly_3=\,\pro(\w\nh\cap\q)\,$ are three distinct lines in $\,\py$. We may 
thus choose a vector $\,v\in L\smallsetminus\{0\}\,$ such that 
$\,\hs\pro\,v\notin\ly_1\cup\ly_2\cup\ly_3$, while, for 
$\,\ly_4=\,\pro(\bbR v)\,$ and some Euclidean norm in $\,\py$, an arc of the 
unit circle around $\,0\,$ in $\,\py\,$ intersects each of the four lines 
$\,\ly_j$ just once, in this order: $\,\ly_1,\dots,\ly_4$. Let us refer to the 
direction of $\,v\,$ in $\,\w$ as \hbox{{\it vertical\hh}.} For every 
$\,C^\infty$ function $\,\varphi:\p\to\bbR\hs$, the embedding $\,h:\p\to V$ 
given by $\,h(y)=y+\varphi(y)\hs v\,$ is \tr, since so the surface $\,h(\p)$. 
(As $\,L$ contains the vertical vector $\,v$, it cannot be a tangent plane of 
a graph.) To get (a), (b) and $\,h(\p)\cap\q=\{0\}$, we now set 
$\,\varphi=\varphi_1\varphi_2$, where $\,\varphi_1$ is a $\,C^\infty$ function 
$\,\p\to[\hs0,1]\,$ with 
$\,\varphi_1(y)=1\,$ whenever $\,|y|\le c\,$ and $\,\varphi_1(y)=0\,$ whenever 
$\,|y|\ge\cj\hh'\nnh$, while $\,\varphi_2$ is a linear functional on $\,\p\,$ 
having the graph $\,R$. (Such $\,\varphi_2$ exists since the vertical vector 
$\,v\,$ does not lie in $\,R$.) The condition $\,h(\p)\cap\q=\{0\}\,$ follows 
as $\,\hs\pro(h(\p))\cap\,\pro(\w\nnh\cap\q)=\{0\}$ in $\,\py\,$ and 
$\,\hs\pro:\w\nnh\cap\q\to\py\,$ is injective, while the orientation and inner 
product in $\,R\,$ are chosen so as to make $\,\varPsi:\p\to R$, with 
$\,\varPsi y=y+\varphi_2(y)$, an o\-ri\-en\-ta\-tion-pre\-serv\-ing linear 
isometry, and the intersection of $\,R\,$ and $\,\q$ then is negative due to 
our choice of $\,v$.

Using the fact that our claim holds when $\,\pi(\p)\cap\pi(\q)=\hs\emp\hs$, we 
now establish it under the weaker additional assumption $\,\pi(\p)\ne\pi(\q)$. 
Specifically, $\,\w$ and $\,R\,$ required above are selected by fixing 
$\,R\in\text{\rm TR}\hh(V)\,$ for which $\,\pi(R)\cap\pi(\q)=\hs\emp\hs\,$ and 
$\,\pi(\p)\cap\hs\pi(R)\,$ has two elements. We may also assume that 
$\,\dim\hh(\p\cap\nh R)=1\,$ (which is achieved by replacing $\,R$ with 
$\,zR$ for a suitable complex number $\,z\ne0$), and then we set 
$\,\w\nh=\hs\spanr(\p\cup\nnh R)$. Thus, $\,\p\cap R\cap L=\{0\}\,$ for 
$\,L\subset \w$ as above: one element of $\,\pi(\p)\cap\hs\pi(R)\,$ is 
$\,\hs\spanc(\p\cap\nh R)$, and so $\,L\,$ must be its other element, spanned 
over $\,\bbR\,$ by a line in $\,\p\,$ and a line in $\,R\,$ different from 
$\,\p\cap\nh R$. In view of the reduction argument, our assertion, valid for 
$\,R\,$ and $\,\q\,$ (since $\,\pi(R)\cap\pi(\q)=\hs\emp$), hold for 
$\,\p,\q\,$ as well.

Similarly, using our conclusion for arbitrary $\,\p,\q\,$ with 
$\,\pi(\p)\ne\pi(\q)$, we now prove it in the remaining case 
$\,\pi(\p)=\pi(\q)$. Namely, we choose $\,\w$ and $\,R\,$ to be any real 
$\,3$-space in $\,V$ containing $\,\p\,$ and, respectively, any plane in 
$\,\w$ with $\,\p\cap R\cap L=\{0\}$, for the complex line $\,L\subset \w\nh$. 
Thus, $\,\pi(R)\ne\pi(\p)=\pi(\q)$, so that our claim is 
true for $\,R,\q\,$ (see the last paragraph), and, consequently, also for 
$\,\p,\q$.

\section{Connected sums}\label{cs}
\setcounter{equation}{0}
The main result of this section is Theorem~\ref{cnsum}, which, essentially, 
establishes the existence of a con\-nect\-ed-sum operation in the category of 
totally real surfaces immersed in a fixed almost complex surface $\,\y$. The 
totally real immersed cylinder joining the two surfaces is of the type known 
as {\it Whitney's umbrella\/}; it has a unique self-in\-ter\-sec\-tion (which 
is transverse and negative, in the sense of Section~\ref{rt}). Since there 
always exist arbitrarily many pairwise disjoint totally real tori embedded in 
$\,\y\,$ (Corollary~\ref{tretn}), a self-in\-ter\-sec\-tion cannot, in 
general, be avoided: otherwise, connected sums of such tori would constitute 
\tr\ \compact\ orientable surfaces of all genera $\,s\ge2$, embedded in any 
given $\,\y$, in patent contradiction to Corollary~\ref{codie}.
\begin{thm}\label{cnsum}Let\/ $\,\x^\pm$ be two disjoint totally real 
surfaces embedded in an almost complex surface $\,\y$, both endowed with the 
subset topology and closed as subsets of\/ $\,\y$. Also, let\/ 
$\,x^\pm\in\x^\pm\nnh$. Then, for any sufficiently small closed\/ 
$\,2$-disks\/ $\,D^\pm\nh$ \feicly\ embedded in\/ $\,\x^\pm\nh$ and 
containing\/ $\,x^\pm\nh$ as interior points, there exist a totally real 
surface $\,\x\,$ immersed in $\,\y\,$ and an immersion 
$\,h:[-1,1]\times S^1\to\y\,$ such that
\begin{enumerate}
  \def\theenumi{{\rm\roman{enumi}}}
\item $\x\,$ is the disjoint union of\/ $\,\x^+\nh\smallsetminus D^+\nnh$, 
$\,\x^-\nh\smallsetminus D^-$ and the image of $\,h$,
\item $h\,$ sends the circles $\,\{\pm1\}\times S^1$ onto 
$\,\partial D^\pm\nnh$,
\item $\x\,$ has just one self-in\-ter\-sec\-tion\/{\rm:} a double point of\/ 
$\,h\,$ in $\,(-1,1)\times S^1\nnh$,
\item the self-in\-ter\-sec\-tion of\/ $\,\x\,$ is transverse and negative, 
cf.\ Section~{\rm\ref{rt}},
\item the union of $\,D^+\nh\cup\hs D^-\nh$ and the image of\/ $\,h\,$ is 
contained in an open subset of\/ $\,\y\hs$ \feic\ to $\,\bbR^4\nnh$.
\end{enumerate} 
\end{thm}
We prove Theorem~\ref{cnsum} in Section~\ref{pc}. First, we need three lemmas, 
in which $\,(\,\,)\dot{\,}=\,d/dt$, $\,(\,\,)'=\,d/ds$, 
$\,(\,\,)_t=\,\partial/\partial t\,$ and $\,(\,\,)_s=\,\partial/\partial s\,$ 
stand for ordinary or partial derivatives with respect to real variables\/ 
$\,\hs t\,$ and $\,s$.
\begin{lem}\label{zaxes}If\/ $\,c,\ve\in(0,\infty)$, while 
$\,t\mapsto\varphi(t)\,$ and\/ $\,s\mapsto\psi(s)\,$ are $\,C^\infty$ 
functions of\/ $\,t\in(-\hs c,c)\,$ and\/ $\,s\in(-\hs\ve,\ve)$, valued in a 
fi\-nite\diml\ \vs\/ $\,V\nnh$, then the following two conditions are 
equivalent\/{\rm:}
\begin{enumerate}
  \def\theenumi{{\rm\alph{enumi}}}
\item $\varphi(t)=F_s(t,0)\,$ and $\,\psi(s)=F_t(0,s)\,$ for all\/ 
$\,(t,s)\in(-\hs c,c)\times(-\hs\ve,\ve)$ and some $\,C^\infty$ function 
$\,F:(-\hs c,c)\times(-\hs\ve,\ve)\to V$ having the property that\/ 
$\,F(t,0)=F(0,s)=0\,$ whenever $\,t\in(-\hs c,c)\,$ and\/ 
$\,s\in(-\hs\ve,\ve)$,
\item $\varphi(0)=\psi(0)=0\,$ and\/ $\,\dot\varphi(0)=\psi\hh'(0)$.
\end{enumerate}
\end{lem}
In fact, (a) implies (b): 
$\,\dot\varphi(0)=F_{st}(0,\nh0)=F_{ts}(0,\nh0)=\psi\hh'(0)$. Conversely, 
assuming (b), we may set 
$\,F(t,s)=ts\hs[\hh\alpha(t)+\beta(s)-\dot\varphi(0)\hh]\,$ for $\,C^\infty$ 
functions $\,\alpha,\beta\,$ with $\,\varphi(t)=t\hs\alpha(t)\,$ and 
$\,\psi(s)=s\hh\beta(s)$ (so that $\,\alpha(0)=\beta(0)=\dot\varphi(0)$).
\begin{lem}\label{crvcp}Given an \acsu\ $\,\y$, a point\/ $\,y\in\y$, and a 
nonzero vector $\,u\in T_y\y$, there exists an $\,\y\nnh$-val\-ued\/ 
$\,C^\infty$ embedding\/ $\,f\nnh\,$ of a \nb\ of\/ $\,\hs0\,$ in\/ 
$\,\bbC\,$ such that, for all real\/ $\,t,s\nh\,$ near $\,0$, the 
differentials $\,\df_t$ and\/ $\,\df_{is}$ are com\-plex-lin\-e\-ar, while 
$\,f(0)=y\,$ and\/ $\,\df_0$ sends the vectors $\,1,i\in\bbC=T_0\bbC\,$ to\/ 
$\,\hs u\nh\,$ and $\,Ju$.
\end{lem}
\begin{proof}Let us fix a complex plane $\,V$ and a \feic\ 
identification of a \nb\ of $\,y\,$ in $\,\y\,$ with a \nb\ of 
$\,0\hs$ in $\,V\nnh$, under which $\,y=0\,$ and the almost complex structure 
$\,J\,$ of $\,\y\,$ equals, at $\,0$, the standard structure of $\,V$ 
(represented by the symbol $\,i$). Thus, $\,u,iu\in V\nh=T_y\y$. For real 
$\,t,s\,$ near $\,0\,$ we now set $\,f(t+is)=(t+is)u+t^2v+F(t,s)$, where 
$\,v\in V$ and $\,F\,$ is a $\,C^\infty$ function 
$\,F:(-\hs c,c)\times(-\hs\ve,\ve)\to V\nnh$, with sufficiently small 
$\,c,\ve$, such that $\,F(t,s)=0\,$ whenever $\,ts=0$. Clearly, $\,\df_t$ (or, 
$\,\df_{is}$) sends the real basis $\,\{1,i\}\,$ of $\,\bbC=T_t\bbC\,$ (or, 
$\,\bbC=T_{is}\bbC$) to $\,\{u+2tv,iu+F_s(t,0)\}$ (or, respectively, to 
$\,\{u+F_t(0,s),iu\}$).

Thus, to ensure com\-plex-lin\-e\-ar\-i\-ty of all $\,\df_t$ and $\,\df_{is}$, 
we need to choose $\,F\,$ and $\,v\,$ so that $\,F_s(t,0)=\varphi(t)\,$ and 
$\,F_t(0,s)=\psi(s)\,$ for all $\,t,s$, where $\,\varphi(t)=Ju+2tJv-iu\,$ and 
$\,\psi(s)=-\hs u-Jiu$. Note that $\,u,iu,v\,$ are vectors (i.e., constant 
vector fields on $\,V$), while $\,J\,$ denotes here the real-lin\-e\-ar 
operator $\,J_x:V\to V\nnh$, depending on $\,x=f(t)\,$ or, respectively, 
$\,x=f(is)$. Finally, $\,F\,$ with $\,F_s(t,0)=\varphi(t)\,$ and 
$\,F_t(0,s)=\psi(s)$ exists, by Lemma~\ref{zaxes}, if and only if the sum 
$\,d_u(Ju)+d_{iu}(Jiu)\,$ of directional derivatives at $\,0\,$ equals 
$\,-\hs2iv$. With this choice of $\,v$, our assertion follows.
\end{proof} 
Let $\,\x\,$ be a \tr\ surface embedded in an \acsu\ $\,\y$, and contained in 
a fixed $\,3$\diml\ real submanifold $\,\n\hs$ of $\,\y$. The {\it 
characteristic foliation\/} of $\,\x\,$ in $\,\n\hs$ is the $\,1$\diml\ 
foliation on $\,\x$ tangent, at every point $\,x\in\x$, to the {\it 
characteristic direction\/} $\,L_x\cap T_x\x$, for the unique complex line 
$\,L_x$ in $\,T_x\y\,$ with $\,0\in L_x\subset T_x\n\nh$. By a {\it 
characteristic curve\/} of $\,\x\,$ in $\,\n\hs$ we mean any leaf of its 
characteristic foliation.
\begin{rem}\label{chdir}Given $\,\y,\n\hs$ and $\,\x\,$ as above, let $\,\k\,$ 
be another \rsu\ embedded in $\,\y\,$ with $\,\k\subset\n\nh$. If 
$\,y\in\k\cap\x$, while $\,T_y\x\cap T_y\k\,$ is one\diml\ and different from 
the characteristic direction, at $\,y$, of $\,\x\,$ in $\,\n\nh$, then 
$\,\k\,$ is \tr\ at $\,y$, in the sense that $\,y\,$ is not a complex point of 
$\,\k$. (In fact, if $\,T_y\k\,$ were a complex line, the characteristic 
direction of $\,\x\,$ at $\,y\,$ would, by definition, be $\,T_y\x\cap T_y\k$.)
\end{rem}
\begin{lem}\label{spcoo}Let\/ $\,\y,\x^\pm\nnh,x^\pm$ satisfy the hypotheses 
of Theorem~\ref{cnsum}, and let\/ $\,U_\ve$, with 
$\,\ve\in(0,\infty)$, be the set of $\,(z,w)\in\bbC^2$ such that\/ 
$\,|\text{\rm Re}\,z|-1,\hs|\text{\rm Im}\,z|$, $\,|\text{\rm Re}\,w|\,$ and\/ 
$\,|\text{\rm Im}\,w|\,$ are all less than $\,\ve$. If the submanifolds 
$\,\n\subset\hs U_\ve$ and $\,\x_t\subset\n\nh$, for $\,t\in(-1-\ve,1+\ve)$, 
are given by $\,\n=\{(z,w)\in\hs U_\ve:\text{\rm Im}\,w=0\}$ and\/ 
$\,\x_t=\{(z,w)\in\hs U_\ve:\text{\rm Re}\,z=t\,\mathrm{\ and\ }\,
\hs\text{\rm Im}\,w=0\}$, then the following conclusions 
hold for some $\,\ve>0\,$ and some $\,C^\infty\nnh$-\feic\ 
identification of\/ $\,\,U_\ve$ with an open set in $\,\y$.
\begin{enumerate}
  \def\theenumi{{\rm\alph{enumi}}}
\item $x^\pm\nh=(\pm1,0)\in\hs U_\ve\subset\bbC^2$ and\/ 
$\,\x^\pm\nh\cap\hs U_\ve=\x_{\pm1}$.
\item Each $\,\x_t$ is a totally real surface embedded in $\,\y\hs$ and its 
characteristic foliation in the ambient\/ $\,3$\mfd\/ $\,\n$ consists of line 
segments parallel to the\/ $\,\hs\text{\rm Im}\,z\,$ coordinate direction.
\item The surface\/ $\,Y\nnh=\{(z,w)\in\n:\text{\rm Im}\,z=0\}\,$ is totally 
real in $\,\y\,$ and its characteristic foliation relative to $\,\n$ consists 
of line segments parallel to the $\,\hs\text{\rm Re}\,z\,$ coordinate 
direction.
\item The almost complex structure\/ $\,J\,$ of\/ $\,\y\,$ coincides with the 
standard complex structure of the open set\/ $\,\,U_\ve\subset\bbC^2$ at all 
points $\,(z,w)\in\hs U_\ve$ such that\/ {\rm\hs(i)} 
$\,\hs\text{\rm Im}\,z=w=0$, or\/ {\rm\hs(ii)} $\,\hs\text{\rm Re}\,z=w=0$.
\item The coordinate vector fields\/ $\,u\,$ and\/ $\,v\,$ corresponding to 
the $\,\hs\text{\rm Re}\,z\,$ and\/ $\,\hs\text{\rm Re}\,w\,$ directions are 
$\,J$-lin\-e\-ar\-ly independent at every point of\/ $\,\,U_\ve$.
\end{enumerate}
\end{lem}
\begin{proof}Let us fix a point $\,y\in\y\smallsetminus(\x^+\nh\cup\x^-)\,$ 
and a nonzero vector $\,u\in T_y\y$, then choose an embedding $\,f\,$ with the 
properties listed in Lemma~\ref{crvcp} and three coordinate systems, on 
\nb s of $\,x^+\nnh,x^-$ and $\,y$, under which $\,x^+\nnh,x^-$ and 
$\,y\,$ all have zero coordinates, while some \nb s of $\,x^\pm$ in 
$\,\x^\pm$ correspond to \nb s of $\,(0,\nh0,\nh0,\nh0)\,$ in the 
submanifold $\,\rto\times\{(0,\nh0)\}\,$ of $\,\bbR^4\nnh$, and $\,f\,$ 
appears, in the coordinates chosen at $\,y$, as the mapping 
$\,\bbC\ni z\mapsto(z,0)\in\bbC^2\nh=\bbR^4$ restricted to a \nb\ of 
$\,0\in\bbC\hs$.

Let $\,\hs\exp\hs\,$ now be the exponential mapping of a Riemannian metric on 
$\,\y\,$ chosen so that, on some \nb s of $\,x^+\nnh,x^-$ and $\,y$, 
the metric is flat and the above coordinate mappings are isometries onto 
\nb s of $\,(0,\nh0,\nh0,\nh0)\,$ in $\,\bbR^4$ with the Euclidean  
metric. (Such a metric on $\,\y\,$ can easily be obtained using a finite 
partition of unity.)

A $\,C^\infty$ \feo\ $\,H:U_\ve\to H(U_\ve)\,$ onto an open set 
$\,H(U_\ve)\subset\y$, providing the required \feic\ identification, 
will be constructed in four successive steps, each of which consists in 
defining $\,H\,$ just on some subset of $\,\,U_\ve$ (and making $\,\ve\,$ 
smaller, if necessary). The four subsets are: the interval $\,I_\ve$ of all 
$\,(t,0)\in\bbC^2$ such that $\,t\in\bbR\,$ and $\,|t|<1+\ve$, the rectangle 
$\,R_\ve=\{(z,0)\in\bbC^2:|\hh\text{\rm Re}\,z|<1+\ve,\,\,
|\hh\text{\rm Im}\,z|<\ve\}$, the $\,3$\mfd\ $\,\n_\ve$ appearing in the 
statement of the lemma (where it is denoted by $\,\n$), and $\,\,U_\ve$ 
itself. Thus, $\,I_\ve\subset R_\ve\subset\n_\ve\subset U_\ve$.

We begin by restricting $\,f\,$ chosen earlier in this proof to a small 
\nb\ of $\,0\,$ in $\,\bbR\subset\bbC\,$ and then extending this 
restriction to a $\,C^\infty$ embedding $\,t\mapsto x(t)\in\y\,$ defined for 
real $\,t\,$ with $\,|t|<1+\ve$, such that $\,x(\pm1)=x^\pm$ and 
$\,J\dot x(\pm1)\,$ is, for either sign $\,\pm\hs$, tangent to $\,\x^\pm$ at 
$\,x^\pm$. (To realize the latter condition, it suffices to prescribe the 
velocity $\,\dot x(\pm1)\,$ as a vector in $\,J(T_{x(\pm1)}\x^\pm)$.) We now 
choose $\,H:I_\ve\to\y\,$ by setting $\,H(t,0)=x(t)$.

Making $\,\ve\,$ smaller, we extend $\,H\,$ from the interval $\,I_\ve$ to the 
rectangle $\,R_\ve$ by setting $\,H(z,0)=\hs\exp_{\hs x(t)}sJ\dot x(t)$, 
where $\,t,s\in\bbR\,$ and $\,z=t+is$. Note that $\,H\,$ is injective on 
$\,I_\ve$, and so is the differential of $\,H:R_\ve\to\y\,$ at each point of 
$\,I_\ve$. Hence $\,H:R_\ve\to\y\,$ is, for small $\,\ve>0$, an embedding.

Next, we choose a $\,C^\infty$ vector field 
$\,I_\ve\ni(t,0)\mapsto v(t)\in T_{x(t)}\y$, nowhere tangent to the 
embedded surface $\,H(R_\ve)$, and such that $\,v(\pm1)\,$ is tangent to 
$\,\x^\pm$ at $\,x^\pm\nnh$. Its $\,C^\infty$ extension 
$\,R_\ve\ni(z,0)\mapsto\tilde v(z,0)\in T_{H(z,0)}\y$, obtained using parallel 
transports along the curves $\,s\mapsto H(t+is,0)$, for $\,t,s\,$ as before, 
gives rise to the mapping $\,\tilde H:\n\to\y\,$ with 
$\,\tilde H(z,r)=\hs\exp_{\hs H(z,0)\hh}r\tilde v(z,0)$, where $\,\n=\n_\ve$. 
As in the last paragraph, $\,\tilde H\,$ is an embedding for sufficiently 
small $\,\ve>0$. Also, for small $\,\ve$, the $\,\tilde H$-im\-ages of the 
surfaces $\,Y\subset\n\hs$ and $\,\x_t\subset\n\hs$ with $\,|t|<1+\hs\ve$, 
defined in the lemma, are all \tr\ in $\,\y$, while $\,\tilde H(I_\ve)\,$ (or, 
$\,\{\tilde H(is,0):s\in(-\ve,\ve)\}$) is a characteristic curve of 
$\,\tilde H(Y)\,$ (or, respectively, of $\,\tilde H(\x_0)$) in 
$\,\tilde H(\n)$, and the $\,\tilde H$-im\-age of the 
$\,\hs\text{\rm Im}\,z\,$ coordinate direction at $\,(t,0)\in\n\hs$ is, for 
each $\,t\,$ with $\,|t|<1+\ve$, the characteristic direction, at 
$\,\tilde H(t,0)$, of $\,\tilde H(\x_t)\,$ in $\,\tilde H(\n)$.

In fact, $\,\tilde H\,$ restricted to $\,R_\ve$ coincides with $\,H$, while, 
if $\,\ve\hh>0\,$ is sufficiently small, the mapping $\,z\mapsto H(z,0)\,$ 
defined for $\,z\in\bbC\,$ such that $\,|\hh\text{\rm Re}\,z|<1+\ve$ and 
$\,|\hh\text{\rm Im}\,z|<\ve$, has a com\-plex-lin\-e\-ar differential at 
every $\,z\,$ with $\,\hs\text{\rm Re}\,z=0$ or $\,\hs\text{\rm Im}\,z=0$. 
(This is clear from our definition of $\,H(z,0)\,$ and the 
com\-plex-lin\-e\-ar\-i\-ty assertion in Lemma~\ref{crvcp}, since 
$\,H(z,0)=f(z)\,$ for $\,z\,$ near $\,0$ in $\,\bbC\,$ due to our choice of 
coordinates and the metric on a \nb\ of $\,y$.) Each of the directions 
in question, including those of the two curves, is characteristic, as it 
constitutes the intersection of the plane tangent to the respective surface 
with the complex line forming the image of one of the differentials just 
mentioned; the surfaces $\,Y$ and $\,\x_t$, being totally real at the points 
of the two curves (since the tangent-plane intersections are $\,1$\diml), 
become totally real everywhere when $\,\ve\,$ is made smaller.

Let $\,\pro:\tilde H(\n)\to\tilde H(R_\ve)\,$ be the mapping with 
$\,\pro(\tilde H(z,r))=\tilde H(z,0)$, for $\,\n=\n_\ve$. If $\,\ve>0\,$ is 
sufficiently small, we extend $\,H:R_\ve\to\y\,$ to an embedding 
$\,H\nnh:\nh\n\to\y\,$ by requiring that, for real $\,t,s,r$ with 
$\,|t|<1+\hs\ve\,$ and $\,s,r\in(-\hs\ve,\ve)$, the image of the mapping 
$\,t\mapsto\nh H(t,r)\,$ (or, $\,s\mapsto H(t+is,r)$) be a characteristic 
curve of $\,\tilde H(Y)\,$ (or, respectively, of $\,\tilde H(\x_t)$) in 
$\,\tilde H(\n)$, and $\,\hs\pro(H(t+is,r))=H(t+is,0)$. Such a $\,C^\infty$ 
extension of $\,H\,$ exists: its defining conditions mean that 
$\,t\mapsto\nh H(t,r)\,$ (or, $\,s\mapsto\nh H(t+is,r)$) is an integral curve 
of the unique vector field on $\,\tilde H(Y)\,$ (or, on each 
$\,\tilde H(\x_t)$) which is tangent to the characteristic foliation and 
projects under $\,\pro\,$ onto the restriction to $\,Y\nh\cap R_\ve$ of the 
coordinate vector field $\,u\,$ for the $\,\hs\text{\rm Re}\,z$ direction 
(or, respectively, onto the restriction of $\,iu\,$ to $\,\x_t\cap R_\ve$).

Finally, we extend $\,H:\n\to\y\,$ to an embedding $\,H:U_\ve\to\y$, for small 
$\,\ve$, by setting 
$\,H(z,w)=\hs\exp_{\hs H(z,\hs\text{\rm Re}\,w)\hs}[(\text{\rm Im}\,w)Jv]$, 
which clearly yields (d) and (e). As the remainder of our assertion is an 
obvious consequence of how we chose $\,H\,$ on $\,\n\nh$, this completes the 
proof.
\end{proof}

\section{Proof of Theorem~\ref{cnsum}}\label{pc}
\setcounter{equation}{0}
We use the notations of Lemma~\ref{spcoo}, treating $\,\,U_\ve$ as an open 
subset of both $\,\y\,$ and $\,\bbC^2\nnh$. First, let us fix functions 
$\,\rho\,$ and $\,\my\,$ of the variable $\,t\in[-1,1]$ such that $\,\my\,$ 
is of class $\,C^\infty$ on $\,[-1,1]$, has a compact support contained in 
$\,(-1,1)$, and $\,0\le\my\le\my(0)<\ve\,$ on $\,[-1,1]$, with 
$\,\my(0)>0$, while $\,\rho\,$ is even, continuous on $\,[-1,1]$, of class 
$\,C^\infty$ on $\,(-1,1)$, with $\,\ddot\rho(0)>0\,$ and $\,\dot\rho>0$ on 
$\,(0,1)$, the inverse of $\,\rho:[\hs0,1]\to[\rho(0),\rho(1)]\,$ has 
derivatives of all orders equal to $\,0\,$ at $\,\rho(1)$, and 
$\,0<\rho(0)\le\rho\le\rho(1)<\ve\,$ on $\,[-1,1]$. (We write 
$\,(\,\,)\dot{\,}\,$ for $\,d/dt$). Suppressing in  our notation the 
dependence of $\,\rho\,$ and $\,\my\,$ on $\,t$, we define a mapping 
$\,h:[-1,1]\times S^1\to\,U_\ve$, by
\begin{equation}\label{hte}
h(t,e^{i\theta})=(z,w)\hs\hskip4.5pt\mathrm{with}\hskip5.5pt(z,w)
=(t+i\rho\cos\theta\hs,t\rho\sin\theta+i\hs\rc\hh\my\sin2\theta)\hs.
\end{equation}
Thus, $\,h\,$ depends on a parameter $\,\rc\in[-1,1]\,$ and has the partial 
derivatives
\begin{equation}\label{hst}
\begin{array}{l}
h_t\,=\,(1+i\dot\rho\cos\theta\hs,\,(\rho+t\dot\rho)\sin\theta
+i\hs\rc\hh\dot\my\sin2\theta)\hs,\\
h_\theta\,=\,(-i\rho\sin\theta\hs,\,t\rho\cos\theta
+2\hh i\hs\rc\hh\my\cos2\theta)\hs.
\end{array}
\end{equation}
Also, $\,h\,$ is injective except for one double point 
$\,(0,\nh0)=h(0,\pm\hs i)$, when $\,\rc>0$, or a curve of double points 
$\,h(0,e^{\pm i\theta})$, $\,0<\theta<\pi$, for $\,\rc=0$. (By (\ref{hte}), 
$\,t=\,\text{\rm Re}\,z$, and $\,e^{i\theta}=\alpha/|\alpha|\,$ for 
$\,\alpha=\,\text{\rm Im}\,z+i\,\text{\rm Re}\,w/t\hh$, if $\,t\ne0$.)

If $\,\rc>0$, the self-in\-ter\-sec\-tion of $\,h\,$ at 
$\,(0,\nh0)=h(0,\pm\hs i)\,$ is transverse and negative, as one easily sees 
using the bases, provided by (\ref{hst}), of the two real planes tangent to 
the image of $\,h\,$ at $\,(0,\nh0)$. (Note that $\,\dot\rho(0)=\dot\my(0)=0$.)

If $\,\rc=0\,$ and $\,(t,e^{i\theta})=(0,\pm1)$, while $\,\rho,\ddot\rho\,$ 
stand for $\,\rho(0),\ddot\rho(0)$, then
\begin{equation}\label{hti}
\begin{array}{rl}
\mathrm{i)}&\hskip8pth_t=(1,0)\hs,\hskip4pth_{tt}
=(\pm\hs i\ddot\rho,0)\hs,\hskip4pth_{\theta t}
=(0,\pm\hs\rho)\hs,\hskip4pth_{\theta\theta}=(\mp\hs i\rho,0)\hs,\\
\mathrm{ii)}&\hskip8pth_\theta=h_{\theta tt}=h_{\theta\theta t}
=h_{\theta\theta\theta}=(0,\nh0)
\end{array}
\end{equation}
from (\ref{hst}) with $\,\dot\rho(0)=0$. Also, by (\ref{hte}) -- (\ref{hst}), 
for all $\,(t,\theta,\rc)$,
\begin{equation}\label{imj}
\text{\rm Im}\,\mathcal{J}=2\hh\rc\my
+(\rho^2-4\hh\rc\my)\sin^{\hs2\nnh}\theta+t\rho\hs\dot\rho\hs,
\hskip6pt\mathrm{where}\hskip5pt\mathcal{J}\nh=z_tw_\theta-z_\theta w_t\hs.
\end{equation}
(Thus, $\,\mathcal{J}\,$ is the Jacobian of $\,h$.) Hence 
$\,\hs\text{\rm Im}\,\mathcal{J}\ge0\,$ at any $\,(t,\theta)\,$ and for any 
$\,\rc\in[\hs0,\rc_{\text{\rm max}}]$, where 
$\,\rc_{\text{\rm max}}=\hs\text{\rm min}\hs(1,[\rho(0)]^2/[4\my(0)])$. (In 
fact, the three terms 
$\,2\hh\rc\my,\,(\rho^2-4\hh\rc\my)\sin^{\hs2\nnh}\theta\,$ and 
$\,t\rho\hs\dot\rho\,$ are all nonnegative.) Moreover, the strict inequality 
$\,\hs\text{\rm Im}\,\mathcal{J}>0\,$ clearly holds unless $\,\rc=t=0\,$ and 
$\,\theta\in\bbZ\pi$.

Therefore (cf.\ Example~\ref{trsub}), $\,h\,$ with any fixed 
$\,\rc\in(0,\rc_{\text{\rm max}}]\,$ is a \tri\ for the standard complex structure 
of $\,\,U_\ve\subset\bbC^2\nnh$. This is also true if $\,\rc=0$, provided that 
one excludes the two points with $\,t=0\,$ and $\,\theta\in\bbZ\pi\,$ (at 
which $\,dh\,$ is not even injective, cf.\ (\ref{hti})). We will now show that 
the same conclusions hold for the original almost complex structure $\,J$ of 
$\,\y$, as long as one chooses $\,\rho\,$ more carefully, and $\,\rc\,$ is 
sufficiently small.

First, suppose that $\,\rc=0$. As we just saw, at $\,(t,e^{i\theta})\,$ with 
$\,t=0$ and $\,\theta\notin\bbZ\pi\,$ the vectors $\,h_t,h_\theta$ are 
linearly independent in $\,\bbC^2\nnh$. Hence, by Lemma~\ref{spcoo}(d) (case 
(ii)), they are also $\,J$-lin\-e\-ar\-ly independent. Now let 
$\,y=h(t,e^{i\theta})\,$ with $\,t\ne0$. The (disconnected) surface 
$\,\k\subset\y$, formed by all such $\,y$, is contained in the 
$\,3$\diml\ submanifold $\,N\hs$ defined in Lemma~\ref{spcoo}. Thus, $\,\k\,$ 
is $\,J$-to\-tal\-ly real at $\,y\,$ in view of Remark~\ref{chdir} applied to 
a suitable \tr\ surface $\,\x\subset\n\nh$, which intersects $\,\k\,$ 
transversely in $\,\n\hs$ at $\,y\,$ along a non-char\-ac\-ter\-is\-tic 
direction; namely, $\,\x=Y\nh$ if $\,e^{i\theta}=\pm\hs i$ and $\,\x=\x_t$ 
otherwise. (Notation of Lemma~\ref{spcoo}.) That the intersection is 
non-char\-ac\-ter\-is\-tic follows from (b) -- (c) in Lemma~\ref{spcoo}: 
$\,\k\cap\x_t$, for each $\,t\ne0$, is an ellipse in standard position 
relative to the $\,\hs\text{\rm Im}\,z\,$ and $\,\hs\text{\rm Re}\,w$ 
coordinate axes, while the curve $\,\k\nh\cap Y\nh$ is pa\-ram\-e\-trized by 
$\,t\mapsto(z,w)=(t,\pm\hs t\rho)$, with $\,(t\rho)\dot{\,}>0\,$ as 
$\,t\dot\rho\ge0\,$ and $\,\rho>0$.

The next five paragraphs deal with the remaining case $\,\rc>0$.

According to Lemma~\ref{spcoo}(e) , the $\,J$-com\-plex exterior product of 
$\,u$ and $\,v\,$ trivializes the line bundle $\,[\tm]^{\wedge2}$ restricted 
to $\,\,U_\ve\subset\y$, allowing us to treat the $\,J$-com\-plex exterior 
product of any vector fields $\,\xi,\xi'$ on $\,\,U_\ve$ as a function 
$\,\,U_\ve\to\bbC\hh$. The imaginary part $\,B(\xi,\xi')\,$ of that function 
depends on $\,\xi\,$ and $\,\xi'$ only through their values at points of 
$\,\,U_\ve$ (since so does the function itself); we may therefore view $\,B\,$ 
as a differential $\,2$-form on $\,\,U_\ve$, i.e., as a mapping 
$\,B:U_\ve\to\mathcal{X}\,$ valued in the space $\,\mathcal{X}\,$ of 
skew-sym\-met\-ric bi\-lin\-ear forms $\,\bbC^2\nnh\times\bbC^2\to\bbR\hs$. Note 
that $\,B(u,v)=0\,$ everywhere in $\,\,U_\ve$.

Let us now set $\,\vd=B(h_t,h_\theta)$, treating $\,\vd\,$ as a real-val\-ued 
function of $\,(t,\theta,\rc)\,$ (where $\,\vd\hs$ depends on $\,\rc\,$ since 
$\,h\,$ does). If $\,\rho,\ddot\rho,\my\,$ denote $\,\rho(0)$, 
$\,\ddot\rho(0)$ and $\,\my(0)$, then, for the value and partial derivatives 
of $\,\vd$ at the points $\,\op_+,\op_-\in\rtr$ that have the 
$\,(t,\theta,\rc)$ coordinates $\,(0,0,0)\,$ and $\,(0,\pi,0)$, we get
\begin{equation}\label{itt}
\begin{array}{rl}
\mathrm{a)}\hskip0pt&\vd\,=\,\vd_t\,=\,\vd_{\theta}\,=\,\hs0\hs,\hskip20pt
\vd_{\theta t}\,=\,\mp\rho\hs(d_uB)(u,iu)\hs,\\
\mathrm{b)}\hskip0pt&\vd_\rc\,=\,2\hs\my>0\hs,\hskip9pt\vd_{tt}
\,=\,2\rho\ddot\rho>0\hs,\hskip9pt\vd_{\theta\theta}\,=\,2\rho^2\nh>0\hs,\\
\mathrm{c)}\hskip0pt&B_\theta\,=\,B_\rc\,=\,\hs0\hskip7pt\mathrm{and}\hskip7pt
B_t\,=\,d_uB\hskip7pt\mathrm{at}\hskip6pt\op_\pm\hs.
\end{array}
\end{equation}
Here $\,B_t,B_\theta,B_\rc$ are the $\,\mathcal{X}$-val\-ued partial derivatives 
of $\,B\circ h$, while $\,u$ denotes the $\,\hs\text{\rm Re}\,z\,$ coordinate 
vector field, $\,d_uB:U_\ve\to\mathcal{X}\,$ is the corresponding directional 
derivative of $\,B:U_\ve\to\mathcal{X}$, and $\,i\,$ in $\,iu\,$ refers to the 
complex structure of $\,\bbC^2\nnh$.

In fact, we have (\ref{itt}.c) and $\,\vd\nh=0\,$ at $\,\op_\pm$ as 
$\,h_\theta=h_\rc=0\,$ and $\,h_t=u$ at $\,\op_\pm$ by (\ref{hti}) and 
(\ref{hte}). Now let $\,A\,$ be the analogue of $\,B\,$ obtained by using, 
instead of $\,J$, the complex structure of $\,\bbC^2\nnh$. Hence $\,A\,$ is 
constant as a function $\,U_\ve\to\mathcal{X}$, while, by 
Lemma~\ref{spcoo}(d), $\,B=A\,$ at all points of $\,\,U_\ve$ that have the 
form $\,(t,0)\,$ or $\,(is,0)\,$ with $\,t,s\in\bbR\hs$, including, when 
$\,\rc=0$, the $\,h$-im\-ages $\,(0,\pm\hs\rho(0))\,$ of 
$\,(t,e^{i\theta})=(0,\pm1)$. (Also, 
$\,A(h_t,h_\theta)=\hs\text{\rm Im}\,\mathcal{J}$, cf.\ (\ref{imj}.)) To 
evaluate the partial derivatives of $\,\vd\hs$ at $\,\op_\pm$ in (\ref{itt}), 
we differentiate the three-fac\-tor ``product'' $\,B(h_t,h_\theta)$ via the 
Leibniz rule. The partial derivatives in question differ from those of 
$\,A(h_t,h_\theta)\,$ just by the {\it extra terms\/} in which the $\,B\,$ 
factor is differentiated.  However, by (\ref{hti}.ii), such an extra term can 
only be nonzero if exactly one differentiation falls on the $\,h_\theta$ 
factor. This, combined with (\ref{itt}.c), shows that the extra terms are all 
zero, except, possibly, those in $\,\vd_{tt}$ and $\,\vd_{\theta t}$, equal to 
$\,\pm\hs2\rho\hh(d_uB)(u,v)$ and, respectively, $\,\mp\rho\hs(d_uB)(u,iu)$. 
(By (\ref{hti}.i), $\,h_t=u$, $\,h_{\theta t}=\pm\hs\rho v\,$ and 
$\,h_{\theta\theta}=\mp\hs\rho\hs iu\,$ at $\,\op_\pm$.) Now (\ref{itt}) 
follows, since $\,B(u,v)=0\,$ identically on $\,\,U_\ve$, and so 
$\,(d_uB)(u,v)=d_u[B(u,v)]=0$, while the partial derivatives of 
$\,A(h_t,h_\theta)\,$ at $\,\op_\pm$ are easily found using the Leibniz rule 
and (\ref{hti}): for instance, 
$\,[A(h_t,h_\theta)]_t=A(h_t,h_{\theta t})=\pm\hs\rho(0)\hs A(u,v)=0$.

Our functions $\,\rho\,$ and $\,\my\,$ were so far subject to the specific 
conditions listed above, but otherwise arbitrary. However, a special choice of 
$\,\rho\,$ gives
\begin{equation}\label{dti}
\vd_{tt}\vd_{\theta\theta}\,-\,(\vd_{\theta t})^2\,
>\,0\hskip8pt\mathrm{at\ the\ two\ points}\hskip6pt(t,\theta,\rc)=\op_\pm\hs.
\end{equation}
Namely, by (\ref{itt}), $\,\vd_{tt}\vd_{\theta\theta}-(\vd_{\theta t})^2
=4\rho^3\nh\ddot\rho-[\rho\hs(d_uB)(u,iu)]^2\nh$. Although the function 
$\,(d_uB)(u,iu):U_\ve\to\bbR\,$ does not depend on the choice of $\,\rho$, its 
value used here does, as it is taken at the point 
$\,h(0,\pm1)=(\pm\hs i\rho(0),0)$. Modifying $\,\rho\,$ so as to keep 
$\,\rho(0)>0\,$ unchanged while making $\,\ddot\rho(0)$ so large that 
$\,4\rho(0)\ddot\rho(0)>[(d_uB)(u,iu)]^2$ at $\,(\pm\hs i\rho(0),0)$, we 
obtain (\ref{dti}).

As $\,\vd_\rc>0=\vd\,$ at $\,\op_\pm$ (see (\ref{itt})), the implicit function 
theorem applied to the function $\,(t,\theta,\rc)\mapsto\vd\,$ shows that, in 
some \nb\ of either point $\,\op_\pm$ in $\,\rtr\nh$, the equation 
$\,\vd=0\,$ describes the graph of a $\,C^\infty$ function 
$\,(t,\theta)\mapsto\rc$, equal to $\,0\,$ at $\,(0,\nh0)\,$ or, respectively, 
at $\,(0,\pi)$. Partial differentiations of the equality $\,\vd=0\,$ with 
respect to $\,t,\theta\,$ show that, as a consequence of (\ref{itt}) and 
(\ref{dti}), $\,0\,$ is a nondegenerate lo\-cal-max\-i\-mum value of the 
function $\,(t,\theta)\mapsto\rc\,$ at $\,(0,\nh0)\,$ (or, $\,(0,\pi)$). Thus, 
we may choose \nb s of $\,(0,\nh0)\,$ and $\,(0,\pi)\,$ on which 
$\,\vd=B(h_t,h_\theta)\,$ is nonzero everywhere for all sufficiently small 
$\,\rc>0$. Due to the definition of $\,B\,$ and Example~\ref{trsub}, this shows 
that $\,h$, with any such fixed $\,\rc$, is a $\,J$-to\-tal\-ly real immersion 
when restricted to the two corresponding \nb s of $\,(0,\pm1)$ in 
$\,[-1,1]\times S^1\nh$. On the other hand, removing the latter two 
\nb s from $\,[-1,1]\times S^1\nh$, we are left with a compact set 
on which $\,h\,$ is $\,J$-to\-tal\-ly real for any small $\,\rc>0\,$ (since that 
is the case for $\,\rc=0$, as we saw earlier).

Thus, our $\,h$, with any fixed $\,\rc>0\,$ close to $\,0$, is a 
$\,J$-to\-tal\-ly real immersion with a single, transverse, negative 
self-in\-ter\-sec\-tion, while $\,\,U_\ve$ is an obvious choice of the open 
subset required in (v). The only claim that still needs proving is smoothness 
of the union $\,S\,$ of the image of $\,h$ with the complements of the disks 
$\,D^\pm$ of radius $\,\rho(1)$, centered at $\,(\pm1,0)$, in the Euclidean 
planes $\,\{(z,w)\in\bbC^2:\text{\rm Re}\,z\mp1=\,\text{\rm Im}\,w=0\}$. To 
this end, let us fix $\,\delta\in(0,1)\,$ with 
$\,\hs\text{\rm supp}\,\my\subset[-\hs\delta,\delta\hs]$, and denote by 
$\,\n'$ the disconnected $\,3$\mfd\ formed by all $\,(z,w)\in\bbC^2$ with 
$\,|\text{\rm Re}\,z|>\delta\,$ and $\,\hs\text{\rm Im}\,w=0$. Smoothness of 
$\,S\,$ now follows from smoothness of the (disconnected) surface 
$\,\x\hs'\subset\n'$ which consists of $\,(z,w)\in\n'$ such that, setting 
$\,t=\,\text{\rm Re}\,z\,$ and defining $\,\er\ge0\,$ by 
$\,\er^2=(\text{\rm Im}\,z)^2+\hs(\text{\rm Re}\,w)^2\nnh$, one has either 
$\,|t|=1\,$ and $\,\er\ge\rho(1)$, or $\,\delta<|t|<1\,$ and $\,\er=\rho(t)$. 
In fact, $\,\x\hs'$ is a smooth surface in the Euclidean $\,3$-space 
$\,\bbC\times\bbR\subset\bbC^2$ with the coordinates 
$\,\hs\text{\rm Re}\,z,\,\text{\rm Im}\,z,\,\text{\rm Re}\,w$, since it is 
obtained by revolving, about the $\,\hs\text{\rm Re}\,z\,$ axis, a 
(disconnected) $\,C^\infty$ curve lying in the half-plane 
$\,\hs\text{\rm Im}\,z>0=\,\text{\rm Re}\,w$. The curve is the intersection 
of $\,\x\hs'$ with the half-plane, so that its smoothness becomes obvious if 
one pa\-ram\-e\-trizes either of its components using $\,\er\,$ (the distance 
from the $\,\hs\text{\rm Re}\,z\,$ axis) as the parameter, and recalls the 
assumptions made earlier about the derivatives of the inverse of $\,\rho$. 
Finally, the set $\,S\hs'\nh=S\cap\n'$ is relatively open in $\,S\,$ (as 
$\,S\hs'=\{(z,w)\in S:|\text{\rm Re}\,z|>\delta\}$), while $\,S\hs'$ is also 
the image of $\,\x\hs'$ under the \feo\ $\,\n'\nh\to\n'$ sending 
$\,(z,w)\,$ to $\,(z,(\text{\rm Re}\,z)w)$.

The proof of Theorem~\ref{cnsum} is now complete.

\section{Totally real tori and Klein bottles in $\,\bbC^2$}\label{tk}
\setcounter{equation}{0}
We begin this section with some simple examples of \tri s and embeddings of 
the torus $\,T^2$ and \kb\ $\,K^2$ in $\,\bbC^2$. 

For $\,C^1$ functions $\,x,y:U\to\bbC\,$ on a open set $\,\,U\,$ in the 
$\,(t,s)$-plane $\,\rto$, we define the function 
$\,\,\mathcal{J}(x,y):U\to\bbC\,$  by $\,\mathcal{J}(x,y)=x_sy_t-x_ty_s$ (cf.\ 
Example~\ref{trsub}); the subscripts stand for the partial derivatives.
\begin{example}\label{exbdd}Let $\,x,y,h:\rto\to\bbC\,$ be doubly 
$\,2\pi$-periodic $\,C^1$ functions of the variables $\,t,s\,$ such that 
$\,|\mathcal{J}(x,y)|\,$ is bounded on $\,\rto\,$ and $\,h_s=0$ identically, 
while $\,|x_sh_t|\ge\ve\,$ for some real number $\,\ve>0$. The mapping 
$\,(x,y+\ax h):\rto\to\,\bbC^2\nnh$, with any constant $\,\ax\in\bbR\hs$, then 
descends to the torus $\,T^2=\hs[\bbR/2\pi\bbZ]\times\hs[\bbR/2\pi\bbZ]\,$ 
and, for large $\,|\ax|$, it produces a \tri\ $\,f:T^2\to\bbC^2$. If, in 
addition, $\,x,y,h\,$ are all invariant under the transformation 
$\,(t,s)\,\mapsto\,(t+\pi,-s)\,$ of $\,\rto$, then $\,(x,y+\ax h)\,$ further 
descends to a totally real immersion $\,f:K^2\to\bbC^2$ of the \kb\ 
$\,K^2=\rto\nnh/\hs\Gamma=T^2/\bbZ_2$, where $\,\Gamma\,$ is the 
transformation group generated by $\,\varPhi$ and $\,\varPsi\,$ with 
$\,\,\varPhi(t,s)=(t+\pi,-s)\hs\,$ and $\,\,\varPsi(t,s)=(t,s+2\pi)$. 

In fact, $\,(x,u):\rto\to\bbC^2$ is a \tri\ \iff\ 
$\,\mathcal{J}(x,u)\ne0$ everywhere in $\,\rto$ (see Example~\ref{trsub}). 
Clearly, for any $\,r\in\bbR\hs$,
\begin{equation}\label{jrh}
\mathcal{J}(x,y+\ax h)\,=\,\mathcal{J}(x,y)\,+\,\ax x_sh_t\,,
\end{equation}
and so $\,|\mathcal{J}(x,y+\ax h)|\,\to\,\infty\,$ as $\,\ax\to\infty$, 
uniformly on $\,\,U$. Therefore, if $\,|\ax|$ is sufficiently large, 
$\,\mathcal{J}(x,y+\ax h)\ne0\,$ everywhere in $\,\rto\nh$.
\end{example}
\begin{example}\label{trtkb}The assumptions listed in Example~\ref{exbdd} 
are obviously satisfied by the functions $\,x(t,s)=e^{i\ek t}(\sin s+i\sin2s)$, 
$\,y(t,s)=e^{ilt}\cos s$, $\,h(t,s)=e^{ilt}$, where $\,\ek\,$ and $\,l\,$ are 
fixed integers with $\,l\ne0$. The mapping $\,(x,y+\ax h):\rto\,\to\,\bbC^2$ 
thus descends, for large $\,\ax$, to a \tri\ $\,f=f^{\ek,\hs l}:T^2\to\bbC^2$ of 
the $\,2$-torus; in the case where $\,\ek\,$ is odd and $\,l\,$ is even, it 
similarly descends to a \tri\ $\,f=f^{\ek,\hs l}:K^2\to\bbC^2$ of the \kb. 

Furthermore, if $\,\ek,l\,$ are integers and either 
\begin{enumerate}
  \def\theenumi{{\rm\alph{enumi}}}
\item $\,l\,=\,1$, while $\,\x\,$ is the $\,2$-torus $\,T^2\nh$, or 
\item $\,l\,=\,2$, while $\,\ek\,$ is odd and $\,\x\,$ is the \kb\ $\,K^2\nh$, 
\end{enumerate}
then, for all sufficiently large $\,\ax>1$, the mapping 
$\,f^{\ek,\hs l}:\x\to\bbC^2$ defined above is a totally real {\it embedding}.

In fact, injectivity of $\,f^{\ek,\hs l}$ follows since 
$\,(x,u)=(x(t,s),y(t,s)+\ax h(t,s))$ determines 
$\,(\alpha,\beta)=(e^{is},e^{it})\,$ either uniquely (case (a)), or up to the 
involution $\,(\alpha,\beta)\mapsto(\overline{\alpha},\,-\,\beta)\,$ (case 
(b)). Specifically, $\,\cos s=|u|-\ax$, $\,e^{ilt}=u/|u|$ and 
$\,\sin s=xe^{-i\ek t}(1+2i\cos s)^{-1}\nnh$.
\end{example}

\section{Proofs of Theorem~\ref{thdue} and Corollaries~\ref{cotre} -- 
\ref{conov}}\label{pf}
\setcounter{equation}{0}
Let $\,\y\,$ be an \acsu. Corollary~\ref{sfinm} implies that $\,S^2$ admits a 
\tri\ in $\,\y$, while \tre s $\,T^2\to\y\,$ and $\,K^2\to\y\,$ exist in view 
of Proposition~\ref{zoopr} combined with Example~\ref{trtkb}.

If \rsu s $\,\x\,$ and $\,\x\hs'$ admit \tri s $\,f\,$ and $\,f'$ in $\,\y$, 
then to\-tal\-ly-real immersibility of $\,\x\,\#\,\x\hs'$ in $\,\y\,$ is an 
obvious consequence of Theorem~\ref{cnsum} applied to 
$\,\y\smallsetminus[f(\partial D)\cup f'(\partial D\hs')]\,$ rather than 
$\,\y$, and $\,x^\pm\nnh\in\y\,$ chosen so that 
$\,f(x)=x^+\nnh\ne x^-=f'(x')\,$ for some 
$\,x\in\nh\x,\,\hs x'\in\nh\x\hs'\nnh$, along with 
$\,\x^+=f(\text{\rm Int}\hs D)\,$ and $\,\x^-=f'(\text{\rm Int}\hs D\hs')$, 
where $\,D,D\hs'$ are small closed disks embedded in $\,\x,\x\hs'$ so as to 
contain $\,x,x'$ as interior points.

Next, let $\,\x^+$ be a \compact\ \tr\ surface embedded in $\,\y$. A \tr\ 
$\,2$-torus $\,\x^-$ embedded in $\,\y\smallsetminus\x^+$ exists by 
Corollary~\ref{tretn}; in view of Theorem~\ref{cnsum}, $\,\x^+\#\,\x^-$ admits 
a \tri\ in $\,\y\,$ with just one self-in\-ter\-sec\-tion, in the form of a 
double point at which the self-in\-ter\-sec\-tion is transverse and, if 
$\,\x\,$ is orientable, also negative.

The statement following Lemma~\ref{ctout}, applied to the immersed image of 
$\,\x^+\#\,\x^-\approx\hs\x\,\#\,T^2$ (rather than $\,\x$), gives rise to a 
\tre\ of $\,\x\,\#\,T^2\#\,K^2$ in $\,\y$, thus proving Theorem~\ref{thdue}. 

Corollaries~\ref{cotre} -- \ref{cocin} now follow, and Corollary~\ref{cosei} 
is easily derived from Corollary~\ref{cocin} using a real form 
$\,\rp^2\subset\cp^2$ (see (vi) in Section~\ref{se}).

Corollary~\ref{coset} for $\,T^2$ and $\,S^2$ is clear from 
Theorem~\ref{thdue},  Lemma~\ref{blwup}(a) and 
Example~\ref{prbcp} (or, Example~\ref{rdsph}) combined with 
Remark~\ref{third}. To prove 
Corollary~\ref{coset} for nonorientable \compact\ surfaces $\,\x$, note that 
the operation $\,\x\,\mapsto\,\x\,\#\,T^2\#\,K^2$ (where $\,K^2$ is the \kb) 
reduces the Euler characteristic by $\,4$, and so, in view 
Theorem~\ref{thdue}, one needs only to show the existence of a \tre\ in 
$\,\y\,=\,\cp^2\,\#\,\,\mk\hs\overline{\cp^2}$, $\,\mk\ge1$, of any \compact\ 
nonorientable surface $\,\x\,$ such that $\,\chi(\x)\in\{-2,-1,0,1\}$. If 
$\,\chi(\x)\in\{-2,0\}$, this follows from Theorem~\ref{thdue} and 
embeddability of $\,S^2\nnh$, as $\,\x=S^2\#\,T^2\#\,K^2$. Now let 
$\,\chi(\x)=\pm1$. Thus, $\,\x=\rp^2$ or $\,\x=3\hskip.8pt\rp^2$, and our 
claim is obvious from (vi) in Section~\ref{se}, Lemma~\ref{blwup}(a) and 
Corollary~\ref{trerp}.

Finally, Corollary~\ref{conov} is immediate from Corollary~\ref{cotre}, 
(\ref{wco}.a) and (\ref{wus}).

\section{Obstructions for embedded orientable surfaces}\label{os}
\setcounter{equation}{0}
Let $\,\y\,$ be a compact \acsu. If a \compact\ oriented real 
surface $\,\x\,$ admits a totally real embedding $\,f:\x\to\y\,$ and 
$\,\sgm\in H^2(\y,\bbR)$ denotes the (real) Poincar\'e dual of 
$\,f_*[\x\hs]\in H_2(\y,\bbR)$, setting $\,\chi=\chi(\x)$ and 
$\,c_{\hs1}=c_{\hs1}(\y)\,$ we can rewrite (\ref{dot}), with $\,n=2$, and 
(\ref{wco}.b) as
\begin{equation}\label{cup}
\sgm\smallsmile\sgm\,=\,-\,\chi\,,\qquad\quad\sgm\smallsmile c_{\hs1}\,
=\,\hs0\hs.
\end{equation}
For any \compact\ oriented $\,4$-man\-i\-fold $\,\y$, the cup product 
$\,\,\smallsmile\,\,$ is nondegenerate as a real-val\-ued symmetric 
bi\-lin\-ear form in $\,H^2(\y,\bbR)$, and its sign pattern consists of 
$\,b^+\nh$ pluses and $\,b^-$ minuses, for some $\,b\hh^\pm\nnh\in\bbZ\,$ with 
$\,b_2(\y)=b^+\nnh+\hs b^-\nnh$. 
\begin{prop}\label{btwpl}Let\/ $\,\y\,$ be a \compact\ almost complex 
surface for which\/ $\,b^+\nh=1$, 
$\,c_{\hs1}\hskip-1.3pt\smallsmile c_{\hs1}\ge0\,$ and\/ $\,c_{\hs1}\ne0\,$ 
in\/ $\,H^2(\y,\bbR)$, and let a \compact\ orientable surface\/ 
$\,\x\,$ admit a totally real embedding\/ $\,f:\x\to\y$.
\begin{enumerate}
  \def\theenumi{{\rm\alph{enumi}}}
\item $\x\,$ must then be \feic\ to the torus\/ $\,T^2$ or the sphere\/ 
$\,S^2\nnh$.
\item If, in addition, $\,\y\,$ has\/ $\,b^-\nh=0$, then\/ $\,\x\,$ is 
\feic\ to\/ $\,T^2\nnh$.
\item If $\,\x\,$ is \feic\ to\/ $\,T^2$, then either 
$\,c_{\hs1}\hskip-1.3pt\smallsmile c_{\hs1}>0\,$ and\/ $\,f_*[\x\hs]=0$ in 
$\,H_2(\y,\bbR)$, or\/ $\,c_{\hs1}\hskip-1.3pt\smallsmile c_{\hs1}=0\,$ and\/ 
$\,f_*[\x\hs]\,$ is a 
real multiple of the Poin\-car\'e dual of\/ $\,c_{\hs1}$ in\/ 
$\,H_2(\y,\bbR)$. 
\end{enumerate}
\end{prop}
\begin{proof}If $\,\chi\le0$, relations (\ref{cup}) imply that the 
cup-product form $\,\,\smallsmile\,\,$ is positive semidefinite on the 
subspace $\,W\subset H_2(\y,\bbR)\,$ spanned by $\,\sgm\,$ and $\,c_{\hs1}$. 
Since $\,\,\smallsmile\,\,$ has the Lorentzian sign pattern $\,+-\ldots-$, 
this shows that $\,\dim W=1$. Using (\ref{cup}), we now obtain (c). Also, 
as the inequality $\,\chi<0\,$ would, by (\ref{cup}), 
make $\,\sgm\,$ and $\,c_{\hs1}$ linearly independent (and hence yield 
$\,\dim W=2$), we see that $\,\chi=\chi(\x)\ge0$, which proves (a). Finally, 
condition $\,b^-\nnh=0\,$ implies that $\,c_{\hs1}$ spans $\,H^2(\y,\bbR)$, so 
that (\ref{cup}) gives $\,\sgm=0\,$ and $\,\chi=0$, which yields (b). This 
completes the proof.
\end{proof} 
\begin{cor}\label{tonly}The torus $\,T^2$ is the only \compact\ orientable 
\rsu\ that admits a \tre\ in $\,\bbC^2$, $\,\cp^2$ or the \csu\ obtained by 
blowing up a point in $\,\cp^2$.
\end{cor}
\begin{proof}A \tre\ of $\,T^2$ exists according to (v) in Section~\ref{se}. 
Concerning its nonexistence for other real surfaces, the case of $\,\bbC^2$ 
follows from that of $\,\cp^2$ via the inclusion $\,\bbC^2\subset\cp^2$ (or, 
directly from (\ref{dot}) with $\,H_2(\bbC^2\nnh,\bbZ)=\{0\}$). As for 
$\,\y\,=\,\cp^2$ or $\,\y=\,\cp^2\,\#\,\,\overline{\cp^2}\nnh$, the only other 
possibility left by Proposition~\ref{btwpl}(a) is that of a \tr\ $\,2$-sphere 
embedded in $\,\y$. This in turn is excluded by Corollary~\ref{notrs}. (For 
$\,\cp^2\nnh$, we may also use Proposition~\ref{btwpl}(b).) 
\end{proof} 
\noindent{\em Proofs of Corollaries~\ref{codie} and~\ref{cound}.\/} 
In both corollaries, the 
nonexistence part is obvious from Proposition~\ref{btwpl} and 
Corollary~\ref{conov}. (Note that condition $\,\mk\le9\,$ in 
Corollary~\ref{cound} amounts to 
$\,c_{\hs1}\hskip-1.3pt\smallsmile c_{\hs1}\ge0$, cf.\ Section~\ref{su}.) 
Corollary~\ref{cound} is now immediate from Corollary~\ref{coset}. As for 
Corollary~\ref{codie}, the existence assertion 
for $\,T^2$ is clear from (v) in Section~\ref{se}, and for $\,S^2$ it is provided by 
(vii) in Section~\ref{se}. Finally, let $\,\x\,$ be a nonorientable \compact\ 
surface. Since the operation $\,\x\,\mapsto\,\x\,\#\,T^2\#\,K^2$ (where 
$\,K^2$ is the \kb) reduces the Euler characteristic of any \compact\ surface 
by $\,4$, by Theorem~\ref{thdue} we just need to show the existence of a 
\tre\ $\,\x\to\y\,=\,\cp^1\!\times\cp^1$ under the additional assumption that 
$\,\chi(\x)\in\{-2,0\}$. Now, if $\,\chi(\x)=0$, this follows from 
Corollary~\ref{coqtr}. For $\,\chi(\x)=-\,2$, we have 
$\,\x=\,S^2\#\,T^2\#\,K^2\nnh$, and so we may use the embeddability of 
$\,S^2$ along with Theorem~\ref{thdue}.

\section{Embeddings of nonorientable surfaces}\label{ns}
\setcounter{equation}{0}
The self-in\-ter\-sec\-tion formula (\ref{dtf}) often fails to detect 
non-embeddability of {\it nonorientable\/} surfaces. For instance, the genus 
$\,3\,$ surface $\,\rp^2\#\,\rp^2\#\,\rp^2$ admits no \tre\ in $\,\cp^2\nnh$, 
yet this fact cannot be derived from (\ref{dtf}). To obtain useful 
obstructions for \tre s of nonorientable \compact\ manifolds one needs more 
subtle in\-ter\-sec\-tion-theo\-retic tools, such as the {\it Pontryagin 
square\/} operation, applicable to this case via a result of Massey described 
below. Our presentation follows \cite{massey}. 

Given a manifold $\,\y\,$ and an integer $\,\ek\ge0$, the {\it Pontryagin 
square\/} \cite{massey} is a natural cohomology operation 
$\,H^\ek(\y,\bbZ_2)\to H^{2\ek}(\y,\bbZ_4)\,$ which, applied to 
$\,\hs\text{\rm mod}\hskip4pt2\,$ reductions of integral classes, assigns the 
value $\,[\hs\xi\smallsmile\xi\hskip6pt\text{\rm mod}\hskip4pt4\hs]$ to 
$\,[\hs\xi\,\text{\rm mod}\hskip4pt2\hs]$, for $\,\xi\in H^\ek(\y,\bbZ)$. The 
symbol $\,\sgm^2$ will stand for the Pontryagin square of 
$\,\sgm\in H^n(\y,\bbZ_2)\,$ in the case of a \compact\ oriented manifold 
$\,\y\,$ of dimension $\,2n\,$ (such as a \compact\ \acm\ with $\,\dimc\y=n$). 
Thus, $\,\sgm^2\in\bbZ_4$. Also, by definition, $\,\sgm^2=0\in\bbZ_4$ if 
$\,\y\,$ is \feic\ to $\,\bbC^n\nnh$. We begin with the following 
special case of a result of Massey \cite[Theorem~1]{massey}:
\begin{lem}\label{masle}Given an embedding 
$\,f:\x\to\y\,$ of a \compact\ manifold\/ $\,\x\,$ in a \compact\ oriented 
manifold\/ $\,\y\,$ with $\,\dim\y=2\hs\dim\x=2n$, where $\,n\,$ is even and\/ 
$\,n\ge2$, let\/ $\,\sgm\in H^n(\y,\bbZ_2)\,$ be the Poincar\'e dual of\/ 
$\,f_*[\x\hs]\in H_n(\y,\bbZ_2)$. If\/ $\,\chi(\nu)\in\bbZ\,$ and\/ 
$\,\sgm^2\in\bbZ_4$ are the twisted Euler number of the normal bundle\/ 
$\,\nu\,$ of\/ $\,f\,$ and the Pontryagin square of\/ $\,\sgm$, then 
\begin{equation}\label{msf}
\sgm^2\,=\,\,[\chi(\nu)\hskip1.5pt\text{\rm mod}\hskip4pt4\hs]\,
+\,2\left[w_1\smallsmile w_{n-1}\right]\in\,\bbZ_4\hs,\hskip7pt\text{\rm 
where}\hskip5ptw_j=\hs w_j(\x)\hs,
\end{equation}
$2\left[w_1\smallsmile w_{n-1}\right]$ being the image of\/ 
$\,w_1\smallsmile w_{n-1}$ under the nontrivial coefficient homomorphism 
$\,\bbZ_2\to\bbZ_4$.
\end{lem}
\begin{proof}Any embedding $\,f:\x\to\y\,$ gives rise to the associated 
mapping $\,F\,$ from $\,\y\,$ into the Thom space 
$\,\text{\rm Th}\hskip.6pt(\nu)\,$ (that is, a one-point compactification) of 
the normal bundle $\,\nu\,$ of $\,f$. The Thom space is an oriented 
pseu\-do\-man\-i\-fold, and with the usual orientation conventions, $\,F\,$ 
induces the identity homomorphism between the top (co)homology groups. 
Moreover, $\,\sgm\,$ is the $\,F$-pull\-back of 
$\,[\hskip1.2pt\mathcal{U}\,\text{\rm mod}\hskip4pt2\hskip1pt]
\in H^n(\text{\rm Th}\hskip.6pt(\nu),\bbZ_2)$, that is, of 
the (well-defined) $\,\hs\text{\rm mod}\hskip4pt2\,$ reduction of 
the Thom class $\,\mathcal{U}$. (The Thom class itself lives in the cohomology 
of another pair with twisted $\,\bbZ\,$ coefficients.) As shown by 
Massey \cite[Theorem~1]{massey}, the Pontryagin square 
$\,[\hskip1.2pt\mathcal{U}\,\text{\rm mod}\hskip4pt2\hskip1pt]^2$ is the 
image, under the Thom isomorphism, of 
$\,[\hs\eu(\nu)\,\text{\rm mod}\hskip4pt4\hs]\,
+\,2\left[w_1\smallsmile w_{n-1}\right]$, where 
$\,[\hs\eu(\nu)\,\text{\rm mod}\hskip4pt4\,]\,$ denotes the
$\,\hs\text{\rm mod}\hskip4pt4\,$ reduction of the twisted Euler class of
$\,\nu$. Now (\ref{msf}) follows since naturality of the Pontryagin square 
allows us to compute the right-hand side of (\ref{msf}) in the cohomology of 
either $\,\text{\rm Th}\hskip.6pt(\nu)\,$ or $\,\y$, with both sides treated 
as elements of $\,\bbZ_4$.
\end{proof}
\begin{cor}\label{sigsq}Let $\,f:\x\to\y\,$ be a totally real embedding 
of a \compact\ \rmf\/ $\,\x\,$ in a \compact\ almost complex 
manifold $\,\y\,$ with $\,\dimr\x=\dimc\y=n$, where $\,n\,$ is 
even, and let $\,\sgm\in H^n(\y,\bbZ_2)\,$ be the Poincar\'e dual of the 
homology class $\,f_*[\x\hs]\in H_n(\y,\bbZ_2)$. The Pontryagin square 
$\,\sgm^2\in\bbZ_4$ then is characterized by\/ {\rm(\ref{msf})}. For\/ 
$\,n=2$, this becomes
\begin{equation}\label{psq}
\sgm^2\,=\,\,[\hs\chi(\x)\hskip2pt\text{\rm mod}\hskip4pt4\hs]\,.
\end{equation}
Equality {\rm(\ref{psq})} remains valid, with $\,\sgm^2=0$, if\/ $\,\y$, 
rather than being compact, is assumed \feic\ to $\,\bbC^2\nnh$.
\end{cor}
This is obvious from (\ref{tnu}) and Wu's formula (\ref{wus}); if 
$\,\y\approx\bbC^2\nnh$, we may use Lemma~\ref{masle} with $\,\y\,$ replaced 
by its one-point compactification $\,S^4\nnh$.

A further consequence is 
\begin{cor}\label{chmod}Suppose that a \compact\ \rsu\/ $\,\x\,$ 
admits a totally real embedding in the \csu\/ $\,\y\,=\,\cp^2$ or\/ 
$\,\y\,=\,\cp^2\,\#\,\,\overline{\cp^2}$. Defining\/ 
$\,\chi_4\in\{0,1,2,3\}\,$ by\/ 
$\,\chi_4\equiv\chi(\x)\hskip6pt\text{\rm mod}\hskip4pt4$, we then have\/ 
$\,\chi_4\in\{0,1\}\,$ if\/ $\,\y\,=\,\cp^2$ and\/ 
$\,\chi_4\in\{0,1,3\}\,$ if\/ $\,\y\,=\,\cp^2\,\#\,\,\overline{\cp^2}$. 
\end{cor}
In fact, as in the argument for Proposition~\ref{trchi}, this is immediate 
from (\ref{psq}), since, for integers $\,p,q$, the 
$\,\hs\text{\rm mod}\hskip4pt4\,$ congruence class of $\,p^2$ must contain 
$\,0\,$ or $\,1$, while that of $\,p^2\nnh-q^2$ must contain $\,0$, $\,1\,$ or 
$\,-1$.
\begin{proof}[Proofs of Corollaries~\ref{codod} --~\ref{coqtd}] The 
nonexistence assertion for orientable $\,\x\,$ other than the torus is 
immediate from Corollary~\ref{tonly}, while the case of $\,T^2$ is covered by 
Corollary~\ref{coqtr}. Now, let $\,\x\,$ be nonorientable. The operation 
$\,\x\,\mapsto\,\x\,\#\,T^2\#\,K^2$ (with $\,K^2$ denoting the \kb) 
reduces $\,\chi(\x)\,$ by $\,4\hskip.6pt\text{\rm;}$ thus, by 
Theorem~\ref{thdue}, for the {\it existence} assertions in 
Corollaries~\ref{codod} -- \ref{coqtd} it suffices to show that a \tre\ exists 
if $\,\chi(\x)\,$ is zero, or $\,\chi(\x)\in\{0,1\}\,$ or, respectively, 
$\,\chi(\x)\in\{-1,0,1\}$. For $\,\chi(\x)=0\,$ in all three corollaries, 
or $\,\chi(\x)=1\,$ in Corollary~\ref{cotrd}, or $\,\chi(\x)=\pm1\,$ in 
Corollary~\ref{coqtd}, this existence statement is clear from 
Corollary~\ref{coqtr}, or (vi) in Section~\ref{se} or, respectively, (vi) in 
Section~\ref{se} combined with Corollary~\ref{trerp}, as $\,\chi(\x)=-1$ for 
$\,\x=\,\rp^2\#\,\rp^2\#\,\rp^2\nnh$.

Finally, the nonexistence statements for nonorientable $\,\x\,$ can be 
established as follows. Let $\,\x\,$ admit a totally real embedding in 
$\,\y$. Assume first that $\,\y\,=\,\cp^2$ or 
$\,\y\,=\,\cp^2\,\#\,\,\overline{\cp^2}$. By Corollary~\ref{chmod}, 
$\,\chi(\x)$ has the required remainder $\,\,\text{\rm mod}\hskip4pt4$. Next, 
if $\,\y=\bbC^2\nnh$, this last fact yields 
$\,\chi(\x)\equiv0\hskip6pt\text{\rm mod}\hskip4pt4$, since 
$\,\bbC^2\subset\cp^2$ and $\,\chi(\x)\,$ is even by Corollary~\ref{conov}.
\end{proof}

\section{The sets $\,\hs\fri_q(\x)\,$ and $\,\hs\depspm(\y)\,$ in special 
cases}\label{ts}
\setcounter{equation}{0}
Given a \rmf\ $\,\x\,$ and $\,q\in\{2,4,\dots,\infty\}$, let 
$\,\hs\fri_q(\x)\,$ be the set defined by (\ref{iqs}), with 
$\,\bbZ_\infty=\bbZ$. Thus, as $\,H^1(S^n,\bbZ_q)=\{0\}\,$ for $\,n\ge2$,
\begin{equation}\label{isn}
\fri_q(S^n)\,=\,\{0\}\qquad\text{\rm whenever}\quad n\ge2\quad\text{\rm and}
\quad q\in\{2,4,6,\dots,\infty\}\hs.
\end{equation}
Since $\,\fri_q(\x)\,$ is either empty or forms a coset of the subgroup of 
$\,H^1(\x,\bbZ_q)=\,\h(\pi_1\x,\bbZ_q)\,$ consisting of all homomorphisms 
valued in the even subgroup $\,2\bbZ_q$ (the image of $\,\bbZ_q$ under the 
homomorphism $\,\xi\mapsto2\xi$), we also have
\begin{equation}\label{itn}
\fri_q(T^n)\,=\,(2\bbZ_q)^n\,\subset\,(\bbZ_q)^n\,=\,H^1(T^n\nnh,\bbZ_q)\hs,
\end{equation}
$(\,)^n$ in $\,(2\bbZ_q)^n$, $\,(\bbZ_q)^n$ being the $\,n$th Cartesian power, 
with $\,H^1(T^n\nnh,\bbZ_q)\,=(\bbZ_q)^n$ due to the standard identification 
resulting from (\ref{coh}) (where $\,\pi_1T^n=\bbZ^n$). Next, one easily 
verifies that, for all $\,q\in\{2,4,6,\dots,\infty\}$, 
\begin{equation}\label{irp}
\fri_q(\rp^n)\,=\,\begin{cases}{\{\varphi\}\,\,}&\text{\rm if $\,q/2\,$ is 
finite and odd}$$\hskip.4pt,\cr
{\,\emp\,}&\text{\rm if $\,q=\infty\,$ or $\,q/2\,$ is even}$$\hskip.4pt,
\cr\end{cases}
\end{equation}
where $\,\varphi\in H^1(\rp^n,\bbZ_q)=\,\h(\bbZ_2,\bbZ_q)\,$ (cf.\ 
(\ref{coh})) is the unique nonzero homomorphism $\,\bbZ_2\to\bbZ_q$. (Note 
that $\,\pi_1[\rp^n]=\bbZ_2$, and $\,w_1(\rp^n)\,$ with (\ref{ori}) equals 
$\,\text{\rm Id}:\bbZ_2\to\bbZ_2$.)

For an \acm\ $\,\y$, a fixed sign $\,\pm\hs$, and $\,\ve\in\{0,1\}$, let 
$\,\hs\depspm(\y)\,$ now be defined by (\ref{dmk}). If $\,\y\,=\,\cp^2$, we 
thus have
\begin{equation}\label{dcp}
\donem(\cp^2)\nnh=\nh\emp\hs,\hskip4.5pt
\donep(\cp^2)\nnh=\dzerp(\cp^2)\nnh=\nnh\{0\},\hskip4.5pt
\dzerm(\cp^2)\nnh=\nnh\{[\rp^2]\},\hskip-7pt
\end{equation}
$[\rp^2]\,$ being the $\,\bbZ_2$-homology class of $\,\x=\rp^n$ in (vi) of 
Section~\ref{se} for $\,n=2$. (In fact, $\,c_{\hs1}(\y):H_2(\y,\bbZ)\to\bbZ\,$ 
and $\,w_2(\y):H_2(\y,\bbZ_2)\to\bbZ_2$ then are {\it isomorphisms}.) Finally, 
since $\,H_2(\bbC^2,\bbZ_{[2]})=0$,
\begin{equation}\label{dco}
\donep(\bbC^2)\,=\,\{0\}\,,\hskip10pt
\dzerp(\bbC^2)\,=\,\{0\}\,,\hskip10pt
\dzerm(\bbC^2)\,=\,\,\donem(\bbC^2)\,=\,\emp\,.
\end{equation}

\section{Closed surfaces and cohomology}\label{ch}
\setcounter{equation}{0}
We begin with some known facts gathered here for easy reference. First, the 
fundamental group $\,\Gamma=\pi_1K^2$ of the \kb\ 
$\,K^2=\,\rp^2\#\,\rp^2\nnh$, treated as a group of deck transformations in 
$\,\rto\nnh$, has the generators $\,\varPhi,\varPsi\,$ described in 
Example~\ref{exbdd}. Since 
$\,\varPsi^2=\varPsi\varPhi\varPsi^{-1}\varPhi^{-1}\nh$, the Abelianization 
$\,H_1(K^2\nnh,\bbZ)\,$ of $\,\Gamma\,$ may from now on be identified with the 
direct product $\,\bbZ\times\bbZ_2$ (the factor groups being generated by 
$\,\varPhi\,$ and $\,\varPsi$). As $\,\hs\varPsi\nh\,$ is 
o\-ri\-en\-ta\-tion-pre\-serv\-ing while $\,\varPhi\,$ is not, 
$\,w_1(K^2):\pi_1K^2\to\bbZ_2=\{0,1\}\,$ sends $\,\varPhi\,$ to $\,1\,$ and 
$\,\varPsi\,$ to $\,0\,$ (cf.\ (\ref{ori})). The homomorphism 
$\,\bbZ\times\bbZ_2\to\bbZ_2$ arising as $\,w_1(K^2)\,$ descends to 
$\,H_1(K^2\nnh,\bbZ)\,$ therefore sends $\,(\ek,\ve)\,$ to 
$\,[\hs\ek\hskip2.5pt\text{\rm mod}\hskip4pt2\hskip1pt]$. Finally, the 
transformation $\,(t,s)\mapsto\,(-t,s)\,$ in $\,\rto$ commutes with 
$\,\varPsi\,$ and conjugates $\,\varPhi\,$ with $\,\varPhi^{-1}\nnh$, so that 
it descends to a diffemorphism $\,K^2\to K^2$ of $\,K^2=\rto\nnh/\hs\Gamma$, 
which acts in $\,H_1(K^2,\bbZ)=\bbZ\times\bbZ_2$ by 
$\,(\ek,\ve)\mapsto(-\hs\ek,\ve)$.
\begin{rem}\label{srfcs}For a \compact\ \rsu\ $\,\x\,$ of genus $\,s\,$ we 
have, up to a \feo, one of three cases:
\begin{enumerate}
  \def\theenumi{{\rm\alph{enumi}}}
\item $\x\,=\,\,T^2\#\,\ldots\,\#\,T^2$ ($s\,$ summands),
\item $s\,$ is even and $\,\x=T^2\#\,\ldots\,\#\,T^2\#\,K^2$, with $\,s/2\,$ 
summands,
\item $s\,$ is odd and $\,\x=T^2\#\,\ldots\,\#\,T^2\#\,\rp^2$, with 
$\,(s+1)/2\,$ summands.
\end{enumerate}
This includes $\,\x=S^2$ (case (a), with $\,s=0$), the \kb\ $\,\x=K^2$ (case 
(b) with $\,s=2$), and $\,\x=\rp^2$ (case (c), for $\,s=1$).
\end{rem}
\begin{rem}\label{abgps}Let $\,G,G_1,\dots,G_\ek$ be Abelian groups and let 
$\,\x\,$ be a manifold with $\,H_1(\x,\bbZ)=G_1\times\ldots\times G_\ek$. Then, 
by (\ref{coh}) and Remark~\ref{reduc},
\begin{enumerate}
  \def\theenumi{{\rm\roman{enumi}}}
\item $H^1(\x,G)=\,\h(G_1,G)\times\ldots\times\,\h(G_\ek,G)$,
\item the coefficient-reduction homomorphism $\,H^1(\x,G)\to H^1(\x,G')\,$ 
corresponding to any homomorphism $\,h:G\to G'$ of Abelian groups is the 
Cartesian product of the composition homomorphisms for the factor groups, with 
both $\,H^1(\x,G)\,$ and $\,H^1(\x,G')\,$ decomposed as in (i).
\end{enumerate} 
\end{rem}
Given an Abelian group $\,G$, we denote by $\,G\ortwo$ the subgroup of $\,G\,$ 
consisting of zero and all elements of order $\,2$. Thus,
\begin{equation}\label{hkb}
H^1(K^2,G)\,=\hs\,G\times G\ortwo\hs,
\end{equation}
as one sees combining Remark~\ref{abgps}(i) with the relation 
$\,H_1(K^2,\bbZ)=\bbZ\times\bbZ_2$ at the beginning of this section. Obviously,
\begin{equation}\label{ord}
\text{\rm$(\bbZ_q){}\ortwo=\{0,q/2\}\,$\ for\ any\ even\ positive\ integer\ 
$\,q$\hs.}
\end{equation}
Also, given \compact\ \rsu s $\x,\x\hs'$ such that $\,\x\hs'$ is orientable,
\begin{equation}\label{hss}
\aligned
\text{\rm i)}&\hskip5pt
H_1(\x\,\#\,\x\hs',\bbZ)\,=\,H_1(\x,\bbZ)\times H_1(\x\hs',\bbZ)\hs,\\
\text{\rm ii)}&\hskip5pt
H^1(\x\,\#\,\x\hs'\hskip-1.5pt,G)
=H^1(\x,G)\times H^1(\x\hs'\hskip-1.5pt,G)\hs,\\
\text{\rm iii)}&\hskip5pt
\fri_q(\x\#\x\hs')=\,\fri_q(\x)\times\hs\fri_q(\x\hs')\hskip10pt\text{\rm for}
\hskip8ptq\in\{2,4,6,\dots,\infty\}\hs,\endaligned
\end{equation}
where $\,\fri_q(\x)\,$ is defined by (\ref{iqs}), with $\,\bbZ_\infty=\bbZ$. 
Namely, the well-known isomorphic identification (\ref{hss}.i) combined with 
Remark~\ref{abgps}(i) and (\ref{coh}) gives (\ref{hss}.ii). Next, 
$\,w_1(\x\,\#\,\x\hs')=(w_1(\x),w_1(\x\hs'))\,$ (in terms of (\ref{hss}.ii) 
with $\,G=\bbZ_2$), since, by (\ref{ori}), $\,w_1(\x\,\#\,\x\hs')$, acting on 
$\,H_1(\,\cdot\,,\bbZ)\,$ rather than $\,\pi_1(\,\cdot\,)$, coincides with the 
first-factor projection in (\ref{hss}.i) followed by $\,w_1(\x)$. Using 
Remark~\ref{abgps}(ii) and (\ref{iqs}), we now get (\ref{hss}.iii).

On the other hand, for the \kb\ $\,K^2$ and $\,q\in\{2,4,6,\dots,\infty\}$,
\begin{equation}\label{ikb}
\fri_q(K^2)\,=\,(\bbZ_q\smallsetminus2\bbZ_q)\times
[(\bbZ_q){}\ortwo\cap2\bbZ_q]\,.
\end{equation}
In fact, combining the description of $\,w_1(K^2)\,$ at the beginning of this 
section with (\ref{hkb}) for $\,G=\bbZ_2$ we get 
$\,w_1(K^2)=(1,0)\in\bbZ_2\times\bbZ_2$ (as $\,(\bbZ_2){}\ortwo=\bbZ_2$). Now 
(\ref{ikb}) is immediate from Remark~\ref{abgps}(ii) and (\ref{hkb}).

Given a \compact\ manifold $\,\x\,$ and an even positive 
integer $\,q$, the $\,\hs\text{\rm mod}\hskip4ptq$ reduction homomorphism 
$\,H^1(\x,\bbZ)\to H^1(\x,\bbZ_q)\,$ sends the set 
$\,\hs\fri_\infty(\x)$ into $\,\hs\fri_q(\x)$. (Notation of (\ref{iqs}), 
with $\,\bbZ_\infty=\bbZ$.) This is clear from Remark~\ref{reduc}: as $\,q\,$ 
is even, reduction $\,\hs\text{\rm mod}\hskip4ptq\,$ followed by reduction 
$\,\hs\text{\rm mod}\hskip4pt2$ results in reduction 
$\,\hs\text{\rm mod}\hskip4pt2$. We thus obtain the 
$\,\hs\text{\rm mod}\hskip4ptq\,$ {\it reduction mapping} 
\begin{equation}\label{inf}
\fri_\infty(\x)\mathop{\,-\!\!\!\!-\!\!\!-\!\!\!-\!\!\!\!\longrightarrow\,}
\limits^{\text{\rm mod}\hskip2ptq}\fri_q(\x)\,.
\end{equation}
\begin{lem}\label{surfq}Given a \compact\ \rsu\/ $\,\x\,$ and an even 
integer\/ $\,q>0$,
\begin{enumerate}
  \def\theenumi{{\rm\roman{enumi}}}
\item the mapping\/ {\rm(\ref{inf})} is surjective \iff\ 
either\/ $\,\x\,$ is orientable, or\/ $\,\x\,$ is nonorientable and\/ 
$\,\chi(\x)-q/2\,$ is odd,
\item if\/ $\,q=4\,$ and\/ $\,\x\,$ is the \kb\ $\,K^2$, 
the image of\/ {\rm(\ref{inf})} is the $\,2$-el\-e\-ment subset\/ 
$\,\{1,3\}\times\{0\}\,$ of the $\,4$-el\-e\-ment set\/ 
\begin{equation}\label{ifr}
\fri_4(K^2)\,=\,\{1,3\}\times\{0,2\}\,,
\end{equation}
which itself is a subset of\/ 
$\,H^1(K^2,\bbZ_4)=\{0,1,2,3\}\times\{0,2\}\subset\bbZ_4\times\bbZ_4$, cf.\ 
{\rm(\ref{hkb})}, with\/ $\,\bbZ_4=\{0,1,2,3\}$. 
\end{enumerate}
\end{lem}
\begin{proof}Succesively applying (\ref{hss}.iii) we obtain, from (\ref{isn}) 
-- (\ref{irp}) and (\ref{ikb}), with superscripts denoting 
Cartesian powers, $\,\fri_q(\x)=(2\bbZ_q){}^{2g}\,$ in case (a) of 
Remark~\ref{srfcs}, 
$\,\fri_q(\x)=(2\bbZ_q){}^{g-2}\times(\bbZ_q\smallsetminus2\bbZ_q)\times[
(\bbZ_q){}\ortwo\cap2\bbZ_q]\,$ in case (b), and 
$\,\fri_q(\x)=(2\bbZ_q){}^{g-1}\times[(\bbZ_q){}\ortwo
\smallsetminus2\bbZ_q]\,$  in case (c) (as 
$\,\fri_q(\rp^n)=(\bbZ_q){}\ortwo\smallsetminus2\bbZ_q$ for 
$\,q\in\{2,4,6,\dots,\infty\}$, cf.\ (\ref{irp}) and (\ref{ord})). In view of 
to Remark~\ref{abgps}(ii), 
the mapping (\ref{inf}) acts by $\,\hs\text{\rm mod}\hskip4ptq\,$ reduction in 
each Cartesian factor set just listed. Also, since $\,q\,$ is even, 
the reduction mappings $\,\bbZ\to\bbZ_q$ and 
$\,\bbZ\smallsetminus2\bbZ\to\bbZ_q\smallsetminus2\bbZ_q$ are surjective. 
Therefore, surjectivity of (\ref{inf}) always holds in case (a) of 
Remark~\ref{srfcs}, while in case (b) or (c) it is equivalent to 
$\,(\bbZ_q){}\ortwo\cap\hs2\bbZ_q=\{0\}\,$ or, respectively, 
$\,(\bbZ_q){}\ortwo\smallsetminus2\bbZ_q=\,\emp$. Now (i) is immediate from 
(\ref{ord}), as $\,\chi(\x)=2-g\,$ in cases (b) and (c). Finally, (\ref{hkb}) 
-- (\ref{ord}) give 
$\,H^1(K^2,\bbZ_4)=\{0,1,2,3\}\times\{0,2\}$, so that (\ref{ifr}) follows from 
(\ref{ikb}) for $\,q=4$, and (ii) is immediate as 
$\,\bbZ{}\ortwo\cap\hs2\bbZ=\{0\}$.
\end{proof} 

\section{More on tori and Klein bottles in $\,\bbC^2$}\label{mt}
\setcounter{equation}{0}
We will now evaluate the Maslov index $\,\hs\ind\hs\,$ and degree 
$\,\hs\dg\hs\,$ for the \tri s $\,f^{\ek,\hs l}$ constructed in 
Example~\ref{trtkb}.

When $\,\x\,$ is the torus $\,T^2$ or the \kb\ $\,K^2$, (\ref{itn}) and 
(\ref{ikb}) with $\,n=2\,$ and $\,q=\infty\,$ become 
$\,\,\fri_\infty(T^2)=2\bbZ\times2\bbZ\subset\bbZ\times\bbZ=H^1(T^2,\bbZ)\,$ 
and 
$\,\,\fri_\infty(K^2)=(\bbZ\smallsetminus2\bbZ)\times\{0\}\subset
\bbZ\times\{0\}=H^1(K^2,\bbZ)\,$ (cf.\ (\ref{hkb}), (\ref{ord})). 

A \tri\ $\,(x,u):\rto\to\bbC^2$ which is doubly $\,2\pi$-periodic (or, in 
addition, also invariant under the transformation $\,\varPhi\,$ of 
Example~\ref{exbdd}) descends to a \tri\ $\,f:\x\to\bbC^2$ with $\,\x=T^2$ 
(or, respectively, $\,\x=K^2$). Let us now set, for 
$\,\phi=\mathcal{J}(x,u)\,$ defined as in Example~\ref{trsub},
\begin{equation}\label{tkb}
p=\hs\frac1{2\pi i}\,\int_0^{2\pi}\nh\frac{\phi_s}{\phi}\,ds\hs,\hskip7pt
\rx=\hs\frac1{2\pi i}\,\int_0^{2\pi}\nh\frac{\phi_t}{\phi}\,dt\hs,\hskip7pt
m=\hs\frac1{\pi i}\,\int_0^\pi\nnh\left[\frac{\phi_t}{\phi}\right]
\underset{{\hskip-2.6pts\hs=\hs0}}{\phantom{_jjj}}dt\hs,\hskip-2pt
\end{equation}
where the subscripts represent partial derivatives. Thus, $\,p,\rx,m\in\bbZ\,$ 
and the Maslov index of $\,f\,$ (see Section~\ref{sr}) is
\begin{equation}\label{ind}
\ind(f)\,=\,\begin{cases}{(2p,2\rx)\,\in\,2\bbZ\times2\bbZ\,
=\,\,\fri_\infty(\x)\,\,}
&\text{\rm if}$$\quad\x\,=\,T^2\nnh,\cr
{(m,0)\,\in\,(\bbZ\smallsetminus2\bbZ)\times\{0\}\,
=\,\fri_\infty(K^2)\,\,}&\text{\rm if}$$\quad\x\,=\,K^2\nnh.\cr
\end{cases}
\end{equation}
\begin{prop}\label{tkind}Let\/ $\,\x\,$ stand for the\/ $\,2$-torus\/ $\,T^2$ 
or the \kb\ $\,K^2\nh$, and let\/ $\,f:\x\to\bbC^2$ be a totally real 
immersion obtained as in Example~\ref{exbdd} from some\/ $\,x(t,s),\,y(t,s)\,$ 
and\/ $\,h(t)\,$ satisfying the hypotheses of Example~\ref{exbdd}. Then\/ 
$\,f\,$ has the Maslov index\/ {\rm(\ref{ind})} with\/ 
$\,p=(2\pi i)^{-1}\int_0^{2\pi}(x_{ss}/x_s)\,ds$, 
$\,\rx=(2\pi i)^{-1}\int_0^{2\pi}[(x_{st}/x_s)+(h_{tt}/h_t)]\,dt$, and\/ 
$\,m=(\pi i)^{-1}\nnh\int_0^\pi[(x_{st}/x_s)\nh+\nh(h_{tt}/h_t)]_{s\hs=\hs0}\,
dt$, the subscripts being partial derivatives.
\end{prop}
In fact, (\ref{tkb}) holds for $\,\phi=\phi^{[\ax]}$ with 
$\,\phi^{[\ax]}=\ax^{-1}\,\mathcal{J}(x,y+\ax h)\,$ and large $\,\ax>0\,$ 
(notation of Example~\ref{trsub}). Homotopy invariance of the degree then 
guarantees that the integrals (\ref{tkb}) do not depend on $\,\ax$. Taking 
their limits as $\,\ax\to\infty$, we obtain our assertion, since 
$\,\lim_{\,\ax\to\infty}\phi^{[\ax]}=x_sh_t$ by (\ref{jrh}).
\begin{prop}\label{tkbtr}There exist totally real embeddings of the\/ 
$\,2$-torus\/ $\,T^2$ and the \kb\ $\,K^2$ in\/ $\,\bbC^2$ which 
realize any prescribed Maslov index\/ $\,\hs\ind(f)\,$ with\/ {\rm(\ref{ind})}.
\end{prop}
\begin{proof}To exhibit such embeddings, we first apply 
Proposition~\ref{tkind} to the \tre s $\,f^{\ek,1}:T^2\to\bbC^2$ for 
$\,\ek\in\bbZ\,$ and $\,f^{\ek,2}:K^2\to\bbC^2$ for odd $\,\ek$, described in 
Example~\ref{trtkb}. The mapping $\,F:\bbR/2\pi\bbZ\,\to\,S^1$ defined by the 
assignment $\,\bbR\ni s\mapsto x_s/|x_s|\,$ (with any fixed $\,t$) now is of 
degree zero (since $\,x_s(t,s+\pi)\,=\,-\,\overline{x_s(t,s)}$, and so $\,F\,$ 
is homotopic to its composite with the conjugation). Therefore, by 
Proposition~\ref{tkind}, we have $\,p=0$, that is, $\,f^{\ek,1}$ and $\,f^{\ek,2}$ 
have the Maslov indices $\,(2p,2\rx)=(0,2\ek+2)\,$ and, respectively, 
$\,(m,0)=(\ek+2,0)$. We thus have 
realized all index values for the \kb; to obtain an arbitrary Maslov 
index $\,(2p,2\rx)\,$ for the torus, we set $\,\ek=-1\,$ (when $\,p=\rx=0$) or, 
when $\,(p,\rx)\ne(0,\nh0)$, use the composite of our $\,f$, for $\,\ek=d-1$, with 
the group automorphism of $\,T^2$ corresponding to a suitable matrix 
$\,\mathfrak{A}\in\,\text{\rm SL}(2,\bbZ)$. Specifically, we fix 
$\,b,c\in\bbZ\,$ such that $\,b\hs\rx-c\hs p=d$, where $\,d\,$ is the greatest 
common factor of $\,p\,$ and $\,\rx$, and then choose $\,\mathfrak{A}\,$ to 
have the rows $\,[\hs b\hskip8ptc\hs]\,$ and $\,[\hs p/d\hskip10pt\rx/d\hs]$.
\end{proof}

\section{The integer $\hs\,q\hs\,$ in \,{\rm(\ref{pzq})}}\label{tq}
For a \sc\ \acm\ $\,\y$, let $\,q\,$ be as in (\ref{pzq}).
\begin{enumerate}
  \def\theenumi{{\rm\alph{enumi}}}
\item $w_2(\y)\,$ is the $\,\hs\text{\rm mod}\hskip4pt2\,$ reduction of 
$\,c_{\hs1}(\y)$.
\item The connecting homomorphism in (\ref{seq}) for $\,E=E^+(\y)\,$ is, up to 
a sign, the composite $\,\pi_2\y\nnh\to\nnh H_2(\y,\bbZ)\nnh\to\nnh\bbZ\,$ of 
the Hurewicz isomorphism with $\,c_{\hs1}(\y)$, and so its image, that is, the 
kernel of $\,\bbZ\to\pi_1E\,$ in (\ref{seq}), is trivial if $\,q=\infty$, and 
generated by $\,q/2\,$ if $\,q<\infty$.
\item $\y\,$ is spin if and only if $\,q=\infty\,$ or $\,q\,$ is finite and 
divisible by $\,4$. Similarly, $\,q=\infty\,$ if and only if 
$\,c_{\hs1}(\y)=0$.
\item If $\,\y\,$ is not spin, or $\,c_1(M)=0$, the 
$\,\hs\text{\rm mod}\hskip4pt2\,$ reduction homomorphism 
$\,H_2(\y,\bbZ)\to H_2(\y,\bbZ_2)\,$ maps 
$\,\hs\text{\rm Ker}\,[c_{\hs1}(\y)]\,$ {\it onto\/} 
$\,\hs\text{\rm Ker}\,[w_2(\y)]$.
\item If $\,\y=\bbC^2\nnh$, or $\,\y=\cp^2\nnh$, or 
$\,\y\,=\,\cp^1\!\times\cp^1\nnh$, then $\,q=\infty$, or $\,q=6$, or,
respectively, $\,q=4$.
\item $q=2\,$ whenever $\,\y\,$ is obtained by blowing up a point in a \sc\ 
\csu.
\end{enumerate}
In fact, (a) and (b) are well known, while (c) follows from (a), (b) and 
Remark~\ref{reduc}. As $\,H_2(\y,\bbZ)\to H_2(\y,\bbZ_2)\,$ is surjective, 
(a) yields (d) in the case where $\,c_1(M)=0\,$ (cf.\ Remark~\ref{reduc}). 
Again by (a), if $\,\y\,$ is not spin, $\,c_{\hs1}(\y):H_2(\y,\bbZ)\to\bbZ\,$ 
must assume some odd values; the homomorphism sending 
$\,H_2(\y,\bbZ)/\text{\rm Ker}\,[c_{\hs1}(\y)]\,\approx\,\bbZ\,$ into 
$\,H_2(\y,\bbZ_2)/\text{\rm Ker}\,[w_2(\y)]\,\approx\,\bbZ_2$, induced by the 
$\,\hs\text{\rm mod}\hskip4pt2\,$ reduction, is therefore surjective. Every 
element of $\,\hs\text{\rm Ker}\,[w_2(\y)]\,$ thus has a pre\-im\-age in 
$\,H_2(\y,\bbZ)\,$ with an even image in 
$\,H_2(\y,\bbZ)/\text{\rm Ker}\,[c_{\hs1}(\y)]$, which proves (d). Next, 
(a) implies (e), as the values assigned by $\,c_{\hs1}(\y)\,$ to the standard 
generator(s) of $\,H_2(\y,\bbZ)\,$ for $\,\y=\cp^2$ or 
$\,\y\,=\,\cp^1\!\times\cp^1$ are $\,3\,$ or, respectively, $\,2\,$ and $\,2$. 
Finally, (f) follows from (a), since $\,c_{\hs1}(\y)\,$ yields the value 
$\,1\,$ when integrated over the exceptional divisor resulting from the 
blow-up (cf.\ also formula (\ref{con})).

\section{Index and degree after modifications}\label{am}
\setcounter{equation}{0}
We now determine how the zooming and con\-nect\-ed-sum operations of 
Sections~\ref{zg} and~\ref{cs} affect the Maslov index and degree.
\begin{lem}\label{indeg}Given a totally real im\-mer\-sion/em\-bed\-ding\/ 
$\,f\,$ of a \compact\ \rsu\/ $\,\x\,$ in\/ $\,\bbC^2\nnh$, and a \sc\ \acsu\/ 
$\,\y$, let a totally real im\-mer\-sion/em\-bed\-ding\/ $\,f':\x\to\y\,$ be 
obtained from\/ $\,f\,$ by a zooming procedure described in\/ {\rm 
Proposition~\ref{zoopr}}. Its Maslov index and degree, introduced in\/ 
Section~{\rm\ref{sr}}, then are
\begin{equation}\label{ifp}
\ind(f')\,=\,[\hs\ind(f)\,\text{\rm mod}\hskip4ptq\hskip1pt]\,\in\,
\fri_q(\x)\,\subset\,H^1(\x,\bbZ_q),\qquad\dg(f')\,=\,0\,,
\end{equation}
where\/ $\,q\,$ is determined by\/ $\,\y\hs$ via\/ {\rm(\ref{pzq})} and\/ 
$\,[\hs\ind(f)\,\text{\rm mod}\hskip4ptq\hskip1pt]\,$ denotes the 
$\,\hs\text{\rm mod}\hskip4ptq\,$ reduction of\/ $\,\ind(f)\in H^1(\x,\bbZ)$, 
that is, its image under the mapping\/ {\rm(\ref{inf})}.
\end{lem}
\begin{proof}Let $\,\,U\subset\bbC^2$ be a ball containing $\,f(\x)$. Since 
$\,\,U\hs$ is contractible, $\,\dg(f')=0$. Also, $\,\hs\ind(f')\,$ is the 
homomorphism (\ref{hfa}) of fundamental groups induced by 
$\,\mathfrak{M}(f')$. Thus, $\,\hs\ind(f')\,$ equals $\,\hs\ind(f)\,$ followed 
by the homomorphism $\,\pi_1[E(U)]\,\to\,\pi_1[E(\y)]\,$ induced by an 
embedding of $\,\,U\to\y$. As the fibre inclusion 
$\,\hs\text{\rm U}\hs(1)\to E(\y)\,$ represents a generator of 
$\,\pi_1[E(\y)$, the latter homomorphism is nothing else than 
the projection $\,\bbZ\to\bbZ_q$.
\end{proof} 
Lemma~\ref{indeg} can be modified when $\,\y\,$ is a {\it complex\/} surface; 
rather than by zooming, $\,f'$ may then be obtained as the composite of 
$\,f\,$ with a hol\-o\-mor\-phic embedding in $\,\y\,$ of an open ball in $\,\bbC^2$ 
containing $\,f(\x)$. 
\begin{cor}\label{conas}Given a \sc\ \acsu\ $\hs\y$, let us consider the 
following condition imposed on a \compact\ \rsu\ $\,\x\,$ and an element\/ 
$\,\hs\ind\in\fri_q(\x)$, with\/ $\,q\,$ and\/ $\,\hs\fri_q(\x)\,$ defined by 
{\rm(\ref{pzq})} and {\rm(\ref{iqs})}{\rm:} 
\begin{enumerate}
  \def\theenumi{{\rm\roman{enumi}}}
\item[($*$)] {\rm some\ totally\ real\ embedding\ $\,f:\x\to\y\,$\ has\ 
$\,\hs\ind(f)=\hs\ind\,\,$\ and\ $\,\hs\dg(f)=\hs0$.}
\end{enumerate}
Condition\/ {\rm($*$)} is satisfied by 
\begin{enumerate}
  \def\theenumi{{\rm\alph{enumi}}}
\item the torus\/ $\,\x=T^2$ and every $\,\hs\ind\in\,\fri_q(\x)$, 
\item the \kb\ $\,\x=K^2$ and every $\,\hs\ind\in\,\fri_q(\x)$, 
provided that\/ $\,q/2\,$ is either infinite, or finite and odd, 
\item the \kb\ $\,\x=K^2$ and every $\,\hs\ind\hs\,$ in the subset\/ 
$\,\{1,3\}\times\{0\}\,$ of\/ $\,\hs\fri_4(\x)$, cf.\ {\rm(\ref{ifr})}, 
provided that\/ $\,q=4$.
\end{enumerate}
\end{cor}
This is immediate from Lemma~\ref{indeg}, Proposition~\ref{tkbtr} and 
Lemma~\ref{surfq}.
\begin{lem}\label{iifdf}Let\/ $\,f,f'$ be \tri s of \compact\ \rsu s 
$\,\x,\x\hs'$ in a \sc\ \acsu\ $\,\y$. If\/ $\,\x\hs'$ is orientable, then a 
\tri\ $\,\x\,\#\,\x\hs'\to\y\,$ of the connected sum, obtained from $\,f\,$ 
and\/ $\,f'$ as in Theorem~\ref{cnsum}, has the Maslov index and degree
\begin{equation}\label{iff}
\ind\,\,=\,(\ind(f),\,\ind(f'))\,,\qquad
\dg\,=\,\hs\dg(f)\,+\,\hs\dg(f')\,\in\,H_2(\y,\bbZ_{[2]})\,.
\end{equation}
We use here the notational conventions of\/ {\rm(\ref{hss}.iii)} and\/ 
{\rm(\ref{ztw})}, while adding an element of\/ $\,H_2(\y,\bbZ)\,$ to an 
element of\/ $\,H_2(\y,\bbZ_2)\,$ is to be preceded by\/ 
$\,\hs\text{\rm mod}\hskip4pt2\,$ reduction of the former, so that the sum 
lies in\/ $\,H_2(\y,\bbZ_2)$.
\end{lem}
In fact, (\ref{iff}) follows from (\ref{hss}.i), since the new immersion 
coincides with $\,f\,$ or $\,f'$ except on the boundary of a solid cylinder 
immersed in $\,\y\,$ (Theorem~\ref{cnsum}(v))).

Next, for $\,\mathcal{Z}(\x,\y)\,$ defined in the lines following 
Theorem~\ref{thuno}, we have
\begin{lem}\label{cldif}Given a \sc\ \acsu\/ 
$\,\y$, let\/ $\,\mathcal{E}$ be the class of \feo\ types of \compact\ 
\rsu s $\,\x\,$ such that for every\/ $\,(\ind,\dg)\in\mathcal{Z}(\x,\y)\,$ 
there exists a \tri\/ $\,f:\x\to\y$ with\/ $\,\hs\ind(f)=\,\ind\hs\,$ and\/ 
$\,\hs\dg(f)=\,\dg$. Then\/ $\,\mathcal{E}$ is closed under the 
con\-nect\-ed-sum operation applied to two surfaces, of which at least one is 
orientable. 
\end{lem}
\begin{proof}Let $\,(\ind,\dg)\in\mathcal{Z}(\x\,\#\,\x\hs'\nh,\y)$, where 
$\,\x,\x\hs'$ belong to $\,\hs\mathcal{E}\nh\,$ and $\,\x\hs'$ is orientable. 
By (\ref{hss}.iii), $\,\hs\ind\,=(\ind_0,\hs\ind')\,$ with 
$\,\hs\ind_0\in\,\fri_q(\x)\,$ and $\,\hs\ind'\in\,\fri_q(\x\hs')$, for 
$\,q\,$ as in (\ref{pzq}), so that $\,(\ind_0,\dg)\in\mathcal{Z}(\x,\y)\,$ due 
to the definition of $\,\mathcal{Z}(\x,\y)$. Similarly, 
$\,(\ind'\nh,0)\in\,\fri_q(\x\hs')\times\{0\}
\subset\fri_q(\x\hs')\times\donep(\y)=\mathcal{Z}(\x\hs'\nh,\y)$. Thus, 
$\,(\ind_0,\dg)\,$ and $\,(\ind'\nh,0)\,$ are realized as the index-degree 
pairs of some \tri s $\,\x\to\y\,$ and $\,\x\hs'\to\y$. Realizability 
of $\,(\ind,\dg)\,$ now follows from Lemma~\ref{iifdf}, completing the proof.
\end{proof}
For a \sc\ \acsu\ $\,\y\nh$, consider the condition
\begin{equation}\label{ife}
\begin{array}{l}
q/2\,\,\text{\rm\ is\ infinite,\ or\ finite\ and\ even,\ or,\ finally,\ }\,\,
\rp^2\text{\rm\ admits}\\
\text{\rm a\ totally\ real\ immersion\ in\ }\hs\y\nh\text{\rm,\ 
while\ }\,q/2\,\text{\rm\ is\ finite\ and\ odd.}
\end{array}
\end{equation}
with $\,q\,$ defined by (\ref{pzq}). Further conditions imposed on $\,\y\,$ 
are:
\begin{equation}\label{thx}
\begin{array}{l}
\text{\rm totally\ real\ immersions\ }\,\,\,S^2\to\y\,\,\text{\rm\ realize\ 
every\ degree\ }\,\,\hs\dg\\
\text{\rm in\ some\ set\ generating\ the\ group\ } 
\,\text{\rm Ker}\,[c_{\hs1}(\y)]\nh\subset\nnh H_2(\y,\bbZ)
\end{array}\hskip5pt
\end{equation}
($\emp\,$ is such a set when $\,\hs\text{\rm Ker}\,[c_{\hs1}(\y)]=\{0\}$), 
and, for a \compact\ \rsu\ $\,\x$,
\begin{equation}\label{evr}
\begin{array}{l}
\text{\rm every\ }\,\hs\ind\in\,\fri_q(\x)\text{\rm, with\ }\,q\,
\text{\rm\ as\ in\ (\ref{pzq}), equals\ the\ Maslov}\\
\text{\rm index }\,\,\ind(f)\,\,\text{\rm of\ some\ totally\ real\ immersion\ }
\,f:\x\to\y.
\end{array}\hskip16pt
\end{equation}
We use here the notation of Section~\ref{sr}.
\begin{lem}\label{cedif}Suppose that\/ $\,\x\,$ is a \compact\ \rsu\ and\/ 
$\,\y\,$ is a \sc\ \acsu\ such that either\/ $\,\y\,$ is not spin, or\/ 
$\,c_{\hs1}(\y)=0$. If conditions\/ {\rm(\ref{thx})} and\/ {\rm(\ref{evr})} 
are satisfied, then\/ $\,\x\,$ belongs to the class\/ $\,\mathcal{E}$ defined 
in Lemma~\ref{cldif}.
\end{lem}
\begin{proof}Let 
$\,(\ind,\dg)\in\mathcal{Z}(\x,\y)=\,\fri_q(\x)\times\depspm(\y)$. (Cf.\ (c) 
in Section~\ref{tq} and the definition of $\,\mathcal{Z}(\x,\y)\,$ in the 
lines following Theorem~\ref{thuno}.) Choose $\,f$ for this $\,\hs\ind\hs\,$ 
as in (\ref{evr}). As $\,\hs\dg,\hs\dg(f)\in\hs\depspm(\y)\,$ 
(Lemma~\ref{ifidf}(b)), the difference $\,\hs\dg\hs'=\hs\dg\hs-\hs\dg(f)\,$ 
lies in $\,\hs\text{\rm Ker}\,[c_{\hs1}(\y)]\,$ (if $\,\x\,$ is orientable) or 
in $\,\hs\text{\rm Ker}\,[w_2(\y)]$ (if $\,\x\,$ is not orientable). Thus, 
even in the nonorientable case, $\,\hs\dg\hs'$ is the 
$\,\hs\text{\rm mod}\hskip4pt2$ reduction of an element of 
$\,\hs\text{\rm Ker}\,[c_{\hs1}(\y)]\,$ (see (d) in Section~\ref{tq}).

Hence, in view of (\ref{thx}) combined with Lemma~\ref{iifdf} and 
Corollary~\ref{sfinm}, $\,(\ind\hs,\dg)\,$ is the index-degree pair of a \tri\ 
$\,\x\to\y$.
\end{proof}
\begin{thm}\label{clase}Let a \sc\ \acsu\/ $\,\y\,$ satisfy conditions\/ 
{\rm(\ref{ife})} and\/ {\rm(\ref{thx})}. If the class\/ $\,\mathcal{E}$ 
defined in Lemma~\ref{cldif} contains the \kb\ $\,\x=K^2\nnh$, then\/ 
$\,\mathcal{E}$ contains every \compact\ \rsu.
\end{thm}
\begin{proof}Every \compact\ surface either is one of 
$\,S^2\nnh,T^2\nnh,\rp^2\nnh,K^2\nh$, or can be obtained by iterated connected 
summation in which all intermediate summands, except (possibly) the last one, 
are orientable. (Cf.\ Remark~\ref{srfcs}.) Thus, by Lemma~\ref{cldif}, it 
suffices to show that $\,\mathcal{E}\,$ contains $\,S^2\nnh,T^2$ and 
$\,\rp^2\nnh$. This will in turn follow from Lemma~\ref{cedif}, once we show 
that condition (\ref{evr}) holds when $\,\x\,$ is any of these three \rsu s.

First, (\ref{evr}) is obvious for $\,\x=\rp^2\nnh$, since, by (\ref{irp}), 
$\,\hs\fri_q(\rp^2)\,$ is a one-el\-e\-ment set when $\,q/2\,$ is finite and 
odd, and empty otherwise. Similarly, (\ref{evr}) for $\,\x=S^2$ is immediate 
since $\,\hs\fri_q(S^2)\,$ always has just one element (see (\ref{isn})), and 
a \tri\ $\,S^2\to\y\,$ exists by Corollary~\ref{sfinm}. Finally, (\ref{evr}) 
for $\,\x=T^2$ is obvious from Corollary~\ref{conas}(a).
\end{proof}

\section{Totally real Klein bottles in $\,\cp^1\!\times\cp^1$}\label{rk}
\setcounter{equation}{0}
Let $\,\y=\cp^1\!\times\cp^1\nnh$. The standard isomorphic identification
\begin{equation}\label{hts}
H_2(\cp^1\!\times\cp^1\nnh,\bbZ_{[2]})\,=\,\bbZ_{[2]}\times\bbZ_{[2]}\,,\qquad
\text{\rm with}\hskip6pt\bbZ_{[2]}\hskip4pt\text{\rm as\ in\ (\ref{ztw}),}
\end{equation}
in which the generators $\,(1,0),\hs(0,1)\,$ of $\,\bbZ_{[2]}\times\bbZ_{[2]}$ 
correspond to cycles of the form $\,\cp^1\times\{y\}\,$ and 
$\,\{x\}\times\cp^1\nnh$, gives
\begin{equation}\label{abc}
(a,b)\cdot(a',b')\,=\,\,ab'\,+\hs\,ba'\hskip9pt\text{\rm for}\hskip6pt
a,b,a'\nnh,b'\nh\in\bbZ_{[2]}\hs,
\end{equation}
$\cdot\,$ being the $\,\bbZ_{[2]}$-val\-ued intersection form in 
$\,H_2(\y,\bbZ_{[2]})$. In terms of (\ref{hts}),
\begin{equation}\label{abo}
\text{\rm Ker}\,[c_{\hs1}(\y)]\,=\,\{(a,b)\in\bbZ\times\bbZ:a+b=0\}\,
\subset\,H_2(\y,\bbZ)\hs.
\end{equation}
In fact, $\,c_{\hs1}(\y)\,$ treated as a homomorphism 
$\,\bbZ\times\bbZ=H_2(\y,\bbZ)\to\bbZ\,$ acts by $\,(a,b)\mapsto2(a+b)$, 
since it sends both $\,(1,0)\,$ and $\,(0,1)\,$ to $\,2$.

For the \kb\ $\,\x=K^2$ and $\,\y=\cp^1\!\times\cp^1\nnh$, one now easily sees 
that $\,[\{(1,0),(3,0)\}\times\{(0,\nh0),(1,1)\}]
\cup[\{(1,2),(3,2)\}\times\{(1,0),(0,1)\}]\,$ is the set 
$\,\mathcal{Z}(\x,\y)\,$ is 
defined in the lines following Theorem~\ref{thuno}. (This is clear from 
(\ref{ifr}) with $\,H_1(K^2,\bbZ)=\bbZ\times\bbZ_2$, cf.\ Section~\ref{ch}, and 
(\ref{abo}).) Thus, $\,\mathcal{Z}(\x,\y)\,$ is an $\,8$-el\-e\-ment subset of 
the $\,16$-el\-e\-ment set $\,\fri_4(\x)\times\dzerp(\y)$. Here 
$\,\fri_4(\x)\,=\,\{1,3\}\times\{0,2\}\,$ by (\ref{ifr}), $\,q=4\,$ ((e) in 
Section~\ref{tq}), and $\,\dzerp(\y)=H_2(\y,\bbZ_2)=\{0,1\}\times\{0,1\}\,$ 
for $\,\y\,=\,\cp^1\!\times\cp^1\nnh$, by (\ref{dmk}).
\begin{prop}\label{zkbsp}Let\/ $\,\y=\cp^1\!\times\cp^1$ and\/ $\,\x=K^2\nnh$. 
Every element\/ $\,(\ind,\dg)\,$ of\/ $\,\mathcal{Z}(\x,\y)\,$ then equals\/ 
$\,(\ind(f),\dg(f))\,$ for some \tri\/ $\,f:\x\to\y\nh$. If, in addition, 
$\,\dg\hs\ne(1,1)$, then\/ $\,(\ind,\dg)=(\ind(f),\dg(f))\,$ for some \tre\/ 
$\,f:\x\to\y\nh$.
\end{prop}
We will prove Proposition~\ref{zkbsp} at the end of this section. 
\begin{rem}\label{dadsp}A suitably-oriented an\-ti-di\-ag\-o\-nal \tr\ 
$\,2$-sphere $\,\x\,$ in 
$\,\cp^1\!\times\cp^1$ (see (vii) in Section~\ref{se}) has the degree 
$\,\hs\dg\,=(1,-1)\in\bbZ\times\bbZ$, cf.\ (\ref{hts}). In fact, writing 
$\,\hs\dg\,=(a,b)\,$ we get $\,b=\mp1$, as the composite of the mapping 
$\,x\mapsto(\overline x,x)\,$ followed by the projection onto the second 
factor is $\,\hs\text{\rm Id}:\cp^1\to\cp^1\nnh$. Hence $\,a=\pm1\,$ due to 
relations (\ref{wco}.b) and (\ref{abo}).
\end{rem}
\begin{rem}\label{ifbih}Given a \tri\ $\,f:\x\to\y\,$ of a \compact\ \rsu\ 
$\,\x\,$ in a \sc\ \acsu\ $\,\y\,$ and a \feo\ $\,h:\y\to\y\,$ preserving the 
\acst, we have $\,\ind(h\circ f)=\ind(f)$. The reason is that $\,h$, when 
naturally lifted to $\,E(\y)$, acts trivially on $\,\pi_1[E(\y)]$, as it 
commutes with the $\,\hs\text{\rm U}\hs(1)\,$ action which provides a 
generator of $\,\pi_1[E(\y)]$.
\end{rem}
\begin{lem}\label{trbdl}Suppose that\/ $\,\n,\p\,$ are almost complex 
manifolds, $\,\x\,$ is the total space and\/ $\,\,\text{\rm pr}:\x\to\q\,$ is 
the projection of a locally trivial fibre bundle over a \rmf\/ $\,\q$, and\/ 
$\,h:\q\to\n\,$ is a \tr\ \ie. Moreover, let\/ $\,\varPhi:\x\to\p\,$ be a 
$\,C^\infty$ mapping whose restriction to each fibre $\,\x_y$, $\,y\in\q$, is 
a \tr\ \ie. Then\/ $\,f\nh\,$ given by 
$\,f(x)\,=\,(h(\text{\rm pr}(x)),\varPhi(x))$, for\/ $\,x\in\x$, is a \tr\ 
\ie\ of\/ $\,\x\,$ in the product \acm\/ $\,\y\nh=\n\times\p$. 
\end{lem}
In fact, nonzero horizontal and vertical vectors in $\,\x\,$ have 
$\,\df$-images that are linearly independent over $\,\bbC\nh$.
\begin{proof}[Proof of Proposition~\ref{zkbsp}]Again, let 
$\,\y=\cp^1\!\times\cp^1$ and $\,\x=K^2\nnh$. It suffices to show that every 
$\,\hs\dg\,=(a,b)\in\{0,1\}\times\{0,1\}=\bbZ_2\times\bbZ_2=H_2(\y,\bbZ_2)\,$ 
is realized as the degree $\,\hs\dg(f)\,$ of some \tri\ $\,f:\x\to\y$, 
injective if $\,\dg\hs\ne(1,1)$, and having the Maslov index 
$\,\hs\ind(f)\in\,\fri_4(\x)\,=\,\{1,3\}\times\{0,2\}\,$ with the required 
second component $\,2|a-b|$. Namely, when $\,f\,$ is replaced by its 
composite $\,f'$ with the \feo\ $\,K^2\to K^2$ described immediately 
before Remark~\ref{srfcs}, the possible values $\,1,3\,$ of the first 
component of $\,\hs\ind(f)\,$ become interchanged, while the second component 
and the degree remain the same. (In fact, that \feo\ acts by 
$\,(a,b)\,\mapsto\,(-\hs a,b)\,$ in $\,H^1(K^2\nnh,G)=G\times G\ortwo$, cf.\ 
(\ref{hkb}), for any Abelian group $\,G$, due to the way it operates on  
$\,H_1(K^2,\bbZ)$, as described before Remark~\ref{srfcs}; on the other hand, 
$\,3\,$ is the opposite of $\,1\,$ in $\,\bbZ_4$.)

First, Corollary~\ref{conas}(c) and Remark~\ref{dadsp} allow us to choose 
\tre s $\,f:\x\to\y\,$ and $\,f':S^2\to\y\,$ with $\,\hs\dg(f)\,=(0,\nh0)\,$ 
(proving realizability of $\,\hs\dg\,=(0,\nh0)$) and 
$\,\hs\dg(f')\,=(1,-1)\in\bbZ\times\bbZ$. Lemma~\ref{iifdf} now provides a 
\tri\ $\,\x=\x\,\#\,S^2\to\y\,$ with the degree 
$\,\hs\dg\,=(1,1)\in\bbZ_2\times\bbZ_2$ and the Maslov index 
$\,\hs\ind\,=\,\ind(f)\in\,\fri_4(\x)\,$ (where $\,\hs\ind(f')=0\,$ by 
(\ref{isn}), and we identify $\,\fri_4(\x)\times\{0\}\,$ with $\,\fri_4(\x)$). 
Thus, $\,\hs\dg\,=(1,1)\,$ is realized as well.

Next, to realize the degree $\,\hs\dg\,=(0,1)$, let us identify $\,\cp^1$ with 
the unit sphere $\,S^2$ centered at zero in $\,\rtr\nnh$, and define a \tre\ 
$\,f:\x\to\y\,$ as in Lemma~\ref{trbdl}, choosing the base $\,\q\,$ to be any 
circle embedded in $\,\cp^1\nnh$, and the fibres to be the family of all great 
circles in $\,S^2=\cp^1$ that contain a fixed pair of antipodal points 
$\,u,-\hs u\in S^2\nnh$. As a point $\,y\,$ traverses the base circle $\,\q$, 
we require the corresponding fibre $\,\x_y$ to vary by being rotated about the 
axis $\,\bbR u\,$ in $\,\rtr$ in such a way that, after $\,y\,$ has gone all 
the way around $\,\q$, the fibre circle will have undergone a rotation by the 
angle $\,\pi$. That the resulting embedding $\,f:K^2\to\y=\cp^1\!\times\cp^1$ 
is \tr\ follows from Lemma~\ref{trbdl} (the total-reality assumptions being 
satisfied for dimensional reasons, cf.\ (i) in Section~\ref{se}). That  
$\,\hs\dg(f)=(0,1)\,$ is in turn clear if one considers the images of 
$\,\hs\dg(f)\,$ under the factor projections: the projection of $\,f\,$ onto 
the first $\,\cp^1$ factor is not surjective, while $\,f\,$ projected onto the 
second factor is the blow-down projection $\,K^2\to\cp^1$ 
(with $\,K^2$ treated as a two-point 
real blow-up of $\,\cp^1=S^2$).

Finally, for $\,f\,$ chosen as above, with $\,\hs\dg(f)\,=(0,1)$, let $\,f'$ 
be the composite of $\,f\,$ with the switch mapping 
$\,\y\ni(y,z)\mapsto(z,y)\in\y$. The \tre\ $\,f':\x\to\y\,$ obviously has 
$\,\hs\dg(f')=(1,0)$. As $\,\ind(f')=\ind(f)\,$ (see Remark~\ref{ifbih}), this 
completes the proof.
\end{proof} 

\section{Surfaces immersed in $\,\cp^2\,\#\,\,\mk\hs\overline{\cp^2}$}\label{su}
\setcounter{equation}{0}
For $\,\y,\y'$ as in (\ref{htw}) and 
$\,c_{\hs1}(\,\cdot\,):H_2(\,\cdot\,,\bbZ)\to\bbZ$, we have
\begin{equation}\label{con}
[c_{\hs1}(\y')](\xi,q_1,\dots,q_\mk)\,
=\,[c_{\hs1}(\y)](\xi)\,-\,q_1\hs-\,\ldots\,-\,q_\mk\hs.
\end{equation}
whenever $\,\xi\in H_2(\y,\bbZ)\,$ and $\,q_1,\dots q_\mk\in\bbZ$. In fact, 
for the tangent bundle $\,\tau\,$ and normal bundle 
$\,\nu\,$ of any of the $\,\mk\,$ exceptional divisors $\,\x\,$ in $\,\y'$ 
resulting from the blow-up, 
$\,[c_{\hs1}(\y')]([\x\hs])=c_{\hs1}(\tau)+c_{\hs1}(\nu)
=2+c_{\hs1}(\nu)\,$ and $\,c_{\hs1}(\nu)=[\x\hs]\,\cdot\,[\x\hs]=-1$. (Note 
the sign convention used in (\ref{htw}).)

We now discuss a special case, using further notational conventions. First, 
$\,\y\,=\,\cp^2\#\,\mk\hs\overline{\cp^2}$ stands, in the rest of this section, 
for the complex surface obtained by blowing up any ordered set of $\,\mk\ge1\,$ 
distinct points in $\,\cp^2\nnh$, while $\,\langle\,,\rangle\,$ and $\,\bbZ^\mk$ 
are the standard Euclidean inner product of $\,\rmk$ and the additive subgroup 
of $\,\rmk$ generated by the standard basis 
$\,\hs\text{\bf e}_1,\dots,\,\text{\bf e}_\mk$. The identification (\ref{htw}) 
now yields
\begin{equation}\label{htk}
H_2(\cp^2\,\#\,\,\mk\hs\overline{\cp^2}\nnh,\bbZ)\,=\,
\bbZ^{\mk+1}\,=\,\bbZ\hs\times\bbZ^\mk\,\subset\,\bbZ\hs\times\rmk.
\end{equation}
Specifically, the homology classes which 
correspond here to the elements $\,(1,\,\text{\bf0})\in\bbZ\times\rmk$ and 
$\,(0,\,\text{\bf e}_j)\in\bbZ\times\rmk$, $\,j=1,\dots,\mk$, are realized 
by a projective line in $\,\cp^2$ not containing any of the blown-up points 
(with its standard orientation, described at the beginning of Section~\ref{tp}) and, 
respectively, by the $\,\mk\,$ embedded copies of $\,\cp^1$ that replace the 
blown-up points, each of them with the {\it opposite\/} of its standard 
orientation. Therefore,
\begin{equation}\label{dpd}
(d,\text{\smallbf q})\cdot(d\hs'\nnh,\text{\smallbf q}\hh')\,=\,dd\hs'\,
-\,\langle\text{\smallbf q},\text{\smallbf q}\hh'\rangle
\end{equation}
for the intersection form $\,\cdot\,$ in $\,H_2(\y,\bbZ)$, while 
$\,c_{\hs1}(\y):H_2(\y,\bbZ)\to\bbZ\,$ acts by 
$\,(d,\text{\smallbf q})\mapsto3d-q_1-\ldots-q_\mk$, where 
$\,\hs\text{\smallbf q}=(q_1,\ldots,q_\mk)\,$ (see (\ref{con})). Thus,
\begin{equation}\label{krd}
\text{\rm Ker}\,[c_{\hs1}(\y)]
=\{(d\hh,q_1,\dots,q_\mk):q_1\nh+\ldots+q_\mk\nh=3d\}
\subset H_2(\y,\bbZ)\hs.
\end{equation}
\begin{example}\label{unotr}The \tri\ $\,f:S^2\to\cp^2\#\,\overline{\cp^2}$ 
described in Corollary~\ref{untre} for $\,\mk=1\,$ has $\,\hs\dg(f)=(1,3)$. In 
fact, let  $\,\hs\dg(f)=(d,\text{\smallbf q})$. Then $\,d=1\,$ due to the 
intersection equality in Corollary~\ref{untre} and (\ref{dpd}) for 
$\,(d\hs',\text{\smallbf q}\hh')=(1,\text{\bf0})$, so that (\ref{wco}.b) and 
(\ref{krd}) for $\,\mk=1\,$ give $\,(d,\text{\smallbf q})=(1,3)$.
\end{example}
\begin{example}\label{round}A \tr\ embedded $\,2$-sphere 
$\,\x\subset\cp^2\#\,2\,\overline{\cp^2}$ with the degree 
$\,[\x\hs]=(0,1,-1)$, for a suitable orientation of $\,\x$, can be 
obtained as follows. Using a complex automorphism of $\,\cp^2\nnh$, we may 
assume that the blown-up points are 
$\,x^\pm=(0,\pm\hs a)\in\bbC\times\bbR\subset\bbC^2\subset\cp^2\nnh$, with 
$\,a>0$. The $\,2$-sphere $\,\s\,$ of radius $\,a\,$ in $\,\bbC\times\bbR\hs$, 
centered at $\,(0,\nh0)$, contains $\,x^\pm$ as its unique complex points, both 
removable by blow-up (Example~\ref{rdsph}). Blowing them up transforms 
$\,\s\,$ into a totally real embedded sphere 
$\,\x\,$ in $\,\cp^2\#\,2\,\overline{\cp^2}\nnh$. To see that 
$\,[\x\hs]=(0,1,-1)$, set $\,[\x\hs]=(d,\text{\smallbf q})$. Thus, $\,d=0$ 
from (\ref{dpd}) applied to the homology class 
$\,(d\hs',\text{\smallbf q}\hh')=(1,\text{\bf0})$, represented by a projective 
line in $\,\cp^2$ not intersecting $\,\s$. Next, (\ref{dot}) (for $\,n=2$) and 
(\ref{dpd}) yield $\,\langle\text{\smallbf q},\text{\smallbf q}\hs\rangle
=-\hs(d,\text{\smallbf q})\cdot(d,\text{\smallbf q})=\chi(\x)=2$. The two 
components of $\,\hs\text{\smallbf q}\hs\,$ thus have the absolute value 
$\,1$, and opposite signs (by (\ref{wco}.b) combined with (\ref{krd}) for 
$\,d=0\,$ and $\,\mk=2$).
\end{example}
\begin{rem}\label{gnset}Condition (\ref{thx}) is satisfied by 
$\,\y=\cp^2\#\,\mk\hs\overline{\cp^2}\nnh$, $\,\mk\ge1$. In fact, 
$\,\hs\text{\rm Ker}\,[c_{\hs1}(\y)]\,$ is generated by 
$\,(1,3\hs\text{\bf e}_1)\,$ and $\,(0,\,\text{\bf e}_j\nh-\text{\bf e}_1)$, 
$\,j=2,\dots,\mk$, with $\,\hs\text{\bf e}_j$ as above. (See (\ref{krd}).) These 
generators are all realized by \tri s, in view of Examples~\ref{unotr}, 
\ref{round} and Lemma~\ref{blwup}(a).
\end{rem}
\begin{proof}[Proof of Theorem~\ref{maire}]We will verify below that the 
\csu s (\ref{lis}) all satisfy the assumptions, and hence conclusions, of 
Theorem~\ref{clase}. Thus, \tri s of \compact\ real surfaces $\,\x\,$ in the 
\csu s $\,\y\,$ forming the list (\ref{lis}) realize, as their index-degree 
pairs, all elements $\,(\ind\hs,\dg)\,$ of the set 
$\,\mathcal{Z}(\x,\y)\subset\,\fri_q(\x)\times\depspm(\y)\,$ described in the 
lines following Theorem~\ref{thuno}. This will clearly imply 
Theorem~\ref{maire}, since in Section~\ref{sj} we already proved that, conversely, 
all such index-degree pairs lie in $\,\mathcal{Z}(\x,\y)$.

First, (\ref{ife}) holds if $\,\y\,$ is 
$\,\bbC^2\nnh,\,\cp^1\!\times\cp^1\nnh,\,\cp^2\nnh$, or 
$\,\cp^2\#\,\mk\hs\overline{\cp^2}\nnh$, $\,\mk\ge1$, as $\,q/2\,$ then equals, 
respectively, $\,\infty,2,3\,$ or $\,1\,$ (see (e) -- (f) in Section~\ref{tq}), and 
in the last two cases a \tri\ $\,\rp^2\nnh\to\y\,$ exists according to (vi) in 
Section~\ref{se} combined with Lemma~\ref{blwup}(a).

Next, (\ref{thx}) is satisfied by $\,\bbC^2$ and $\,\cp^2$ in view of 
Corollary~\ref{sfinm} (since 
$\,\hs\text{\rm Ker}\,[c_{\hs1}(\cdot)]=\{0\}$); by 
$\,\cp^1\!\times\cp^1\nnh$, as a consequence of (\ref{abo}) and 
Remark~\ref{dadsp}; and by $\,\cp^2\#\,\mk\hs\overline{\cp^2}\nnh$, $\,\mk\ge1$, 
according to Remark~\ref{gnset}.

Finally, the assumption about $\,K^2$ in Theorem~\ref{clase} holds for 
$\,\cp^1\!\times\cp^1$ by Proposition~\ref{zkbsp}, and for the remaining 
\csu s (\ref{lis}) it is immediate from Lemma~\ref{cedif}: (\ref{evr}) is 
satisfied, with $\,\x=K^2\nnh$, by Corollary~\ref{conas}(b), since 
$\,q/2\in\{\infty,3,1\}\,$ (as remarked above), while (\ref{thx}) was verified 
in the last paragraph. This completes the proof.
\end{proof}

\section{Some special cases of Theorem~\ref{maiem}}\label{sc}
\setcounter{equation}{0}
Let $\,\x\,$ and $\,\y\,$ be, respectively, a fixed \compact\ \rsu\ and one of 
the \csu s (\ref{lis}). The pairs $\,(\ind,\dg)\,$ that can be simultaneously 
realized as the Maslov index and degree of a \tri\ 
$\,\x\to\y\,$ then are nothing else than all elements of the set 
$\,\mathcal{Z}(\x,\y)$. 
(See Theorem~\ref{maire} and the discussion following it in Section~\ref{sr}.) 
This is, however, not the case if one replaces the word `immersion' by 
`embedding' since, in view of (\ref{dot}) and Corollary~\ref{sigsq} (for 
$\,n=2$), the degree $\,\hs\dg\,=\,\dg(f)\,$ of any \tre\ $\,f:\x\to\y\,$ 
satisfies the additional constraint (\ref{ddh}).

Theorem~\ref{maiem} states that there are no further constraints, as long as 
$\,\mk\le7\,$ in (\ref{lis}). In other words, the conditions 
$\,(\ind,\dg)\in\mathcal{Z}(\x,\y)\,$ and (\ref{ddh}) hold if and 
only if the pair $\,(\ind,\dg)\,$ is realized by a \tre\ $\,\x\to\y$, where 
$\,\y\,$ is $\,\bbC^2$, $\,\cp^2$, $\,\cp^1\!\times\cp^1$, or 
$\,\cp^2\#\,\mk\hs\overline{\cp^2}\nnh$, $\,1\le\mk\le7$. 

In this section we establish Theorem~\ref{maiem} in the case where $\,\x\,$ is 
orientable and $\,\y\,$ is $\,\bbC^2\nnh$, $\,\cp^2$ or 
$\,\cp^1\!\times\cp^1\nnh$. (For the other cases, see Sections~\ref{rc} and 
\ref{tn}.) Note that 
$\,\mathcal{Z}(\x,\y)=\hs\fri_q(\x)\times\,\text{\rm Ker}\,[c_{\hs1}(\y)]\,$ 
when $\,\x\,$ is 
orientable, according to (\ref{dfd}), (\ref{dmk}) and the lines following 
Theorem~\ref{thuno}.

First, if $\,\x\,$ is orientable and $\,\y=\bbC^2$ or $\,\y\,=\,\cp^2\nh$, 
relations (\ref{dcp}), (\ref{dco}) and (\ref{dfd}) give 
$\,\dg\in\,\text{\rm Ker}\,[c_{\hs1}(\y)]=\{0\}$, and so, by (\ref{ddh}), 
$\,\x\,$ must be the torus $\,T^2\nnh$. Our assertion is now obvious from 
Corollary~\ref{conas}(a).

To verify Theorem~\ref{maiem} for orientable $\,\x\,$ and 
$\,\y\,=\,\cp^1\!\times\cp^1\nnh$, we use the notations of (\ref{hts}) -- 
(\ref{abc}). In view of (\ref{abc}) -- (\ref{abo}), conditions (\ref{ddh}) and 
$\,\hs\dg\in\,\text{\rm Ker}\,[c_{\hs1}(\y)]$, for any given 
\compact\ oriented surface $\,\x$, now read $\,a+b=0\,$ and 
$\,2ab=-\hs\chi(\x)$, where $\,\hs\dg\,=(a,b)$. This gives 
$\,\chi(\x)=-\hs2ab=2a^2\ge0$, that is, $\,\x\,$ can only be the torus 
$\,T^2$ with $\,a=b=0\,$ or the sphere $\,S^2$ with 
$\,(a,b)=(\pm\,1,\mp\,1)$. A \tre\ $\,\x\to\y\,$ realizing the pair 
$\,\hs(\ind,\dg)\hs\,$ for any $\,\hs\ind\,\in\,\fri_4(\x)\,$ (cf.\ (e) in 
Section~\ref{tq}) and such $\,\hs\dg\,=(a,b)\,$ thus exists in view of 
Corollary~\ref{conas}(a) (for $\,\x=T^2$) or Remark~\ref{dadsp} (for 
$\,\x=S^2$). Note that, if $\,\x=S^2\nnh$, {\it every\/} value of 
$\,\hs\ind\,\in\,\fri_4(\x)\,$ is realized, since (\ref{isn}) allows just one 
value, $\,\hs\ind\hs=0$.

\section{A Diophantine equation}\label{de}
\setcounter{equation}{0}
Let $\,\x\,$ be an oriented \compact\ \rsu, and let $\,\y\,$ be the complex 
surface obtained by blowing up any set of $\,\,\mk\,$ distinct points in 
$\,\cp^2\nnh$, for $\,\mk=1,\dots,8$. Thus, 
$\,\mathcal{Z}(\x,\y)=\hs\fri_q(\x)\times\,\text{\rm Ker}\,[c_{\hs1}(\y)]\,$ 
(see Section~\ref{sc}). As a first step toward proving Theorem~\ref{maiem} for such 
$\,\x\,$ and $\,\y$, we will provide in this section an explicit description 
of the subset of $\,\mathcal{Z}(\x,\y)\,$ defined by the additional 
requirement (\ref{ddh}). However, as $\,q=2\,$ ((f) in Section~\ref{tq}) and, by 
(\ref{iqs}), $\,\hs\fri_2(\x)=\{w_1(\x)\}=\{0\}$, the first component 
$\,\hs\ind\hs\,$ of any element $\,\hs(\ind,\dg)\hs$ of that subset is 
uniquely determined, so that we need only to find the corresponding set of the 
second components $\,\hs\dg\hh$. We use the identification (\ref{htk}) to 
treat $\,\hs\dg\hs\,$ as a pair $\,(d,\text{\smallbf q})\,$ with 
$\,d\in\bbZ\,$ and $\,\hs\text{\smallbf q}\,\in\bbZ^\mk\subset\rmk$ or, 
equivalently, as a $\,(\mk+1)$-tuple $\,(d\hh,q_1,\dots,q_\mk)\,$ of integers. 
We also set $\,\chi=\chi(\x)\,$ and $\,\,\text{\bf1}\,=(1,\dots,1)\in\rmk\nh$. 
The conditions imposed on $\,\hs\dg\hs\,$ are 
$\,\hs\dg\in\,\text{\rm Ker}\,[c_{\hs1}(\y)]\,$ and 
$\,\hs\dg\,\cdot\,\dg\,=\hs-\hs\chi$, which, by (\ref{dpd}) -- (\ref{krd}), 
amounts to
\begin{equation}\label{pot}
\langle\text{\smallbf q}\hskip.4pt,\text{\bf1}\rangle\,=\,3d\,,\qquad
|\text{\smallbf q}\hs|^2\,=\,d^{\hs2}\hs+\,\chi\,,\qquad d\in\bbZ\,,\quad
\text{\smallbf q}\in\bbZ^\mk\nh.
\end{equation}
Here $\,\langle\,,\rangle\,$ is the inner product of $\,\rmk\nnh$, and 
$\,|\,\,|\,$ is the corresponding norm.

We now proceed to solve equations (\ref{pot}) with $\,\chi\le2\,$ and 
$\,1\le\mk\le8$.
\begin{example}\label{xmpls}Treating (\ref{pot}) as a system of equations 
imposed on $\,d\,$ and $\,\hs\text{\smallbf q}\,=(q_1,\dots,q_\mk)$, in which 
$\,\mk,\chi\in\bbZ\,$ are fixed parameters with $\,\mk\ge1$, we can rewrite it 
as $\,q_1+\ldots+q_\mk=3d\,$ and $\,q_1^2+\ldots+q_\mk^2=d^{\hs2}+\chi$, the 
unknowns now being $\,d,q_1,\dots,q_\mk\in\bbZ$. Thus, since 
$\,q^2\equiv q\hskip6pt\text{\rm mod}\hskip4pt2\,$ for any $\,q\in\bbZ$, a 
solution to (\ref{pot}) exists only if $\,\chi\,$ is even. Each of the 
following three families of solutions to (\ref{pot}) represents infinitely 
many values of $\,\mk$ (with $\,d\,$ always denoting an integer):
\begin{enumerate}
  \def\theenumi{{\rm\roman{enumi}}}
\item $\,(d,\text{\smallbf q})=(0,\text{\bf0})=(0,\hs0,\dots,0)$, with any 
$\,\mk\ge1\,$ and $\,\chi=0$,
\item $\,(d,\text{\smallbf q})=(d,\text{\bf1})=(d\hh,\hs1,\dots,1)$, 
for $d\ge1$, with $\mk=3d$ and $\chi=(3-d)d$,
\item $\,(d,\text{\smallbf q})=(d\hh,\hs d\nh-\nnh1,1,\dots,1)$, for any 
$\,d\ge0$, with $\,\mk=2d+2\,$ and $\,\chi=2$.
\end{enumerate}
As we will show in Lemma~\ref{solns} below, a solution to (\ref{pot}) with 
$\,\chi<0$ exists only if $\,\mk\ge10$. In this regard, $\,\mk=10\,$ is a 
threshold value: solutions 
$\,(d,\text{\smallbf q})=(d\hh,q_1,\dots,q_\mk)\,$ with $\,\mk=10\,$ not only 
exist for any prescribed even negative integer $\,\chi$, but can also be 
chosen so that $\,d>0\,$ and $\,q_j>0$ for all $\,j$, as illustrated by 
$\,(d,\text{\smallbf q})
=(3\hh c+1+2\hs\ve\hs;\,c\hs,c\hs,c\hs,c\hh+\ve,\dots,c\hh+\ve,3)$, for 
integers $\,c\ge2\,$ and $\,\ve\in\{0,1\}$, with $\,\chi=8-2\hh\ve-6\hh c\hs$, 
or $\,(d,\text{\smallbf q})
=(3\hh c+4\hs;\,c+2,c\hs,c+1,\dots,c+1,3)$, for integers $\,c\ge1$, with 
$\,\chi=4-6\hh c\hs$.
\end{example}
\begin{rem}\label{staso}Equations (\ref{pot}) remain satisfied after any 
permutation of the components $\,q_1,\dots,q_\mk$ of 
$\,\hs\text{\smallbf q}\hs$, as well as after the signs of $\,d\,$ and all 
$\,q_1,\dots,q_\mk$ have been changed; also, a new solution with $\,\mk\,$ 
replaced by $\,\mk'>\mk$ (or, $\,\mk'<\mk$) arises if one inserts additional 
$\,\mk'-\mk\,$ zeros (or, respectively, deletes existing $\,\mk-\mk'$ zeros) among the 
$\,q_1,\dots,q_\mk$. Successive applications of these operations, repeated in 
any order, any number of times, lead to what we call {\it trivial 
modifications\/} of the given solution $\,(d,\text{\smallbf q})\,$ to 
(\ref{pot}).
\end{rem}
\begin{lem}\label{solns}No solutions to\/ {\rm(\ref{pot})} exist when\/ 
$\,\chi<0\,$ and\/ $\,1\le\mk\le9$. If\/ $\,\chi=0$, the only solution\/  
$\,(d,\text{\smallbf q})\,$ with $\,1\le\mk\le8\,$ is\/ $\,(0,\text{\bf0})$, 
while the only solutions with $\,\chi=0\,$ and\/ $\,\mk=9\,$ are\/ 
$\,(3s\hh,\hs s,s,s,s,s,s,s,s,s)\,$ for\/ $\,s\in\bbZ$.
\end{lem}
\begin{proof}The Schwarz inequality $\,\langle\text{\smallbf q}\hskip.4pt,
\text{\bf1}\rangle^2\le|\text{\smallbf q}\hs|^2|\text{\bf1}|^2$ becomes 
$\,(9-\mk)\hs d^{\hs2}\le\mk\chi\,$ for $\,(d,\text{\smallbf q})\,$ with 
(\ref{pot}). If $\,1\le\mk\le9$, this gives $\,\chi\ge0$. If, in addition, 
$\,\chi=0$, our Schwarz inequality yields $\,(9-\mk)d=0$, and so 
$\,\hs\text{\smallbf q}\hs\,$ is a multiple of $\,\hs\text{\bf1}$ (the 
equality case in the Schwarz inequality); thus, either $\,\chi=9-\mk=0$, or 
$\,\mk<9\,$ and $\,\chi=d=0$, which completes the proof.
\end{proof} 
For $\,(d,\text{\smallbf q}),\mk,\chi\,$ with (\ref{pot}) and the greatest 
integer $\,s\,$ with $\,3s\le d+1$,
\begin{equation}\label{dsr}
d\,=\,3s+r\,,\qquad\text{\rm while}\hskip7pts\in\bbZ\hskip7pt\text{\rm and}
\hskip7ptr\in\{-1,0,1\}\,.
\end{equation}
Setting $\,\hs\text{\smallbf s}\,=s\hs\text{\bf1}\,=(s,\dots,s)\in\rmk$ we now 
have $\,|\text{\smallbf q}\hh|^2=d^{\hs2}+\chi=9s^2+6\hskip.6ptrs+r^2+\chi$, 
$\,|\text{\smallbf s}\hs|^2=\mk s^2$ and 
$\,\langle\text{\smallbf q}\hskip.4pt,\text{\smallbf s}\rangle=3sd=3s(3s+r)
=9s^2+3rs$. Hence
\begin{equation}\label{psk}
|\text{\smallbf q}-\text{\smallbf s}\hs|^2=(\mk-9)s^2+r^2+\chi\quad\text{\rm 
and}\quad r^2\in\{0,1\}\hs.
\end{equation}
\begin{lem}\label{oimin}The only solutions 
$\,(d,\text{\smallbf q})=(d\hh,q_1,\dots,q_\mk)\,$ to {\rm(\ref{pot})} with 
$\,\chi=2\,$ and\/ $\,1\le\mk\le8\,$ are
\begin{equation}\label{oii}
(0\hskip.9pt;\,1,-1)\,,\hskip6pt(1;\,1,1,1)\,,
\hskip6pt(2\hskip.4pt;\,1,1,1,1,1,1)\,,
\hskip6pt(3\hskip.4pt;\,2,1,1,1,1,1,1,1)\,,
\end{equation}
and those obtained from them by trivial modifications, cf.\ {\rm 
Remark~\ref{staso}}.

Note that, up to trivial modifications, {\rm(\ref{oii})} are precisely the 
solutions\/ {\rm(iii)} in Example~\ref{xmpls} for\/ $\,d=0,1,2,3$.
\end{lem}
\begin{proof}Let $\,(d,\text{\smallbf q})=(d\hh,q_1,\dots,q_\mk)\,$ satisfy 
(\ref{pot}) with $\,\chi=2\,$ and $\,1\le\mk\le8$. After a trivial 
modification we get
\begin{equation}\label{ppk}
d\,\ge\,0\hskip15pt\text{\rm and}\hskip13ptq_1\ldots q_\mk\,\ne\,0\,.
\end{equation} 
By (\ref{psk}), $\,s^2\le(9-\mk)s^2\le\chi+1=3$, and so $\,s\in\{-1,0,1\}$. 
Since $\,d\ge0$, we have $\,s\in\{0,1\}$, cf.\ (\ref{dsr}). For 
$\,\ell=|\text{\smallbf q}-\text{\smallbf s}\hs|^2\nh$, (\ref{psk}) with 
$\,\chi=2$ gives $\,\ell\in\{2,3\}\,$ (when $\,s=0$), or 
$\,\ell\in\{\mk-7,\mk-6\}\,$ (when $\,s=1$). As $\,\ell\ge0$, it follows that 
$\,\mk\in\{6,7,8\}\,$ if $\,s=1$, while $\,\mk=\ell\in\{2,3\}\,$ and 
$\,|q_1|=\ldots=|q_\mk|=1\,$ if $\,s=0$, as 
$\,\ell=|\text{\smallbf q}\hs|^2\le3\,$ and, by (\ref{ppk}), 
$\,q_1\ldots q_\mk\ne0$. The triple $\,(s,\mk,\ell)\,$ thus must assume one of 
the seven values $\,(0,2,2),(0,3,3),(1,6,0),(1,7,0),(1,7,1),(1,8,1)\,$ and 
$\,(1,8,2)$. When $\,s=1$, the relations 
$\,\ell=\sum_{j=1}^\mk(q_j\nh-1)^2\le2\,$ and (\ref{ppk}) imply that
\begin{equation}\label{klo}
\text{\rm$\ell\,$\ of\ the\ $\hs\mk\hs$\ integers\ $\hs q_1,\dots,q_\mk$\ 
equal\ $\hs2\hs$\ and\ $\hs\mk\nh-\ell\hs$\ of\ them\ equal\ $\hs1$.}
\end{equation}
Thus, if $\,s=1$, we have $\,q_1+\ldots+q_\mk=\mk+\ell$, while, by 
(\ref{pot}), $\,q_1+\ldots+q_\mk=3d$, and so $\,\mk+\ell=3d$, which eliminates 
three of the seven triples $\,(s,\mk,\ell)$, leaving only those with $\,s=0\,$ 
or $\,\mk+\ell\,$ divisible by three: $\,(0,2,2),(0,3,3),(1,6,0)$, and 
$\,(1,8,1)$. These triples lead to the four possibilities listed in 
(\ref{oii}). Namely, $\,(0,2,2)\,$ and $\,(0,3,3)\,$ have $\,\mk\in\{2,3\}\,$ 
and $\,|q_1|=\ldots=|q_\mk|=1$, so that the sequence of two or three integers 
$\,q_j$ with values $\,\pm1\,$ and sum $\,3d\ge0\,$ (cf.\ (\ref{pot}) and 
(\ref{ppk})) must be $\,(\pm1,\mp1)\,$ with $\,d=0$, or $\,(1,1,1)\,$ with 
$\,d=1$, as required. By (\ref{klo}), the triples $\,(1,6,0)\,$ and 
$\,(1,8,1)\,$ correspond in turn to the last two solutions in (\ref{oii}).
\end{proof} 

\section{Deformations of pseu\-do\-hol\-o\-mor\-phic immersions}\label{dp}
\setcounter{equation}{0}
We now present some facts needed in Section~\ref{po} to prove 
Theorem~\ref{thqnd}.

Given \cvs s $\,W,\,V\hs$ and $\,A\in\,\text{\rm Hom}_\bbR(W,V)\,$ let 
$\,A^\pm$ stand for the unique operators $\,W\to V\,$ such that 
$\,A^+$ is $\,\bbC$-lin\-e\-ar, $\,A^-$ is antilinear, and $\,A=A^++\,A^-\nh$. 
Thus, $\,A^\pm=\hs(A\,\mp\,i\hh\circ A\circ i)/2\,$ are the components of 
$\,A$ relative to the decomposition of $\,\hs\text{\rm Hom}_\bbR(W,V)\,$ 
into the $\,\pm1$-eigenspaces of the involution 
$\,A\mapsto -\hs i\circ\hskip.1ptA\circ i$, where $\,\hs i\hs\,$ stands for 
multiplication by $\,\hs i$. Let $\,\nabla\,$ now be a fixed connection in a 
complex vector bundle $\,\eta\,$ over an \acm\ $\,\x$. The {\it 
Cau\-chy-Rie\-mann operator} $\,\db\,$ of $\,\nabla\,$ is the 
linear differential operator that takes any $\,C^1$ section $\,\psi\,$ of 
$\,\eta\,$ to
\begin{equation}\label{cro}
\db\hs\psi\,=\,[\nabla\psi]^-,
\end{equation}
for $\,[\,\,\,]^-$ as in the preceding lines, so that 
$\,\db\hs\psi\,$ is a section of 
$\,\hs\text{\rm Hom}_\bbC(\overline{T\x},\eta)$ (where $\,\overline{T\x}\,$ 
is the complex conjugate bundle of $\,T\x$), and its value 
$\,\db\hs\psi_x$ at any $\,x\in\x\,$ is 
$\,([\nabla\psi]_x)^-\nnh$. Note that $\,\nabla\psi\,$ itself is a section of 
$\,\hs\text{\rm Hom}_\bbR(T\x,\eta)$ sending $\,v\in T_x\x$, for any 
$\,x\in\x$, to $\,\nabla_{\!v}\psi\in\eta_x$. Thus, 
$\,\db\hs\psi_x:T_x\x\to\eta_x$ and 
$\,2\hs\db\hs\psi_xv
=\nabla_{\!v}\psi\,+\,i\nabla_{\!iv}\psi\,$ for $\,v\in T_x\x$.


For a real $\,\ek$\diml\ manifold $\,\x\,$ and an almost complex manifold 
$\,\y$, let $\,[\tm]^{\wedge\ek}$ denote the $\,\ek$th complex exterior 
power of $\,\tm$, and let $\,\,\detr\hs T\x\,$ be as in (\ref{ftm}.ii). Any 
$\,C^\infty$ mapping $\,f:\x\to\y\,$ then gives rise to the vector bundle 
morphism $\,\det\df:\detr\hs T\x\to f^*([\tm]^{\wedge\ek})\,$ uniquely 
characterized by $\,(\det\df)_x(v_1\wedge\ldots\wedge\,v_\ek)=
[\df_xv_1]\wedge\ldots\wedge\,[\df_xv_\ek]\,$ for $\,x\in\x$ and 
$\,v_1,\dots,v_\ek\in T_x\x$. Obviously, $\,f\,$ is a \tri\ if 
and only if $\,\det\df\,$ is nonzero everywhere as a section of 
$\,\text{\rm Hom}_\bbR(\detr\hs T\x,f^*([\tm]^{\wedge\ek}))$. 

Given a \psh\ immersion $\,f\,$ of an oriented \rsu\ $\,\x$ in an \acsu\ 
$\,\y\,$ (cf.\ the lines following Theorem~\ref{thqnd}) and a $\,C^\infty$ 
curve (homotopy) $\,I\ni t\mapsto f^t$ of mappings $\,\x\to\y$, where 
$\,I\subset\bbR$ is an interval, such that $\,0\in I\hs$ and $\,f^{\hs0}=f$, 
let $\,\nu\,$ be the (complex) normal bundle of $\,f\,$ (see (\ref{tnu})). 
Then
\begin{equation}\label{det}
\text{\rm i)}\quad
\det\df^t\,\vrule\underset{{\hskip2ptt\hs=\hs0}}{\phantom{_jjj}}\,=\,\,0\hs,
\qquad\quad\text{\rm ii)}\quad
\frac{d}{dt}\det\df^t\,\vrule\underset{{\hskip2ptt\hs=\hs0}}{\phantom{_jjj}}\,
=\,\,\,\df\wedge\db\hs\psi\hs.
\end{equation}
Here (i) reflects the fact that $\,f^{\hs0}\nnh=f\,$ is \psh; (ii), however, 
requires further explanation. First, $\,\psi\,$ in (ii) is the section of 
$\,\nu\,$ obtained as the image of $\,(\df^t/dt)_{\hh t\hs=\hs0}$ under the 
projection morphism $\,f^*\tm\to\nu$, and $\,\db\,$ is the 
Cau\-chy-Rie\-mann operator for a suitable connection in $\,\nu$, while both 
sides of (ii) are sections of 
$\,\text{\rm Hom}_\bbR(\detr\hs T\x,f^*([\tm]^{\wedge2}))$, and 
$\,\df\wedge\db\hs\psi$ is given by 
$\,(\df\wedge\db\hs\psi)_x(u,v)=
[\hs\df_xu]\wedge[\hs\db\hs\psi_xv]
-[\hs\df_xv]\wedge[\hs\db\hs\psi_xu]\,$ for $\,x\in\x$ and 
$\,u,v\in T_x\x$. Thus, with $\,iu\in T_x\x\,$ defined using the 
$\,f$-pull\-back of the \acst\ of $\,\y\,$ to $\,\x$,
\begin{equation}\label{wiu}
(\df\wedge\db\hs\psi)_x(u,iu)\,
=\,-\hs2\hh i\hs[\hs\df_xu]\wedge[\hs\db\hs\psi_xu]\hs,
\end{equation}
due to com\-plex-lin\-e\-ar\-i\-ty of $\,\df_x$ and an\-ti\-lin\-e\-ar\-i\-ty 
of $\,\db\hs\psi_x$. Secondly, the spaces $\,\lambda_x
=\,\text{\rm Hom}_\bbR([T_x\x]^{\wedge2}\nh,[T_{f(t,x)}\y]^{\wedge2})$, with 
$\,f(t,x)=f^t(x)$, are the fibres of a complex line bundle $\,\lambda\,$ over 
$\,\x\times I\nh$, and $\,(x,t)\mapsto(\det\df^t)_x$ is, in view of (i), a 
section of $\,\lambda$, vanishing along the submanifold $\,\x\times\{0\}$. In 
general, when a $\,C^1$ section $\,h\,$ of a vector bundle 
$\,\eta\,$ over a manifold $\,N$ vanishes at a point $\,y\in N\nh$, one can 
define a linear operator $\,Dh_y:T_yN\to\eta_y$ to be the composite in which 
the differential $\,dh_y:T_yN\to T_{(y,0)}\eta\,$ (of $\,h\,$ treated as a 
mapping $\,N\to\eta\,$ into the total space) is followed by the projection 
$\,T_{(y,0)}\eta\to\eta_y$ coming from the identifications 
$\,T_{(y,0)}\eta=T_yN\oplus T_0\eta_y$ and $\,T_0\eta_y=\eta_y$. (Hence, for 
$\,\eta=T^*\hskip-2ptN\hs$ and $\,h=d\phi\,$ with any $\,C^2$ function 
$\,\phi:N\to\bbR\hs$, the Hessian of $\,\phi\,$ at the critical point $\,y\,$ 
is $\,Dh_y$.) We may now set $\,(dh(y(t))/dt)_{\hh t\hs=\hs0}=dh_yv\,$ 
whenever $\,t\mapsto y(t)\,$ is a $\,C^1$ curve in $\,N$ with $\,y(0)=0$, 
having the velocity $\,v\in T_yN\,$ at $\,t=0$. As a consequence of this for 
$\,\eta=\lambda\,$ and $\,N=\x\times I\nh$, along with any point $\,y=(x,0)\,$ 
in the submanifold $\,\x\times\{0\}$, and the curve $\,y(t)=(x,t)$, the 
left-hand side of (i) is a well-de\-fined section of 
$\,\text{\rm Hom}_\bbR(\detr\hs T\x,f^*([\tm]^{\wedge2}))$.

Equality (\ref{det}.ii) is now easily verified in local coordinates. For 
instance, one may fix a \feo\ between a \nb\ of any given 
point of $\,f(\x)\,$ and a \nb\ $\,\,U\,$ of $\,(0,\nh0)\,$ in 
$\,\bbC^2\nnh$, which makes $\,f\,$ appear as the inclusion 
$\,\,U\cap(\bbC\times\{0\})\to\bbC^2$ and, at all points of 
$\,\,U\cap(\bbC\times\{0\})$, identifies the \acst\ of $\,\y\,$ with the 
standard complex structure of $\,\bbC^2\nnh$. The connection in $\,\nu\,$ then 
is defined only locally (as well as co\-or\-di\-nate-de\-pen\-dent and 
non-u\-nique), but the corresponding Cau\-chy-Rie\-mann operator 
$\,\db\,$ is uniquely characterized by (\ref{wiu}); a globally 
defined connection leading to $\,\db\,$ may now be obtained 
via a finite partition of unity.

\section{Proof of Theorem~\ref{thqnd}}\label{po}
\setcounter{equation}{0}
\begin{lem}\label{measu}Let\/ $\,dx\,$ be a fixed positive smooth measure 
density on a compact manifold\/ $\,\x\,$ with $\,\dim\x\ge1$.
\begin{enumerate} 
  \def\theenumi{{\rm\roman{enumi}}}
\item For any fi\-nite\diml\ \vs\/ $\,\mathcal{V}\,$ of 
real-val\-ued continuous functions on\/ $\,\x$, there exist\/ $\,C^\infty$ 
functions\/ $\,f,h:\x\to\bbR\hs$, both $\,L^2$-or\-thog\-o\-nal to\/ 
$\,\mathcal{V}$, such that $\,|f|+|h|>0\,$ everywhere in\/ $\,\x$.
\item For any finite\diml\ \vs\/ $\,\mathcal{W}\,$ of com\-plex-val\-ued 
continuous functions on\/ $\,\x$, some $\,C^\infty$ function\/ 
$\,\varphi:\x\to\bbC\,$ is\/ $\,L^2$-or\-thog\-o\-nal to\/ $\,\mathcal{W}\hs$ 
and nonzero everywhere in\/ $\,\x$. 
\end{enumerate} 
\end{lem}
\begin{proof}To prove (i), set $\,\ek=\dim\mathcal{V}$, and let 
$\,\delta:\x\to\mathcal{V}^{\,*}$ be the $\,C^\infty$ mapping sending each 
$\,x\,$ to the evaluation functional (Dirac delta) $\,\delta[x]$, acting by 
$\,f\mapsto f(x)$. Its image $\,\{\delta[x]:x\in\x\}\,$ 
spans $\,\mathcal{V}^{\,*}\nnh$, as otherwise it would lie in a proper 
subspace, that is, some $\,f\in\mathcal{V}\smallsetminus\{0\}\,$ would vanish 
at all $\,x\in\x$. Thus, we may choose $\,2\ek\,$ distinct points 
$\,x_1,\dots,x_\ek,y_1,\dots,y_\ek\in\x\,$ such that both 
$\,\delta[x_1],\dots,\delta[x_\ek]\,$ and 
$\,\delta[y_1],\dots,\delta[y_\ek]\,$ are bases of $\,\mathcal{V}^{\,*}$, by 
first picking the $\,x_a$, and then selecting each $\,y_a$ near the 
corresponding $\,x_a$. Let us also fix pairwise disjoint open sets 
$\,\,U_1,\dots,U_\ek,U'_1,\dots,U'_\ek$ in $\,\x\,$ with $\,x_a\in U_a$ and 
$\,y_a\in U'_a$ for $\,a=1,\dots,\ek$. 

There must exist a $\,C^\infty$ function $\,f:\x\to\bbR\,$ which is 
$\,L^2$-or\-thog\-o\-nal to $\,\mathcal{V}\,$ and such that $\,f=1\,$ on 
$\,\x\smallsetminus U$, where $\,\,U=U_1\cup\ldots\cup U_\ek$. In fact, let 
$\,\phi_1,\dots,\phi_\ek$ be the basis of $\,\mathcal{V}\,$ dual to the basis 
$\,\delta[x_1],\dots,\delta[x_\ek]\,$ of $\,\mathcal{V}^{\,*}\nnh$. Thus, 
$\,\phi_a(x_b)=\delta_{ab}$ for $\,a,b=1,\dots,\ek$. The functions 
$\,\phi_1,\dots,\phi_\ek$ are linearly independent, when treated as linear 
functionals acting, via the $\,L^2$ inner product, on the space 
$\,\mathcal{F}\,$ of all $\,C^\infty$ functions $\,\x\to\bbR\,$ with compact 
supports contained in $\,\,U$. (Otherwise, some nontrivial combination of the 
$\,\phi_a$, being $\,L^2$-or\-thog\-o\-nal to $\,\mathcal{F}\nh$, would vanish 
everywhere in $\,\,U\nh$, which is impossible as $\,\phi_a(x_b)=\delta_{ab}$.) 
Hence, given $\,\lambda_1,\dots,\lambda_\ek\in\bbR$, there exists 
$\,\xi\in\mathcal{F}\,$ with $\,\int_\x\phi_a\xi\,dx=\lambda_a$ for 
$\,a=1,\dots,\ek$. Choosing such $\,\xi$ for 
$\,\lambda_a=-\int_\x\phi_a\,\,dx$, $\,a=1,\dots,\ek$, we can now define 
$\,f\,$ by $\,f=\xi+1$.

The same argument may be applied to the $\,y_a$ and $\,\,U'_a$ rather than 
$\,x_a$ and $\,\,U_a$. Thus, there exists a $\,C^\infty$ function 
$\,h:\x\to\bbR\,$ which is $\,L^2$-or\-thog\-o\-nal to $\,\mathcal{V}\,$ and 
such that $\,h=1\,$ on $\,\x\smallsetminus U'\nh$, with 
$\,\,U'=U'_1\cup\ldots\cup U'_\ek$. As 
$\,(\x\smallsetminus U)\cup(\x\smallsetminus U')=\x$, this yields (i). Now 
(ii) follows if we set $\,\varphi=f+ih\,$ with $\,f,h\,$ as in (i) for 
$\,\mathcal{V}=\{\text{\rm Re}\,\chi:\chi\in\mathcal{W}\}$.
\end{proof}
\begin{proof}[Proof of Theorem~\ref{thqnd}]We have 
$\,f^*[\detc\hs\tm]=\tau\otimes\nu=\,\h(\overline{\tau},\nu)\,$ as 
$\,f^*\tm=\tau\oplus\nu$, and so the equivalence of (a) and (b) is obvious. 
Furthermore, (c) implies (a) in view of (\ref{tnu}) with $\,n=2$. Now assume 
(b) and let $\,\db\,$ be the Cau\-chy-Rie\-mann operator 
appearing in (\ref{det}.ii). Since $\,\db\,$ is elliptic, for 
any prescribed $\,C^\infty$ section $\,\phi\,$ of 
$\,\h(\overline{\tau},\nu)\,$ solvability of the equation 
$\,\hs\db\hs\psi\,=\,\phi\,$ with an unknown $\,C^\infty$ 
section $\,\psi\,$ of $\,\nu\,$ is equivalent to 
$\,L^2$-or\-thog\-o\-nal\-i\-ty of $\,\phi\,$ to the kernel of the formal 
adjoint of $\,\hs\db$. (To form the adjoint, one fixes 
Hermitian fibre metrics in $\,\tau\,$ and $\,\nu$, along with a positive 
smooth measure density of on $\,\x$.) Lemma~\ref{measu}(ii) now implies that 
$\,\phi\,$ for which the equation is solvable may be chosen so as to be 
nowhere zero. (In fact, since we assume (b), we may fix a global $\,C^\infty$ 
section of $\,\h(\overline{\tau},\nu)$, unit relative to the Hermitian fibre 
metric naturally determined by those in $\,\tau\,$ and $\,\nu$, and use it to 
treat sections of $\,\h(\overline{\tau},\nu)$, including $\,\phi$, as 
functions $\,\x\to\bbC\hs$.) Choosing a deformation $\,t\mapsto f^t$ of 
$\,f\,$ in the direction of the corresponding solution $\,\psi$, and using 
(\ref{det}), we now obtain (c); note that the section appearing in 
(\ref{det}.ii) is nonzero everywhere in $\,\x\,$ since the vectors within 
square brackets in (\ref{wiu}), one tangent and one normal to $\,\x$, 
are nonzero whenever $\,u\,$ is.
\end{proof}

\section{Proofs of Theorem~\ref{thsed} and Corollaries~\ref{codss} -- 
\ref{codot}}\label{pm}
\setcounter{equation}{0}
We begin by proving Theorem~\ref{thsed} in the special case where the $\,\mk\,$ 
distinct blown-up points, rather than being arbitrary, are selected in a 
particular way: the first $\,c\,$ of them lie in $\,\x$, and the last 
$\,\mk-c\,$ in $\,\y\smallsetminus\x$.

By blowing up the $\,\mk\,$ points we transform $\,\x\,$ into a complex 
submanifold $\,\x\hs'$ of the resulting \csu\ $\,\y'$ with 
$\,[\x\hs'\hs]=([\x\hs],1,\dots,1,0,\dots,0)$, where $\,1\,$ occurs $\,c\,$ 
times; this is clear since the $\,\bbZ^\mk$ component of $\,[\x\hs\hs]\,$ in the 
decomposition (\ref{htw}) is formed by intersection numbers of $\,\x\hs'$ with 
the ordered $\,\mk$-tuple of exceptional divisors. Hence, by (\ref{con}), the 
restriction of $\,c_{\hs1}(\y')\,$ to $\,\x\hs'$ is zero, so that, in view of 
Theorem~\ref{thqnd}, the inclusion mapping $\,\x\hs'\to\y'$ is homotopic to 
a \tre\ $\,\hat h:\x\to\y'\nh$. (We identify $\,\x\hs'$ with $\,\x\,$ 
using the blow-down projection $\,\pi:\y'\nh\to\y$.) If $\,\hat h$ is chosen 
sufficiently $\,C^1$-close to the inclusion $\,\x\hs'\to\y'$ (cf.\ 
Theorem~\ref{thqnd}), $\,\hat h(\x)\,$ will have a single, transverse 
intersection with each of the first $\,c$ exceptional divisors, and will not 
intersect the other $\,\mk-c\,$ of them. The composite 
$\,h=\pi\circ \hat h:\x\to\y\,$ thus is an embedding homotopic to the 
inclusion $\,\x\to\y$, having exactly $\,c\,$ complex points removable by 
blow-up (cf.\ Section~\ref{rb}), located at the first $\,c\,$ original blown-up 
points (that lay on $\,\x$), and no other complex points.

Let $\,y_1,\dots,y_\mk$ now be an arbitrary $\,\mk$-tuple of distinct points in 
$\,\y$. Deforming $\,h\,$ slightly if necessary, we may assume that 
$\,y_j\notin h(\x)\,$ for $\,j=1,\dots,\mk$. (If some $\,y_j$ is one of the 
$\,c\,$ complex points of $\,h(\x)\,$ removable by blow-up, we deform $\,h\,$ 
around $\,y_j$ using the flow of a $\,C^\infty$ vector field on $\,\y\,$ 
supported in a small neigborhood of $\,y_j$, which is hol\-o\-mor\-phic near 
$\,y_j$, so that removability by blow-up is preserved.) Applying the final 
clause of Lemma~\ref{tntcl} to these $\,h\,$ and $\,y_j$, we obtain an 
embedding $\,h':\x\to\y\,$ which, when the points $\,y_1,\dots,y_\mk$ are blown 
up, becomes the required \tre\ $\,f:\x\to\y'\nnh$, completing the proof of 
Theorem~\ref{thsed}.
\begin{proof}[Proof of Corollary~\ref{codss}]This is just Theorem~\ref{thsed} 
with $\,c=3d\,$ and $\,\x$ realized as a nonsingular degree $\,d\,$ curve in 
$\,\cp^2\nnh$. (For instance, the curve given by $\,x^d\nh+y^d\nh+z^d=\hh0\,$ 
in the homogeneous coordinates $\,[x,y,z]$.)
\end{proof}
\begin{proof}[Proof of Corollary~\ref{codot}]If $\,d\le2$, the assertion 
follows from Corollary~\ref{codss} for $\,d\le2\,$ (with a rearrangement of 
the blown-up points if $\,d=1$). Now let $\,d\ge3$, and let $\,\y\,$ be the 
\csu\ obtained by blowing up a point in $\,\cp^2\nnh$. Using a complex 
automorphism of $\,\cp^2\nnh$, we may assume that the blown-up 
point is $\,[\hs0,0,1]$. Equation $\,x^d=y^{d-1}z$, in the homogeneous 
coordinates $\,[x,y,z]$, defines a degree $\,d\,$ singular curve in 
$\,\cp^2\nnh$, and blowing up its unique singularity, at $\,[\hs0,0,1]$, we 
transform it into a nonsingular hol\-o\-mor\-phic curve $\,\x\subset\y$, 
\feic\ to $\,S^2\nh$, with $\,[\x\hs]=(d\hh,\hs d-1)$. (In the 
hol\-o\-mor\-phic local coordinates $\,\xi,\eta\,$ for $\,\y\,$ making the 
blow-down projection appear as $\,(\xi,\eta)\mapsto[\xi,\xi\eta,1]$, equation 
$\,x^d=y^{d-1}z\,$ reads $\,\xi=\eta^{d-1}\nh$, while the exceptional divisor 
is given by $\,\xi=0$.) We can now apply Theorem~\ref{thsed} to these 
$\,\y,\x,$ with $\,c=2d+1\,$ (cf.\ (\ref{con})) and $\,\mk=j-1$.
\end{proof}

\section{Further special cases of Theorem~\ref{maiem}}\label{rc}
\setcounter{equation}{0}
In Section~\ref{sc} we proved Theorem~\ref{maiem} except for two cases, with 
which we deal here and in the next section. In the first remaining case, 
$\,\x\,$ is assumed oriented and $\,\y\,$ is obtained from $\,\cp^2$ by 
blowing up any set of $\,\mk\,$ distinct points, $\,1\le\mk\le8$. As shown in 
Lemmas~\ref{solns} and~\ref{oimin}, $\,\x\,$ then must be either the torus 
$\,T^2\nnh$, with the degree $\,(0,\text{\bf0})$, or the sphere $\,S^2\nnh$, 
with one of the degrees (\ref{oii}) and their trivial modifications. (Cf.\ the 
beginnining of Section~\ref{de}.) We will now describe how each of these 
possibilities is realized.

For $\,(0,\text{\bf0})$, we fix a \tr\ $\,2$-torus in an open ball 
$\,\,U\subset\bbC^2$ (cf.\ (v) in Section~\ref{se}), and then use a 
hol\-o\-mor\-phic embedding $\,\,U\to\y$. That the degrees (\ref{oii}) and 
their trivial modifications are realized by \tr\ embedded $\,2$-spheres with 
$\,\mk\le8\,$ is in turn clear, respectively, from Example~\ref{round} 
(combined with Lemma~\ref{blwup}(a)), Corollary~\ref{untre} for $\,\mk=3\,$ 
(along with Lemmas~\ref{blwup}(a) and~\ref{oimin}), and Corollary~\ref{codot} 
for $\,d=2,3$.

\section{Theorem~\ref{maiem} for nonorientable surfaces}\label{tn}
\setcounter{equation}{0}
Given a \sc\ \acsu\ $\,\y\nh$, a \compact\ \rsu\ $\,\x$, and a fixed element 
$\,\dg\,$ of $\,H_2(\y,\bbZ_{[2]})\,$ (notation of (\ref{ztw})), consider the 
following condition, imposed on $\,\dg$, in which $\,q\,$ is defined by 
(\ref{pzq}), and $\,\mathcal{Z}(\x,\y)\,$ is the set described in the lines 
following Theorem~\ref{thuno}:
\begin{equation}\label{eid}
\begin{array}{l}
\text{\rm Every\ }\,(\ind,\dg)\in\mathcal{Z}(\x,\y)\,\text{\rm\ with\ the\ 
second\ component\ }\,\hs\dg\hs\,\text{\rm\ is\ the}\\
\text{\rm index-degree\ pair\ of\ some\ totally\ real\ embedding\ }\,f:\x\to\y.
\end{array}\hskip-10pt
\end{equation}
\begin{example}\label{eidex}Condition (\ref{eid}) is satisfied by
\begin{enumerate}
  \def\theenumi{{\rm\roman{enumi}}}
\item the torus $\,\x=T^2$ and $\,\y\,$ as above, with $\,\hs\dg=0$,
\item the \kb\ $\,\x=K^2$ and $\,\hs\dg=0$, with any $\,\y\,$ as above for 
which $\,q=\infty\,$ or $\,q\,$ is finite but not divisible by $\,4$,
\item $\y=\cp^1\nnh\times\cp^1$ and $\,\x=K^2\nnh$, with any 
$\,\hs\dg\in H_2(\y,\bbZ_2)\smallsetminus\{(1,1)\}$,
\item $\y=\cp^1\nnh\times\cp^1$ and $\,\x=S^2\nnh$, for  
$\,\hs\dg\hs=(\pm1,\mp1)\in\bbZ^2=H_2(\y,\bbZ)$,
\item $\y=\cp^2$ and $\,\x=\rp^2\nnh$, with 
$\,\hs\dg=1\in\{0,1\}=\bbZ_2=H_2(\y,\bbZ_2)$.
\end{enumerate}
In fact, (i) and (ii) are immediate from (a) -- (b) in Corollary~\ref{conas}, 
(iii) from Proposition~\ref{zkbsp}, and (iv) from Remark~\ref{dadsp} (where 
$\,\hs\ind\hs\,$ must equal $\,0\,$ by (\ref{isn})). Finally, (v) 
follows since $\,\mathcal{Z}(\x,\y)=\fri_6(\rp^2)\times\dzerm(\cp^2)$ has just 
one element (see (e) in Section~\ref{tq}, (\ref{irp}) and (\ref{dcp})), and 
that element is realized by the \tre\ described in (vi) of Section~\ref{se}.
\end{example}
In the following lemma, by the $\,\hs\text{\rm mod}\hskip4pt2\,$ reduction 
of $\,\hs\dg\in H_2(\y,\bbZ_2)\,$ we mean $\,\hs\dg\hs\,$ itself.
\begin{lem}\label{nwemb}Let\/ $\,\y\,$ be a fixed \sc\ \acsu. If condition 
{\rm(\ref{eid})} holds for a \compact\ \rsu\ $\,\x\,$ and\/ 
$\,\dg\in H_2(\y,\bbZ_{[2]})$, then\/ {\rm(\ref{eid})} will remain satisfied 
after $\,\x\,$ and\/ $\,\hs\dg\hs\,$ have been replaced by the connected sum 
$\,\x\,\#\,T^2\#\,K^2$ and the $\,\hs\text{\rm mod}\hskip4pt2\,$ reduction 
of\/ $\,\hs\dg$.
\end{lem}
\begin{proof}Let us assume (\ref{eid}), with our fixed $\,\y$, for some given 
$\,\x\,$ and $\,\dg$. First, $\,\y,\,\x\,\#\,T^2$ and $\,\dg\,$ then satisfy 
a modified version of (\ref{eid}), obtained when the word `embedding' is 
replaced by the phrase {\it immersion having just one double point, at which 
the self-in\-ter\-sec\-tion is transverse and, for orientable $\,\x$, also 
negative in the sense of Section~{\rm\ref{rt}}}. This is clear from 
Example~\ref{eidex}(i) applied, instead of $\,\y$, to a \sc\ open submanifold 
$\,\,U\,$ of $\,\y\smallsetminus\x$, combined with Theorem~\ref{cnsum} and 
(\ref{iff}). Specifically, to represent any given 
$\,(\ind,\dg)\in\mathcal{Z}(\x\,\#\,T^2,\y)\,$ with the second component 
$\,\hs\dg\hs\,$ by a \tri\ $\,f:\x\,\#\,T^2\to\y\,$ having a single 
self-in\-ter\-sec\-tion of the type just described, we write 
$\,\ind\hs=(\ind_{\nh*},\ind'\hh)$, as in (\ref{hss}.iii), and then obtain 
$\,f\,$ from the con\-nect\-ed-sum operation performed on the \tre s 
$\,\x\to\y\,$ and $\,T^2\nh\to U\subset\y\,$ that realize the index-degree 
pairs $\,(\ind_{\nh*},\dg)\,$ and $\,(\ind'\hh,0)$. (By (\ref{hss}), the 
additional condition imposed on elements of $\,\mathcal{Z}(\x,\y)\,$ or 
$\,\mathcal{Z}(\x\,\#\,T^2,\y)\,$ when $\,\x\,$ is nonorientable and 
$\,\chi(\x)\,$ is even, described in the lines following Theorem~\ref{thuno}, 
holds for $\,(\ind,\dg)\,$ if and only if it does for $\,(\ind_{\nh*},\dg)$.)

We now use Lemma~\ref{ctout} to remove the self-in\-ter\-sec\-tion of each 
\tri\ $\,f:\x\,\#\,T^2\to\y\,$ obtained as above, which, as explained in 
Section~\ref{rt}, gives rise to a \tre\ $\,\x\,\#\,T^2\#\,K^2\to\y$. 
Finally, the assertion about homotopy classes in Lemma~\ref{ctout} implies 
(\ref{eid}) for $\,\x\,\#\,T^2\#\,K^2$ and $\,\dg$. In fact, the 
con\-nect\-ed-sum operation in Lemma~\ref{ctout}, resulting in 
$\,\x\,\#\,T^2\#\,K^2\nnh$, may be viewed as involving one orientable summand: 
when $\,\x\,$ is orientable, the summand in question is $\,\x\,\#\,T^2\nnh$, 
while for nonorientable $\,\x\,$ it is $\ T^2$ (since $\,\x\,\#\,T^2\#\,K^2$ 
is then diffeomorphic to $\,\x\,\#\,T^2\#\,T^2$). Thus, if $\,\x\,$ is 
nonorientable, our claim follows from (\ref{hss}), since diffeomorphisms 
$\,T^2\to T^2$ equal to the identity on some disk in $\,T^2$ (and, therefore, 
admitting extensions to $\,\x\,\#\,T^2\#\,T^2$), act on 
$\,H_2(T^2\nnh,\bbZ)=\bbZ^2$ so as to realize every automorphism in 
$\,\text{\rm SL}(2,\bbZ)$. Consequently, the $\,T^2$ Maslov indices 
$\,\bbZ^2\to\bbZ_q$ constructed in Lemma~\ref{ctout} (which are trivial on one 
$\,\bbZ\,$ summand, and arbitrary on the other), represent, after 
re-pa\-ram\-e\-tri\-za\-tion, all Maslov indices $\,\bbZ^2\to\bbZ_q$, as 
required. Similarly, for orientable $\,\x$, we obtain our conclusion using 
diffeomorphisms of $\,K^2$ equal to the identity on a disk.
\end{proof}
We will now prove the last remaining case of Theorem~\ref{maiem}, in which 
$\hs\x\hs$ is assumed nonorientable. First, let $\,\y\,$ be 
$\,\bbC^2\nnh,\hs\cp^2$ or $\,\cp^1\nnh\times\cp^1\nnh$. Every nonorientable 
\compact\ surface admitting a \tre\ in $\,\y\,$ now can be obtained from just 
one or two low-genus {\it primary\/} surfaces by repeated applications of the 
operation $\,\x\,\mapsto\,\x\,\#\,T^2\#\,K^2$ (see Corollaries~\ref{codod}, 
\ref{cotrd} and~\ref{codie}); specifically, for 
$\,\y\nnh=\bbC^2\nnh$, or $\,\y\nh=\cp^2\nnh$, or 
$\,\y\nh=\cp^1\nnh\times\cp^1\nnh$, the primary surfaces are: $\,K^2$ alone, 
or $\,K^2$ and $\,\rp^2$ or, respectively, $\,K^2$ and $\,S^2$. The assertion 
of Theorem~\ref{maiem} for the \csu s $\,\bbC^2\nnh,\,\cp^2$ and 
$\,\cp^1\nnh\times\cp^1\nnh$, with all nonorientable $\,\x$, is now immediate 
from Lemma~\ref{nwemb}, since, in view of Example~\ref{eidex}, it holds for 
the primary surfaces just listed. (If $\,\y=\bbC^2\hs$ or $\,\y=\cp^2\nnh$, 
(e) in Section~\ref{tq} gives $\,q=\infty\,$ or $\,q=6$, and, in both cases, 
$\,\hs\dg=0\,$ for $\,\x=K^2$ by (\ref{ddh}.ii). Similarly, if 
$\,\y\nh=\cp^1\nnh\times\cp^1\nnh$, (\ref{ddh}.ii) and (\ref{abc}) imply that 
$\,\hs\dg\ne(1,1)\,$ for $\,\x=K^2\nnh$.)

Let $\,\y\,$ now be obtained from $\,\cp^2$ by blowing up any set of 
$\,\mk\ge1$ distinct points, and let $\,\x\,$ be a \tr\ \compact\ surface 
embedded in $\,\y$, By the {\it \cml} $\,\hs\text{\rm mod}\hskip4pt2\,$ {\it 
degree\/} of $\,\x\,$ we mean the pair $\,(d\hh,\es)\,$ with 
$\,d\in\{0,1\}$ and $\,\es\in\{0,1,\dots,\mk\}$, characterized as follows. 
The $\,\hs\text{\rm mod}\hskip4pt2\,$ reduction of the degree 
$\,[\x\hs]\in H_2(\y,\bbZ_{[2]})\,$ (equal to $\,[\x\hs]\,$ itself if $\,\x\,$ 
in not orientable) can, in analogy with (\ref{htk}), be treated as a 
$\,(\mk+1)$-tuple $\,(d\hh,q_1,\dots,q_\mk)\,$ of elements of 
$\,\bbZ_2=\{0,1\}$. This gives the value of $\,d\hh$, while $\,\es\,$ is 
defined to be the number of times that $\,1\in\bbZ_2$ occurs among the 
$\,q_j$. Thus, $\,\es\,$ is the sum of the $\,q_j$ treated (and added) as 
integers with $\,q_j\in\{0,1\}$.

For $\,\y,\mk,\x,d\hh,\es\,$ with the properties just listed,
\begin{equation}\label{bdm}
\begin{array}{l}
\text{\rm blowing\ up\ a\ point\ of\ }\,\hs\y\,\text{\rm\ that\ lies\ in\ }
\,\x\hs\,\text{\rm\ replaces\ the}\\
\text{\rm quadruple\ \ }\mk,\x,d,\es\text{\rm\ \ with\ \ }\,\mk+1,
\x\,\#\,\rp^2\nnh,d,\es+1\hs.
\end{array}\hskip16pt
\end{equation}
This is immediate from Lemma~\ref{blwup}(b) for $\,n=\ek=2\,$ and
Lemma~\ref{blwup}(c).

If, in addition, $\,\x\,$ is nonorientable or \feic\ to $\,S^2\nnh$, 
then
\begin{equation}\label{dvf}
d\,-\,\es\,-\,\chi(\x)\,\,\text{\rm\ is\ divisible\ by\ }\,4\hs,
\end{equation}
by (\ref{ddh}) with $\,\hs\dg^{\hskip.2pt2}=d-q_1-\ldots-q_\mk$, for 
$\,q_j\in\{0,1\}\subset\bbZ\,$ as above; the $\,S^2$ case follows if one 
applies the $\,\hs\text{\rm mod}\hskip4pt4\,$ reduction to (\ref{ddh}.i), 
since $\,2=-\hs2\,$ in $\,\bbZ_4$. We can rephrase (\ref{dvf}) in the form of 
a table:

\begin{table}[ht]
\vskip-8pt
\caption{The values, allowed by (\ref{dvf}), of the \cml\  
$\,\hs\text{\rm mod}\hskip4pt2\,$ degree $\,(d\hh,\es)\,$ for 
\tre s with $\,1\le\es\le7\,$ and the four `primary' surfaces listed in the 
top row. Some values of $\,(d\hh,\es)\,$ for $\,S^2$ and for $\,\es>7\,$ are 
listed as well.}
\label{totaldeg}
\vskip-15pt
\hskip6pt
\vbox{\offinterlineskip
\hrule
\halign
{&\vrule#&\strut
\hskip1.46pt\hfil#\hfil\hskip1.46pt\cr
&\hskip10pt\vbox{\vskip-.8pt\hbox{$_{\text{\rm Row}}$}\vskip-1pt}\hskip10pt&
&\vbox{\vskip-1pt\hbox{$_{S^2}$}\vskip-1.1pt}&
&\vbox{\vskip-1pt\hbox{$_{\rp^2}$}\vskip-1.1pt}&
&\vbox{\vskip-1pt\hbox{$_{K^2}$}\vskip-1.1pt}&
&\vbox{\vskip-1pt\hbox{$_{K^2\#\,\rp^2}$}\vskip-2.3pt}&
&\vbox{\vskip-1pt\hbox{$_{K^2\#\hs K^2}$}\vskip-2.3pt}&\cr
\noalign{\hrule}
&{\sl1}&&$\,$&&$\,$&&$(0\hh,0)$&&$(0\hh,1)$&&$(0\hh,2)$&\cr
\noalign{\hrule}
&{\sl2}&&\hskip18pt$(0\hh,2)$\hskip18pt&&\hskip18pt$(0\hh,3)$\hskip18pt&
&\hskip18pt$(0\hh,4)$\hskip18pt&&\hskip18pt$(0\hh,5)$\hskip18pt&
&\hskip18pt$(0\hh,6)$\hskip18pt&\cr
\noalign{\hrule}
&{\sl3}&&$(0\hh,6)$&&$(0\hh,7)$&&$(0\hh,8)$&&$(0\hh,9)$&&$(0\hh,10)$&\cr
\noalign{\hrule}
&{\sl4}&&$(0\hh,10)$&&$(0\hh,11)$&&$\,$&&$\,$&&$\,$&\cr
\noalign{\hrule}
&{\sl5}&&$\,$&&$(1\hh,0)$&&$(1\hh,1)$&&$(1\hh,2)$&&$(1\hh,3)$&\cr
\noalign{\hrule}
&{\sl6}&&$(1\hh,3)$&&$(1\hh,4)$&&$(1\hh,5)$&&$(1\hh,6)$&&$(1\hh,7)$&\cr
\noalign{\hrule}
}}   
\vskip-5pt
\end{table}

Proving Theorem~\ref{maiem} for nonorientable surfaces is now reduced to 
realizing, for each row of Table~2, the leftmost \cml\ 
$\,\hs\text{\rm mod}\hskip4pt2\,$ degree $\,(d\hh,\es)\,$ appearing in that 
row by a \tre\ $\,\x\to\y\,$ (where $\,\y$, again, arises from $\,\cp^2$ by 
blowing up $\,\mk\,$ points, and, this time, $\,\mk=\es$). In fact, for the 
primary surfaces, the remaining \cml\ $\,\hs\text{\rm mod}\hskip4pt2\,$ 
degrees, in each row separately, are then immediately realized by starting 
from the leftmost one, and then successively applying (\ref{bdm}). Next, for 
all the remaining (non-pri\-ma\-ry) nonorientable \compact\ surfaces, 
Theorem~\ref{maiem} follows from Lemma~\ref{nwemb} (which now states that 
realizability for $\,\x\,$ implies the same for $\,\x\,\#\,T^2\#\,K^2\nnh$, 
with $\,\mk\,$ and $\,(d\hh,\es)\,$ unchanged). On the other hand, 
Lemma~\ref{blwup}(a) allows us to replace $\,\mk=\es\,$ with any $\,\mk\ge\es$.

The leftmost \cml\ $\,\hs\text{\rm mod}\hskip4pt2\,$ degree $\,(d\hh,\es)\,$ 
in each row is in turn realized as follows. Row {\sl 1}, by 
Example~\ref{trtkb} and Lemma~\ref{indeg}; Row {\sl 2}, in view of 
Example~\ref{round}; Row {\sl 3}, from Corollary~\ref{codot} for $\,d=2\,$ 
(and $\,j=6$); Row {\sl 4}, by Corollary~\ref{codot} for $\,d=4\,$ (and 
$\,j=10$); Row {\sl 5}, by  (vi) of Section~\ref{se}; and, finally, Row {\sl 
6}, from Corollary~\ref{codss} for $\,d=1\,$ (and $\,\mk=3$).


The $\,\mk\,$ distinct blown-up points in the above argument can be made 
arbitrary, even though the left-to-right steps in every row of Table~2 use 
Lemma~\ref{blwup}, with the blown-up point lying on the original \tr\ surface. 
Namely, by Lemma~\ref{blwup}(c), after a small deformation each of the 
resulting \tr\ surfaces $\,\x\,$ will have just one, transverse intersection 
point with any exceptional divisor that it intersects. Blowing down all such 
divisors gives rise to a finite set of complex points, all removable by 
blow-up, on a new surface $\,\x\hs'\nnh$, and the final clause of 
Lemma~\ref{tntcl} allows us to deform $\,\x\hs'$ so as to move the complex 
points to arbitrary prescribed locations, where they can be blown up again.

\setcounter{section}{1}
\renewcommand{\thesection}{\Alph{section}}
\section*{Appendix. Spheres and tori in 
$\,\cp^2\,\#\,9\,\overline{\cp^2}$}\label{st}
\setcounter{equation}{0}
Our proof, in the last section, of Theorem~\ref{maiem} for {\it 
nonorientable\/} surfaces, works not only if $\,1\le\mk\le8$, but also for 
$\,\mk=9,10,11$. The question of extending the classification of 
Theorem~\ref{maiem} to $\,\mk\in\{9,10,11\}\,$ is thus reduced to determining 
which values of the degree are realized, for such $\,\mk$, by \tre s of 
{\it oriented\/} \compact\ surfaces.

We give here a partial answer to this question for $\,\mk=9$. Namely, let 
$\,\y=\,\cp^2\,\#\,9\,\overline{\cp^2}$ be the \csu\ obtained from $\,\cp^2$ 
by blowing up any ordered set of nine distinct points. By Lemma~\ref{solns}, 
the only \compact\ oriented surfaces admitting \tre s in $\,\y\,$ are $\,T^2$ 
and $\,S^2\nnh$. In the case of $\,T^2\nnh$, Lemma~\ref{solns} also provides a 
list of possible degrees; in Example~\ref{knine} below we realize all those 
degrees by \tre s. Theorem~\ref{ppsss} below provides, in turn, an ``algebraic 
classification'' of the degrees of \tre s $\,S^2\to\y$, analogous to that for 
$\,T^2$ in Lemma~\ref{solns}. We begin with some general remarks.

Since being totally real is an open property, a \tri\ $\,f:\x\to\y\,$ with 
compact $\,\x\,$ and $\,\dimr\x<\dimc\y\,$ leads, via a generic small 
deformation of $\,f$, to totally real {\it embeddings}. The following lemma 
shows that a similar deformation sometimes works when $\,\dimr\x=\dimc\y$.

Given a \compact\ \rmf\ $\,\x$, an integer $\,d\ge2$, and a surjective 
homomorphism $\,\varphi:\pi_1\x\to\bbZ_d$, let $\,\pi:\tilde\x\to\x\,$ be the 
$\,d$-fold covering projection corresponding to the subgroup 
$\,\hs\text{\rm Ker}\,\varphi\,$ of $\,\pi_1\x$. Thus, $\,\tilde\x\,$ is a 
principal $\,\bbZ_d$ bundle over $\,\x$, and we will denote by $\,\lambda\,$ 
the complex line bundle associated with it via the standard action of 
$\,\bbZ_d$ on $\,\bbC\hs$. Note that $\,\lambda^{\otimes d}$ is trivial. 
\begin{lem}\label{trprd}With\/ $\,\x,d,\varphi,\pi,\tilde\x\,$ and\/ 
$\,\lambda\,$ as above, let us suppose that\/ $\,\lambda$ is real-isomorphic 
to a vector subbundle of\/ $\,T\x$. Then, for any totally real embedding\/ 
$\,f:\x\to\y\,$ in an \acm\/ $\,\y$, there exists a totally real embedding\/ 
$\,\tilde\x\to\y\,$ which is $\,C^\infty$ homotopic, through totally real 
immersions, to the composite immersion\/ $\,f\circ\pi$.
\end{lem}
\begin{proof}Since $\,f\,$ is \tr, our assumption about $\,\lambda\,$ allows 
us to choose an embedding $\,F:U\to\y$, where $\,\,U\,$ is a \nb\ of 
the zero section $\,\x\,$ in the total space of $\,\lambda$, such that 
$\,F=f\,$ on $\,\x$. (In fact, we may let $\,F\,$ be the composite of an 
injective real vec\-tor-bun\-dle morphism from $\,\lambda\,$ to the normal 
bundle of $\,f$, which exists in view of (\ref{tnu}), followed by a 
tu\-bu\-lar-neigh\-bor\-hood \feo.) As $\,\tilde\x\,$ is naturally embedded in 
the total space of the unit circle bundle of $\,\lambda\,$ and 
$\,\pi:\tilde\x\to\x\,$ is the restriction to $\,\tilde\x\,$ of the bundle 
projection $\,\pi:\lambda\to\x$, the mappings given by 
$\,\tilde\x\ni\xi\mapsto F(t\hh\xi)\in\y$, each of them depending on a fixed 
parameter $\,t\ge0\,$ close to $\,0$, form a $\,C^\infty$ homotopy between 
$\,f\circ\pi\,$ (with $\,t=0$) and an embedding $\,\tilde\x\to\y\,$ (with any 
$\,t>0\,$ near $\,0$). Being \tr\ is an open property; thus, for $\,t\,$ close 
to $\,0\,$ such embeddings are \tr.
\end{proof}
\begin{example}\label{ttinm}Given a \tre\ $\,f:T^2\to\y\,$ of the $\,2$-torus 
in an \acsu\ $\,\y\,$ and a $\,d$-fold self-co\-ver\-ing 
$\,\pi:T^2\to T^2\nh$, $\,d\ge2$, there exists a \tre\ $\,T^2\to\y\,$ 
homotopic through \tri s to $\,f\nh\circ\pi$. This is clear from 
Lemma~\ref{trprd} for $\,\x=\tilde\x=T^2\nh$. (The corresponding $\,\lambda\,$ 
is trivial, since so is $\,\lambda^{\otimes d}$ and $\,\pi_1\x\,$ is free.)
\end{example}
\begin{example}\label{knine}Let 
$\,\y\hs=\,\cp^2\,\#\,\,9\,\overline{\cp^2}$ be the \csu\ obtained by blowing 
up any set of nine distinct points in $\,\cp^2\nnh$. By Lemma~\ref{solns}, the 
degree of any \tre\ $\,T^2\to\y\,$ equals 
$\,(3s\hh,\hs s,s,s,s,s,s,s,s,s)$ for some $\,s\in\bbZ$. We now show 
that, conversely, every such degree is realized by some \tre\ $\,f:T^2\to\y$. 
First, if $\,s=1$, the embedding $\,f\,$ is provided by  Corollary~\ref{codss} 
with $\,d=3\,$ and $\,\mk=9$. If $\,s>1$, we can use Example~\ref{ttinm} for 
this last $\,f\,$ and $\,d=s$. Finally, in the cases $\,s<0\,$ and $\,s=0\,$ 
it suffices to re-orient $\,\x\,$ or, respectively, invoke 
Corollary~\ref{conas}(a).
\end{example}
For the degrees of \tr\ embedded $\,2$-spheres, we have a partial result:
\begin{thm}\label{ppsss}For every prescribed integer\/ $\,d$, the system\/ 
{\rm(\ref{pot})} with $\,\mk=9$ and\/ $\,\chi=2\,$ has a unique solution\/ 
$\,(d,\text{\smallbf q})=(d\hh,q_1,\dots,q_\mk)\,$ that satisfies the 
normalizing condition\/ $\,q_1\ge\ldots\ge q_\mk$, cf.\ {\rm 
Remark~\ref{staso}}. Explicitly, we have
\[
(q_1,\dots,q_9)= 
\begin{cases}
{(s+1,s,s,s,s,s,s,s,s-1)\,,}&\text{\rm if \ }$$d=3s,\,s\in\bbZ\,,\cr
{(s+1,s+1,s+1,s,s,s,s,s,s)\,,}&\text{\rm if \ }$$d=3s+1,\,s\in\bbZ\,,\cr
{(s,s,s,s,s,s,s-1,s-1,s-1)\,,}&\text{\rm if \ }$$d=3s-1,\,s\in\bbZ\,.\cr
\end{cases}
\]
\end{thm}
In fact, (\ref{psk}) with $\,\mk=9\,$ and $\,\chi=2\,$ yields 
$\,|\text{\smallbf q}-\text{\smallbf s}\hs|^2\in\{2,3\}$. Setting 
$\,p_j=q_j-s\,$ and then decomposing $\,2\,$ or $\,3\,$ into all possible sums 
$\,\sum_{j=1}^9p_j^2$ with $\,p_j\in\bbZ$, $\,p_1\ge\ldots\ge p_9$ and 
$\,\sum_{j=1}^9p_j=3(d-3s)\,$ (which is the first equation in (\ref{pot})), we 
easily obtain the required formula for $\,(q_1,\dots,q_9)$.
%

\end{document}